\documentclass[11pt]{amsart}
\pdfoutput=1
\usepackage{latexsym, graphicx, epsfig, amsmath, amsfonts,amssymb}
\usepackage{multirow}
\usepackage{color}
 \usepackage[percent]{overpic}
\usepackage{algorithm}
\usepackage{diagbox}
\usepackage{algorithmicx}
\usepackage{algpseudocode}
\floatname{algorithm}{Algorithm}
\usepackage{threeparttable}
\usepackage{booktabs}

\def\mb{\mathbf}
\def\R{\mathbb{R}}

\usepackage{epstopdf,mathtools,bbm}
\usepackage[left=1.25in,right=1.25in,top=1in]{geometry}

\def\mb{\mathbf}
\def\R{\mathbb{R}}

\newcounter{Rownumber}

\title{An adaptive surrogate modeling based on deep neural networks for large-scale Bayesian inverse problems}

\author{Liang Yan}
\thanks{Department of Mathematics, Southeast University, Nanjing, China. Email: yanliang@seu.edu.cn. L. Yan is supported by NSF of China (No.11771081), the science challenge project (No. TZ2018001),  Qing Lan project of Jiangsu Province and  Zhishan Young Scholar Program of SEU}

\author{Tao Zhou}
\thanks{LSEC, Institute of Computational Mathematics, Academy of Mathematics and Systems
Science, Chinese Academy of Sciences, Beijing 100190, China. Email: tzhou@lsec.cc.ac.cn. T. Zhou is partially supported the NSFC (under grant numbers 11822111, 11688101, and 11731006), the science challenge project (No. TZ2018001), NCMIS, and the youth innovation promotion association (CAS)}

\date{March 1, 2020}
\begin{document}

\maketitle

\begin{abstract}
In Bayesian inverse problems, surrogate models are often constructed to speed up the computational procedure, as the parameter-to-data map can be very expensive to evaluate. However, due to the curse of dimensionality and the nonlinear concentration of the posterior, traditional surrogate approaches (such us the polynomial-based surrogates) are still not feasible for large scale problems. To this end, we present in this work an adaptive multi-fidelity surrogate modeling framework based on deep neural networks (DNNs), motivated by the facts that the DNNs can potentially handle functions with limited regularity and are powerful tools for high dimensional approximations. More precisely, we first construct offline a DNNs-based surrogate according to the prior distribution, and then, this prior-based DNN-surrogate will be adaptively \& locally refined online using only a few high-fidelity simulations. In particular, in the refine procedure, we construct a new shallow neural network that view the previous constructed surrogate as an input variable -- yielding a composite multi-fidelity neural network approach. This makes the online computational procedure rather efficient. Numerical examples are presented to confirm that the proposed approach can obtain accurate posterior information with a limited number of forward simulations.
\end{abstract}


\pagestyle{myheadings}
\thispagestyle{plain}
\markboth{L. Yan L and T. Zhou}
{ ADNN for BIPs}

\section{Introduction}
Inverse problems arise when one is interested in determining  model parameters or inputs from a set of indirect observations \cite{Evans+Stark2002, Kaipio+Somersalo2005}.   Typically, inverse problems are  ill-posed in the sense that the solution may not exist or may not be unique. More importantly, the parameters may not depend continuously on the observations -- meaning that \textcolor{black}{a small perturbation in the data may cause an enormous deviation in the solution}. The Bayesian approach \cite{Kaipio+Somersalo2005,Stuart2010} is a popular approach for inverse problems \textcolor{black}{which}  casts the solution as a {\it posterior distribution} of the unknowns conditioned on observations, and introduces regularization in the form of \textit{prior} information. By estimating statistic moments according to the posterior distribution, one not only gets point estimates of the parameters, but also obtains a complete description of the uncertainty in model predictions.  However, in practice, the analytical treatment for \textcolor{black}{the} posterior is not feasible in general due to the complexity of the system. Consequently, the posterior is often approximated with numerical approaches such as the Markov chain Monte Carlo (MCMC) method.

In the standard MCMC approach, one aims at generating samples directly from the posterior distribution over the parameters space by using the unnormalized posterior, i.e., the product of the prior and likelihood.  However,  the cost of evaluating the likelihood in the sampling procedure can quickly become prohibitive if the forward model is computationally expensive.  One popular way to reduce the computational cost in the sampling procedure is to replace the original forward model with a cheap surrogate model
\cite{Galbally+Fidkowski+Willcox+Ghattas2010, Jin2008fast,  kennedy2001, Li+Lin2015JCP, Lieberman+Willcox+Ghattas2010, Marzouk+Najm2009, Marzouk+Najm+Rahn2007, Marzouk+Xiu2009, stuart+teckentrup2016, yan+guo2015}. Using a \textcolor{black}{computationally} less expensive, offline constructed, surrogate model can make the online computations very efficient. Furthermore, theoretical analysis shows that if the surrogate converges to the true model in the prior-weighted $L_2$ norm, then the posterior distribution generated by the surrogate converges to the true posterior \cite{Marzouk+Xiu2009,stuart+teckentrup2016,yan+guo2015,Yan+Zhang2017IP}.

Although the surrogate approach can provide significant empirical performance improvements, there are however many challenges for practical applications. First, constructing a sufficiently accurate surrogate over the whole domain of the prior distribution may not be possible for many practical problems. \textcolor{black}{Especially, the posterior distribution often concentrates on a small fraction of the support of the prior distribution}, and a globally prior-based surrogate may not be accurate for online computations \cite{Li+Marzouk2014SISC}. To improve this,  posterior-focused approaches have been suggested recently, where one constructs a sequence of local surrogates in the important region of the posterior distribution, to alleviate the effect of the concentration of posterior \cite{Conrad2016JASA,Cui2014data,Li+Marzouk2014SISC}.

In our previous work \cite{Yan+Zhou19JCP}, we also presented an adaptive multi-fidelity surrogate modelling procedure based on PCEs to speed up the online computations via MCMC. The idea is to begin with a low fidelity PCE-surrogate, and then correct it adaptively using online high fidelity data. Empirical studies on problems of moderate dimension show that the number of high-fidelity model evaluations can be reduced by orders of magnitude, with no discernible loss of accuracy in posterior expectations. Nevertheless, the approaches in \cite{Yan+Zhou19JCP} also admit some limitations: (i) the PCE surrogate has limitations to handle problems with low regularity; (ii) the PCE surrogate suffers from the so called curse of dimensionality. This motivate the present work: we shall present an adaptive multi-fidelity deep neural networks (DNNs) based surrogate modeling for large-scale BIPs, motivated by the facts that DNNs can potentially handle functions with limited regularity and are powerful tools for approximating high dimensional problems (\cite{Han+Jentzen+E2018PNAS,Raissi2019JCP,Schwab+Zech2019AA,Tripathy+Bilionis2018JCP,Zhu+Zabaras2018bayesian}). The key idea is to view the low fidelity surrogate as an input variable into the DNN surrogate of the next iteration -- yielding a composite DNNs that combine two surrogates between two iterations. Another key issue is to adaptively correct the DNN-surrogate locally so that the new surrogate is refined on a more concentrated (posterior) region in the parameter space. By doing this, one can perform the online correction procedure in a very efficient way. We shall present numerical experiments to show the effectiveness of the new approach. To the best of our knowledge, this is the first investigation of the multi-fidelity DNN-surrogate for Bayesian inverse problems.

The rest of the paper is organized as follows. In Section \ref{sec:setup}, we present some preliminaries and provide with a mathematical description of the BIPs. The adaptive multi-fidelity DNNs-surrogate approach is discussed in Section 3.  In Section 4, we use a benchmark elliptic PDE inverse problem to demonstrate the accuracy and efficiency of our approach. Finally, we give some concluding remarks in Section 5.

\section {Background and problem formulation}\label{sec:setup}
In this section, we shall review some basic ideas for the surrogate-based approach in Bayesian inverse problems.
\subsection{Bayesian inverse problems} \label{sec:BIPs}

We consider a discretized system of a mathematical model (such as PDEs) of interest:
\begin{equation}\label{stateq}
F(u_h( z), z) = 0,
\end{equation}
where $u_h(z):  \Xi \rightarrow  \R^{n_h}$ is the discrete solution with $n_h$ being the dimension of the finite-dimensional discretization in the physical domain, and $z \in \Xi \subset \R^{n}$ is an $n$-dimensional parameter vector. The discrete operator $F$ denotes a numerical approximation, e.g., by the finite element  or finite difference method for PDEs.  The goal of  an inverse problem is to estimate the unknown parameter vector $z$ from noisy observations \textcolor{black}{$d \in \R^m$} of the states $u_h$ given by
\begin{equation}\label{dataeq}
d= g(u_h; z)+\xi.
\end{equation}
Here $g$ is a discretized observation operator mapping from the states and parameters to the observable, and $\xi \in \R^m$ is the measurement error (or the noise). The system model (\ref{stateq}) together with the observation model (\ref{dataeq}) define a forward model $y=f(z)$ that maps the unknown parameter to the observable data.

To formulate the inverse problem in a Bayesian framework, we model the parameter $z$ as a random variable (vector), and endow it with a prior distribution $\pi(z)$.  The distribution of the $z$  conditioned on the data $d$, i.e., the posterior distribution $\pi(z|d)$ follows the Bayes' rule:
\begin{eqnarray*}\label{ppdf}
\pi(z|d) \propto \mathcal{L}(z| d,f) \pi(z).
\end{eqnarray*}
In case the density information of the $\xi \sim p_{\xi}$ is given, it follows directly that the likelihood  function can be written as:
\begin{eqnarray}\label{likefun}
\mathcal{L}(z| d,f) = p_{\xi}(d-f(z)).
\end{eqnarray}
Notice that each evaluation of the likelihood function $\mathcal{L}$ requires an evaluate of the forward model $f$. Thus, in sampling schemes such as the MCMC approach one has to perform $\sim 10^5$ forward model simulations, and this is a great challenge if the forward model $f$ represents a large-scale computer model. Consequently, it is popular approach  to construct (in a offline procedure) a cheaper surrogate for the forward model, and use it for online computations.

\subsection{Surrogate-based Bayesian inference}\label{sec:SBM}

Surrogate-based Bayesian inference has received much attention in recently years \cite{Frangos+Marzouk+Willcox2010}.  In this approach, one usually generate a collection of model evaluations (snapshots) $\mathcal{D}:=\{(z,f(z))\}$, and then construct an approximation $\tilde{f}$ based on those snapshots.  Using this approximation $\tilde{f}$, one can obtain an approximated surrogate posterior
\begin{eqnarray*}\label{ppdf_surrogate}
\widetilde{\pi}(z|d) \propto \mathcal{L}(z| d,\tilde{f}) \pi(z).
\end{eqnarray*}

Notice that if the evaluation of the approximation $\tilde{f}$ is inexpensive, then the approximate posterior density $\widetilde{\pi}$  can be evaluated  for a large number of samples, without resorting to additional simulations of the forward model $f$.  Although the surrogate-based Bayesian inference procedure can be quite effective, the big challenge is that in high-dimensional parametric spaces, the number of training points used to build the surrogate (such as the PCE approach) grows fast with respect to the dimension, and this is known as the {\it curse of dimensionality}. To improve this, we shall introduce in the next section a deep neural networks based surrogate modeling which can potentially handle high dimensional Beyesian inference problems.

\section{An adaptive multi-fidelity DNN-based surrogate modeling}\label{sec:method}

In this section, we shall present an adaptive multi-fidelity DNN-based surrogate modeling for Bayesian inverse problems.

\subsection{Feedforward DNN-based surrogate modeling}\label{sec:DNN}

\begin{figure}
\begin{center}
  \begin{overpic}[width=.65\textwidth,trim=20 0 20 15, clip=true,tics=10]{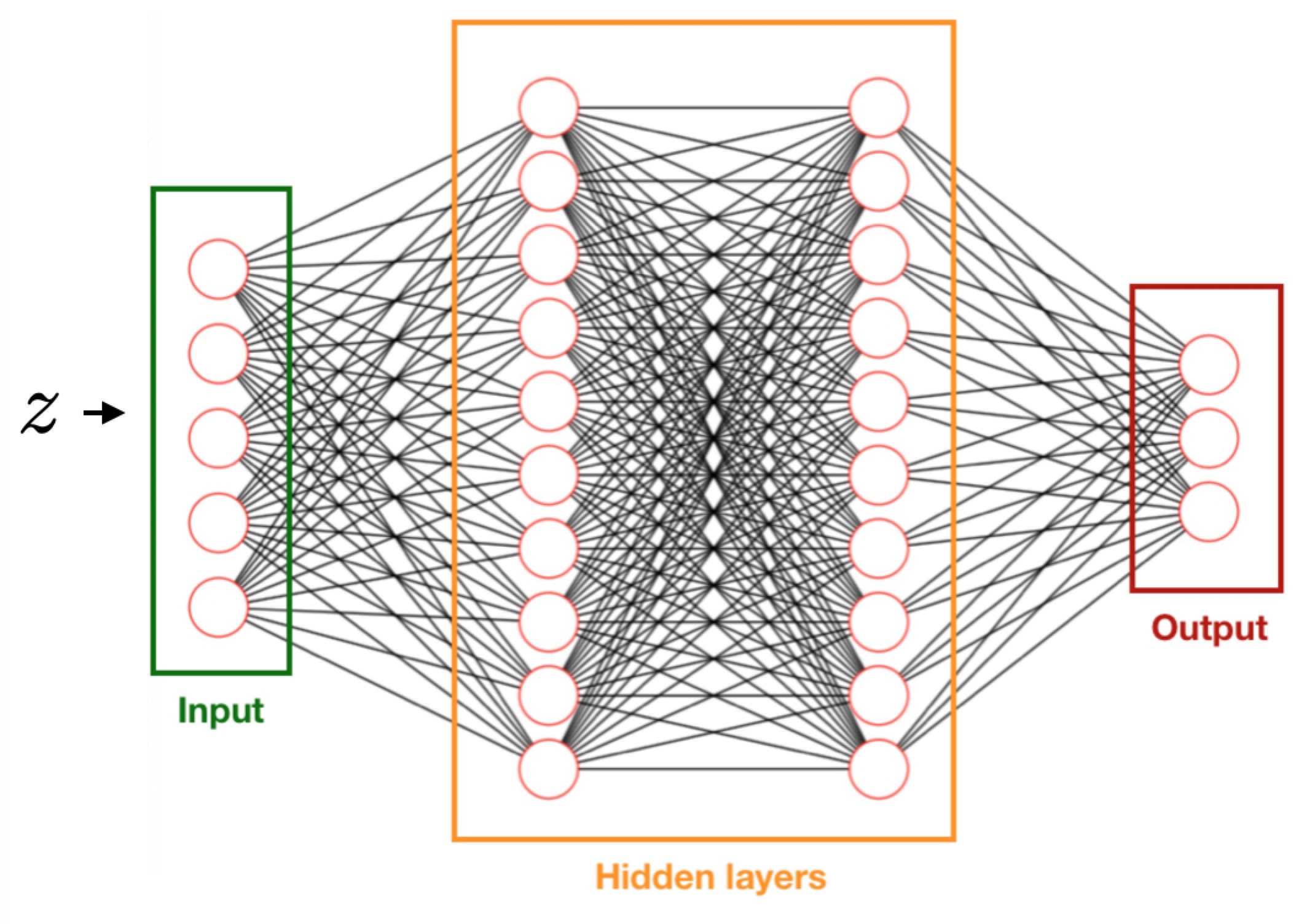}
  \end{overpic}
  \end{center}
\caption{The structure of a two-hidden-layer neural network. The first column of nodes (from left to right) is the input layer, taking an $n$-dimensional vector $z$ as input, and the last column is the output layer that generate an $m$-dimension output $\tilde{z}$. The intermediate columns are the hidden layers.}\label{DNN_structure}
\end{figure}

The basic idea of deep neural networks (DNNs) for surrogate modeling is that one can  approximate  an input-output map $f: \R^n\rightarrow \R^m$ through a hierarchical abstract layers of latent variables \cite{Goodfellow2016DL}. A typical example is the feedforward neural network, which is also called multi-layer perception (MLP). It consists of a collection of layers that include an input layer, an output layer, and a number of hidden layers. The size of the input layer and output layer are fixed and determined by the dimensionality of the input and output.  Figure \ref{DNN_structure} illustrates the structure of a DNN with two hidden layers.  Each circle in the schematic of the DNN is a neuron which calculates a weighted sum of an input vector plus bias and applies a non-linear function to produce an output. Specifically,  given an $n$-dimensional input row vector $z\in \R^n$, \textcolor{black}{we can define a DNN with $L$ hidden layers as following
\begin{align}
&\mathcal{NN}(z)=\mb{W}^{(L)}a^{(L)}+\mb{b}^{(L)},\\
&a^{(k+1)}=\sigma(\mb{W}^{(k)}a^{(k)}+\mb{b}^{(k)}), \quad  k=0,...,L-1.
\end{align}
Here $\mb{W}^{(k)} \in \R^{d_{k+1}\times d_{k}}, \, \mb{b}^{(k)} \in \R^{d_{k+1}}$ are the weights and biases of the network, $d_k$ is the number of neurons in the $k$th layer and $\sigma$ is  the activation function. Notice that here $a^{(0)}$ is the input $z$ and $d_0=n$.} Some popular choices for the activation function include sigmoid,  hyperbolic tangent,  rectified linear unit (ReLU),  to name a few \cite{Goodfellow2016DL,Ramachandran2017}. In the current work, we shall use Swish as the activation function \cite{Ramachandran2017,Tripathy+Bilionis2018JCP}:
\begin{equation*}
\sigma(z)=\frac{z}{1+\exp(-z)}.
\end{equation*}

Once the network architecture is defined, one can resort to optimization tools to find the unknown parameters $\theta =\{ \mb{W}^{(k)},  \mb{b}^{(k)}\}$ based on the training data. Precisely, let $\mathcal{D}: =\{(z_i, y_i)\}^N_{i=1}$ be a set of training data, we can define the following minimization problem:
\begin{equation}\label{thetastar}
\arg\min_{\theta} \frac{1}{N}\sum^N_{i=1} \|y_i-\mathcal{NN}(z_i;\theta)\|^2+\lambda \Omega(\theta),
\end{equation}
where $\mathcal{J}(\theta; \mathcal{D}) =  \frac{1}{N}\sum^N_{i=1} \|y_i-\mathcal{NN}(z_i;\theta)\|^2+\lambda \Omega(\theta)$ is the so called loss function, $\Omega(\theta)$ is a regularizer and $\lambda$ is  the regularization constant.  For our case, the regularizer is set to be \textcolor{black}{$\Omega(\theta) = \|\theta\|^2$}. In practice, the averaging in Eq.(\ref{thetastar}) is performed over a small randomly sampled subset $\mathcal{D}_{M}\subset \mathcal{D}$, at each iteration of the optimization procedure.  Solving this problem is generally \textcolor{black}{achieved by the stochastic gradient descent (SGD)} algorithm \cite{Bottou2010}.   SGD simply minimizes the function by taking a negative step along an estimate  of the gradient  $\nabla_{\theta} \mathcal{J}(\theta; \mathcal{D}_M)$ at iteration $k$. The gradients are usually computed through backpropagation.  \textcolor{black}{At each iteration, SGD updates the solution by
\begin{equation*}
\theta_{k+1}=\theta_k-\epsilon \nabla_{\theta} \mathcal{J}(\theta; \mathcal{D}_M),
\end{equation*}
where $\epsilon$ is the learning rate.} Recent algorithms that offers adaptive learning rates are available, such as Ada-Grad \cite{Zeiler2012Adadelta}, RMSProp \cite{Tieleman+Hinton2012lecture} and Adam \cite{Kingma2014Adam}, ect.  The present work adopts Adam optimization algorithm.

It is clear that after obtaining the parameters $\theta$, we have an explicit functional form $\mathcal{NN}(z;\theta)$. This approximation can be then substituted into the computation procedure of the surrogate posterior, such as in the MCMC framework.  However, we remark that a prior based surrogate might not be accurate enough for online computations, see. e.g. \cite{lu2015JCP,Yan+Zhou19JCP}. Thus, one usually needs to combine the surrogate with additional high fidelity data, yielding a multi-fidelity approach. To this end, we shall present in the next section an adaptive multi-fidelity DNN-surrogate modeling to accelerate the solution of BIPs.

\subsection{An adaptive multi-fidelity DNN surrogate}

In this section, we shall present an adaptive multi-fidelity DNN-based surrogate for BIPs. In particular, we consider to combine our approach within the MCMC framework, and extensions to general approaches  such as Ensemble Kalman inversion \cite{Iglesias+Law+Kody2013ensemble,Yan+Zhou19IJUQ} or RTO \cite{Bardsley2014SISC,Wang+Bardsley2017SISC} are also possible. Our approach is motivated by recent works such as \cite{Lu+Zhu2019MNN,Meng+Karniadakis2019MNN}, where composite DNNs are discussed to deal with multi-fidelity data.
 \begin{figure}
\begin{center}
  \begin{overpic}[width=.8\textwidth,trim=20 0 20 15, clip=true,tics=10]{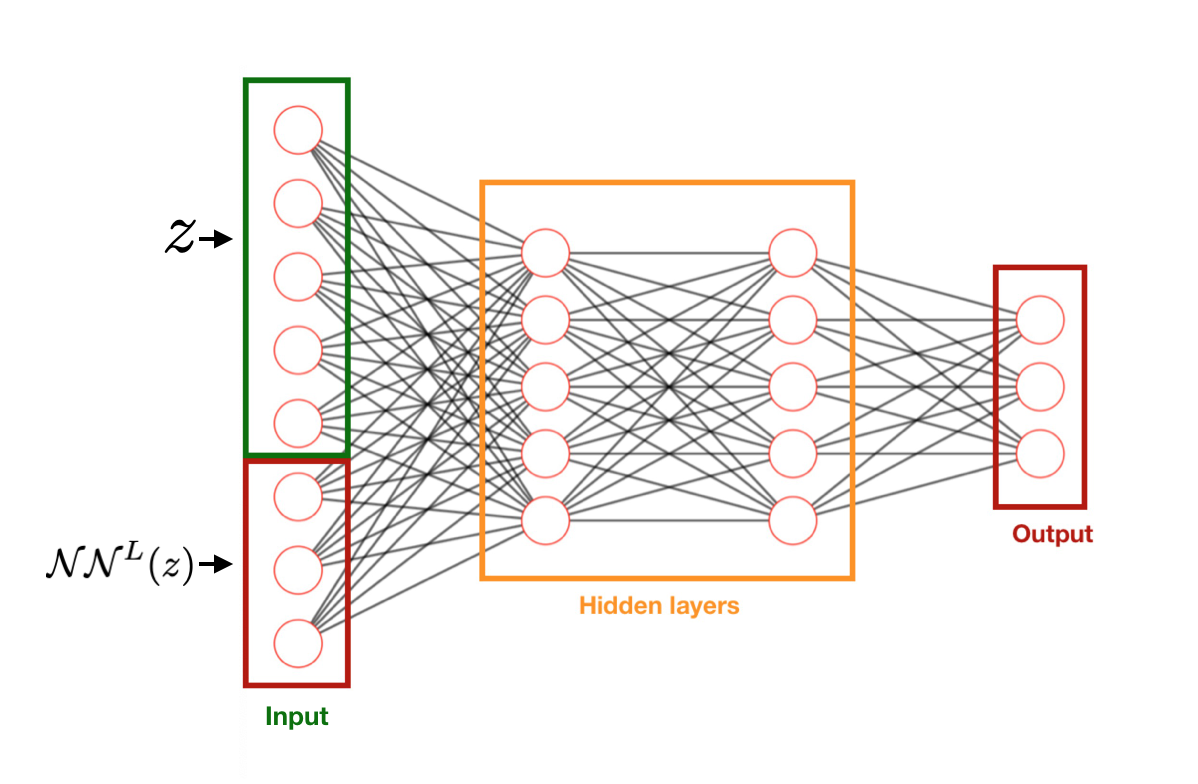}
  \end{overpic}
  \end{center}
\caption{A composite architecture of DNN that combine a low fidelity DNN. The input is the concatenation of the parameter vector $z$ and the low-fidelity data $\mathcal{NN}^L(z)$, and the output is the high-fidelity $f^H(z)$.}\label{MDNN_structure}
  \end{figure}

\begin{algorithm}[t]
  \caption{Multi-fidelity Composite Deep Neural Networks}
  \label{alg:MDNN}
  \begin{algorithmic}[1]
    \Require
    The trained DNN model $\mathcal{NN}^L(z)$; the high-fidelity model $f^H(z)$;

  \State  Choose additional training samples $\{z_k\}^Q_{k=1}$.
  \State  Run the low-fidelity model $\mathcal{NN}^L$ and high-fidelity model $f^H$ for each $z_k$ to obtain $\mathcal{NN}^L(z_k)$  and $f^H(z_k).$
  \textcolor{black}{
\State Construct the training set $\mathcal{D}=\Big\{\Big( \big(z_k, \mathcal{NN}^L(z_k)\big), f^H(z_k)\Big)\Big\}_{k=1}^Q$.
  \State Train a new neural network $\mathcal{NN}^H(z;\theta)$ with input $\big(z, \mathcal{NN}^L(z)\big)$ and output $f^H(z)$ by using the training set $\mathcal{D}$. }
  \end{algorithmic}
\end{algorithm}

The main idea here is to discover and exploit the correlations between low- and high-fidelity data \cite{Peherstorfer2016survey}.  Suppose we have a high-fidelity model $f^H$ (here we simply assume that $f^H$ is the true forward model $f$) and a low-fidelity surrogate $f^L$ (e.g., a trained DNN-surrogate $\mathcal{NN}^L$ in our framework). Then, we aim at learning a nonlinear map $\mathcal{F}$ between the two models:
\begin{equation}\label{eq:nonlinear1}
 f^H(z)=\mathcal{F}\big(z, f^L(z)\big)=\mathcal{F}\big(z, \mathcal{NN}^L(z)\big).
\end{equation}
This can be done by approximating the nonlinear map $\mathcal{F}$ using DNNs, i.e.,
\begin{equation}\label{eq:nonlinear2}
 f^H(z)\approx \mathcal{NN}^H(z; \theta):=\mathcal{NN}\big(z, \mathcal{NN}^L(z); \theta\big).
\end{equation}
The new involved parameters $\theta$ can be trained by using additional high fidelity data $\Big\{\Big( (z_k, \mathcal{NN}^L(z_k)), f^H(z_k)\Big)\Big\}_{k=1}^Q,$ and this yields a composite DNNs as illustrated in Fig. \ref{MDNN_structure}.

The key idea here is to view the low-fidelity model $\mathcal{NN}^L$ as an input into the high fidelity DNN, motivated by the fact that models are highly correlated. By doing this, one can expect to use fewer layers (or even use a shallow neural network) for constructing $\mathcal{NN}^H$ -- resulting much reduced computational complexity in the training procedure. The detailed step for constructing $\mathcal{NN}^H$ are summarized in Algorithm \ref{alg:MDNN}.

At this stage, one may ask why not just include those additional high fidelity data $\Big\{\Big( z_k, f^H(z_k)\Big)\Big\}_{k=1}^Q$ when training $\mathcal{NN}^L$ (--yielding a better surrogate)?  The answer is that, in the beginning, one can only construct the surrogate with prior-based information, \textcolor{black}{and this surrogate may not be accurate enough even if large sample evaluations are used}, as the posterior density is in general concentrate to a small region.  Thus, we aim at adaptively correct the surrogate online using local data. Details will be presented in the next section.

\textsc{remark}. We close this section by remarking that in \cite{Yan+Zhou19JCP}, a correction procedure for the PCE-based surrogate is discussed. In particular, the following correction technique is presented:
\begin{equation}
 f^H(z)\approx f^H_{\textmd{PCE}}(z):= f^L_{\textmd{PCE}}(z) + f_{\textmd{CORR}}(z),
\end{equation}
where $f^L_{\textmd{PCE}}(z)$ is the low-fidelity PCE-surrogate, and $f_{\textmd{CORR}}(z)$ is the correction term that is determined online by additional high fidelity data. Notice that this is a \textit{linear} correction, as $f^H_{\textmd{PCE}}$ and $f^L_{\textmd{PCE}}$ are linearly dependent. While in the present work, the correction procedure (\ref{eq:nonlinear1})-(\ref{eq:nonlinear2}) can learn a nonlinear correlation between models (which is in general the case for practical applications). \textcolor{black}{In Section \ref{sec:tests}, we shall perform the comparison between the current approach and the PCE-based surrogate approach proposed in \cite{Yan+Zhou19JCP} by numerical examples.}

\subsection{An adaptive procedure for correcting the surrogate}

We now propose an adaptive procedure to correct the surrogate in the MCMC framework. The procedure begins with a low fidelity model $\mathcal{NN}^L$ which is constructed offline. Then, for the online computations, an adaptive sampling framework is used to construct and refine the surrogate model $\mathcal{NN}^H$, following the idea in our previous work \cite{Yan+Zhou19JCP}. The procedure consists of the following steps:
\begin{itemize}
\item  Initialization:  build a low fidelity model $\mathcal{NN}^L$. Set $f^L = \mathcal{NN}^L$.
\item Online computations: using the surrogate  $f^L$, run the MCMC  (e.g., a standard Metropolis-Hastings (MH) algorithm) to sample the approximated posterior distribution for a certain number of steps (say 1000 steps). The last state, denoted $z^{-}$, will be used to \textcolor{black}{propose a candidate $z^{+}$ with a proposal density $q$}.
\item Indicator for refinement:  generate an accept sample $\tilde{z}$ using high fidelity information on $z^-$ and $z^+.$  Then compute the difference between the high fidelity model $f^H$ and the surrogate $f^L$ at $\tilde{z}$. If the difference is bigger than a given tolerance, then one generates new high fidelity data to correct the model $\mathcal{NN}^H$ using Algorithm \ref{alg:MDNN}. Set $\mathcal{NN}^H$ as the new surrogate, i.e., $f^L=\mathcal{NN}^H$.
\item  Use the surrogate $f^L$ to accept/reject the proposal $z^+$.
\item  Repeated the above procedure for many times (say at most $I_{max}$ times). Finally the posterior samples can be generated by gathering all the samples in the above procedures.
\end{itemize}

Choosing the indicator for correcting the surrogate $f^L$ is the most critical issue in the above procedure.  As mentioned, we first sampling approximate posterior distribution based on the surrogate model $f^L$ for a certain number of steps using a standard MH algorithm. The goal is to generate several samples  that the initial sample points and the last point are uncorrelated. Similar to the MH algorithm, we compute an acceptance probability using the high fidelity model $f^H$ (or the true forward model ):
\begin{equation}\label{acceppro}
\beta = \min\Big\{1, \frac{\mathcal{L}\big(d, f^H(z^-)\big)\pi(z^-)}{\mathcal{L}\big(d, f^H(z^{+})\big)\pi(z^{+})} \Big\}.
\end{equation}
Using this parameter $\beta$, we can obtain an accept point $\tilde{z}$,  which is expected to be much closer to the posterior region.

\begin{algorithm}[th]
  \caption{Indicator for correcting the surrogate}
  \label{alg:upMNN}
  \begin{algorithmic}[1]
   \Require
  Set an error threshold $tol$ and the radius $R$. Given $z^{-}$ and $z^+$, we do the following steps:

 \State  Compute acceptance probability using high-fidelity model $f^H$
\begin{equation*}
\beta = \min\Big\{1, \frac{\mathcal{L}\big(d, f^H(z^-)\big)\pi(z^-)}{\mathcal{L}\big(d, f^H(z^{+})\big)\pi(z^{+})} \Big\}
\end{equation*}
\State Draw $s \sim \mathcal{U}(0,1)$. If $s<\beta$, let $\tilde{z}=z^-$, otherwise $\tilde{z}=z^{+}$.
 \State Compute the relative error $err(\tilde{z}) $ using Eq. (\ref{conerr})
      \If {$err(y) > \epsilon$}
      we generate $Q$ random points locally $\{z_i\} \in B(\tilde{z},R)$  and correct the surrogate model to get $\mathcal{NN}^H(z)$ using Algorithm \ref{alg:MDNN}
    \EndIf
  \State
\Return Set $\mathcal{NN}^H(z)$ as the new low-fidelity surrogate model.
  \end{algorithmic}
\end{algorithm}

Once we obtain the accept candidate point $\tilde{z}$,  we then compute the relative error
\begin{equation}\label{conerr}
err(\tilde{z}) = \frac{\|f^H(\tilde{z})-f^L(\tilde{z})\|_{\infty}}{\|f^H(\tilde{z})\|_{\infty}}.
\end{equation}
If this error indicator exceeds a user-given threshold $tol$, we shall generate $Q$ random points $\{z_i\}_{i=1}^Q$ locally around $\tilde{z}$ (say, uniformly sampling in a local ball $B(\tilde{z},R):=\{z: \|z-\tilde{z}\|_{\infty} \leq R\}$) to correct the surrogate model $f^L$ using Algorithm \ref{alg:MDNN}. While if the error indicator is smaller than $tol$, it means that the low-fidelity model is still acceptable and we just go ahead. The detailed procedure for updating the surrogate model is summarized in Algorithm \ref{alg:upMNN}.

In the correction procedure, we can not afford too many high fidelity simulations, that is, $Q$ can not be too large. Consequently, to avoid over fitting for the online training procedure, we limit ourselves to use a DNN with at most two hidden layers (or even consider shallow networks).  Algorithm description of the MH approach using locally adapted multi-fidelity DNN-surrogate is  summarized in  Algorithm \ref{alg:AMNN}.

\begin{algorithm}[th]
  \caption{Adaptive multi-fidelity DNN-based MH algorithm}
  \label{alg:AMNN}
  \begin{algorithmic}[1]
  \Require
  Given the initial surrogate $f^L=\mathcal{NN}^L$ and a proposal density $q,$  we fix a stopping indicator $m$ for the MCMC sampling, i.e., we shall stop to check if the correction is needed for each $m$-steps, for instance, we can choose $m=1000.$ We also set a maximum number $I_{max}$ of corrections, that is, we shall at most correct the surrogate for $I_{max}$ time so that one can control the total computational complexity.
 \State  Choose a starting points $z_0$; let $X_0=\{\}$;
 \For {$n=1,\cdots, I_{max}$}
\State Draw $m-1$ samples $\{z_1,\cdots,z_{m-1}\}$ from the approximate posterior based on $f^L$.  Propose $z^*\sim q(\cdot|z_{m-1})$.
  \State  If the surrogate needs refinement near $z^*$ or $z_{m-1}$, then select new samples locally to correct the surrogate to get $\mathcal{NN}^H$ using Algorithm \ref{alg:upMNN}. Set $f^L = \mathcal{NN}^H$.
\State  Compute acceptance probability
\begin{equation*}
\alpha = \min\Big(1, \frac{\mathcal{L}\big(d, f^L(z^*)\big)\pi(z^*)}{\mathcal{L}\big(d,f^L(z_{m-1})\big)\pi(z_{m-1})} \Big)
\end{equation*}
   \State Draw $s \sim \mathcal{U}(0,1)$. If $s<\alpha$, let $z_m=z^*$, otherwise $z_m=z_{m-1}$.
   \State Let $z_0 =z_m$ and $X_n=X_{n-1} \bigcup \{z_1,\cdots,z_m\}$
   \EndFor  \State
    \Return Posterior samples $ X_{I_{max}}$
  \end{algorithmic}
\end{algorithm}
To summary, our approach departs from an offline constructed DNN-surrogate, and then we correct the multi-fidelity  DNN-surrogate adaptively using locally generated high fidelity data. The key idea is to consider the composite DNNs in which the previous trained DNN model $\mathcal{NN}^L$ is viewed as an input variable in the next updated surrogate. The locally generated training data are then expected to concentrate to the high probability region of the posterior density. The online training procedure is also expected to be efficient due to the correlation of two models.

\section{Numerical Examples}\label{sec:tests}

In this section, we present a benchmark elliptic PDE inverse problem to illustrate the accuracy and efficiency of the proposed adaptive multi-fidelity DNN (ADNN) algorithm.  We describe three examples in which ADNN algorithm produce accurate posterior samples using dramatically fewer evaluations of the forward model than the conventional MCMC.

\textcolor{black}{The first example will be chosen for comparison reason, i.e., we shall compare the performance of the current approach with the PCE-based surrogate approach. While two additional examples are constructed to demonstrate the efficiency and the accuracy of the present approach.}

\textcolor{black}{In all our numerical tests, the ADNN approach was performed with a self-written program which was coded by MATLAB. The optimization procedure is carried out by the Adam algorithm as mentioned before. The learning rate is set to be $\epsilon =10^{-3}$, and the hyper-parameter values of Adam are chosen based on default recommendations as suggested in \cite{Kingma2014Adam}.}  For each of these examples, we run the ADNN algorithm for $I_{max} =50$ iterations, with subchain length $m=1,000$.  Unless otherwise specified, we shall use the following parameters $tol=0.1, R=0.2$ in ADNN.
To make a fair comparison, we run the prior-based DNN method described in Section \ref{sec:DNN}  for $50,000$ iterations.  The MCMC simulation using the high-fidelity model is also conducted, and its results are used as the reference to evaluate accuracy and efficiency of the two methods. For both algorithms, the same fixed Gaussian proposal distribution is used, and  the last $30,000$ realizations are used to compute the relevant statistical quantities.  All the computations were performed using MATLAB 2015a on an Intel-i5 desktop computer.

\subsection{Problem setup}

We consider the problem of inferring subsurface permeability from a finite number of noisy pressure head measurements \cite{Cui2014data,Yan+Zhou19JCP}.  More specifically, consider a domain of interest $\Omega=[0,1]^2$, and let  $u(x)$ be pressure head,  is the solution of an elliptic PDE in two spatial dimensions
\begin{eqnarray}\label{2dellip}
\begin{array}{rl}
-\nabla\cdot(\kappa(x) \nabla u(x))&=f(x),\quad x\in \Omega,\\
 u(x)&=0, \quad  \,x\in \partial{\Omega}.
 \end{array}
\end{eqnarray}
The data $d$ is given by a finite set of $u$, perturbed by noise, and the problem is to recover the permeability $\kappa(x)$ from these measurements.  In what follows, we choose the source $f(x) = 100\sin(\pi x_1)\sin(\pi x_2)$.

In the numerical simulation, we  solve the equation (\ref{2dellip}) using a spectral approximations  with  polynomial degree $P=6$.  In order not to commit an 'inverse crime', we generate the data by solving the forward problem using a higher order (P=10) than that is used in the inversion.

\subsection{Example 1: a nine-dimensional inverse problem}

  \begin{figure}
\begin{center}
    \begin{overpic}[width=0.45\textwidth,trim= 20 0 20 15, clip=true,tics=10]{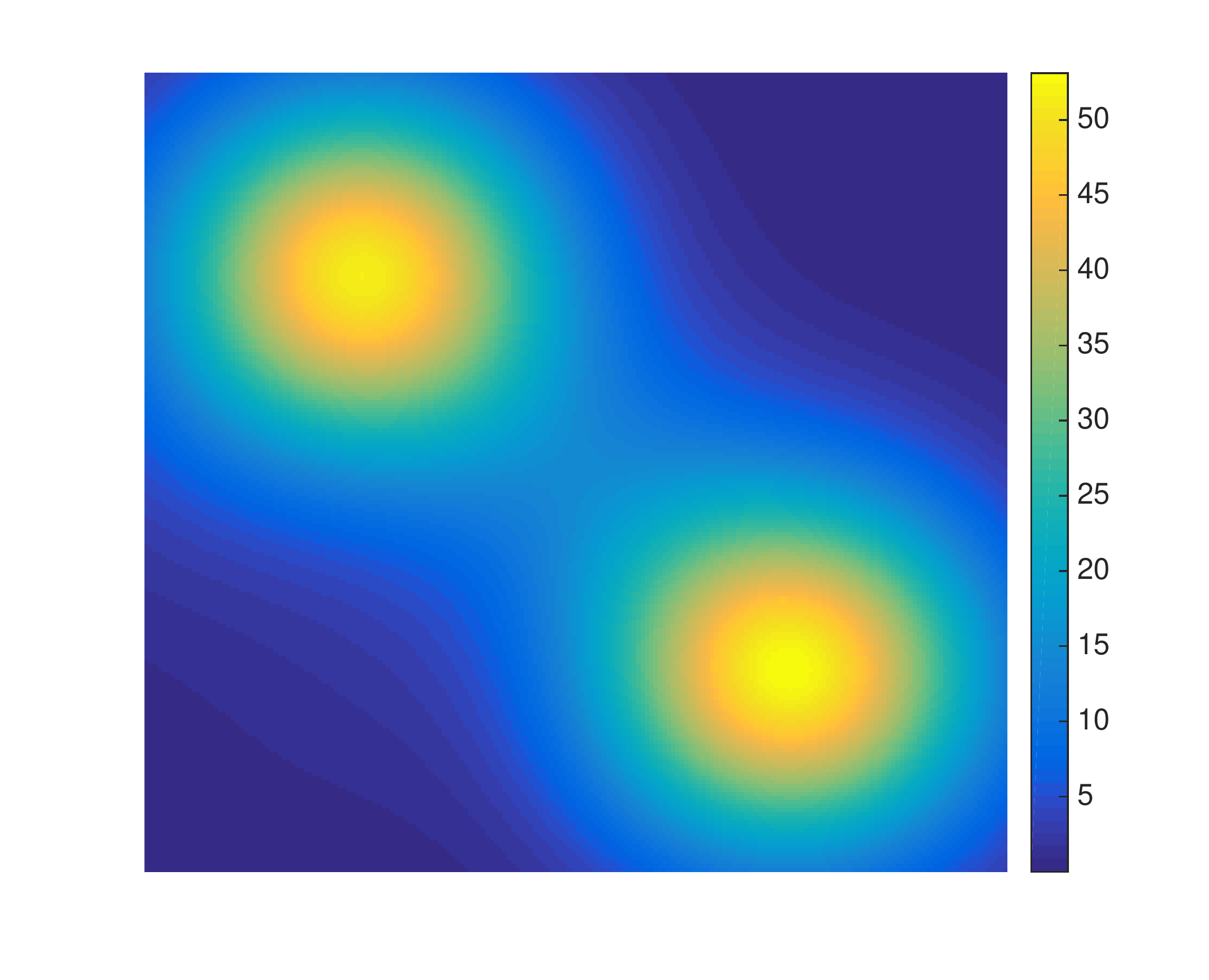}
      \end{overpic}
   \begin{overpic}[width=0.45\textwidth,trim= 20 0 20 15, clip=true,tics=10]{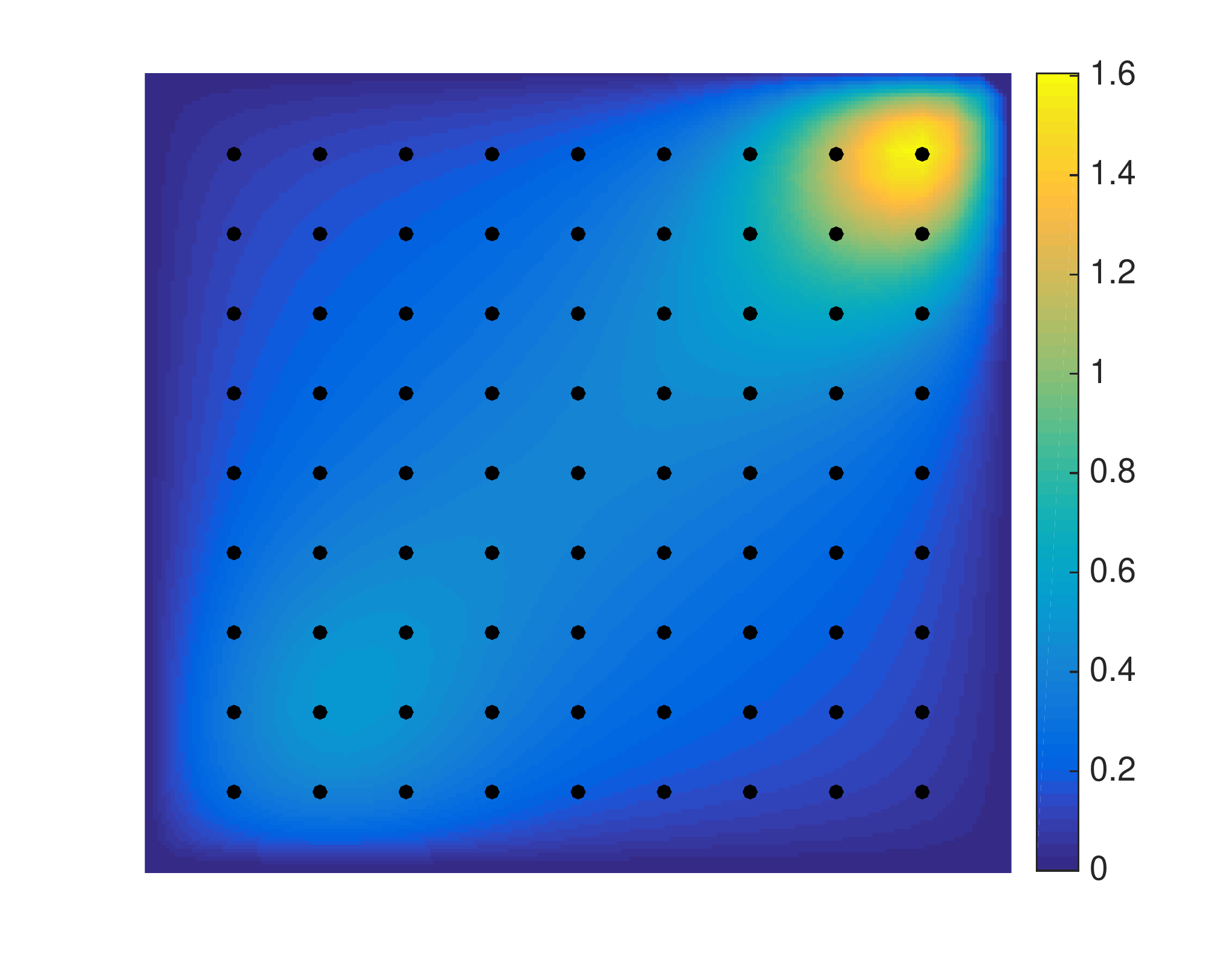}
  \end{overpic}
\end{center}
\caption{Example 1:  Setup of the test case for Example 1. Left: the true permeability used for generating the synthetic data sets. Right: the model outputs of the true permeability.  The figure is adopted from \cite{Yan+Zhou19JCP}. }\label{exact_eg1}
\end{figure}

In the first example, the permeability field $\kappa(x)$ is defined by
\begin{eqnarray*}
\kappa(x)=\sum^{9}_{i=1}\kappa_i \exp(- 0.5\frac{\|x-x_{0,i}\|^2}{0.15^2}),
\end{eqnarray*}
where $\{x_{0,i}\}^{9}_{i=1}$ are the centers of the radial basis function, and the weights $\{\kappa_i\}^9_{i=1}$ are parameters in the Bayesian inverse problem.

This example is a typical benchmark problem considered in Refs. \cite{Cui2014data,Yan+Zhou19JCP} and it is investigated here for comparison purpose.  To this end, we use the same model setup and synthetic data used in \cite{Yan+Zhou19JCP}.  More precisely,  the prior distributions on each of the weights $\kappa_i, i=1,\cdots, 9$ are independent and log-normal; that is, $\log(\kappa_i)\sim N(0,1)$.  The true parameter is drawn from  $\log(\kappa_i)\sim U(-5,5)$, and the true permeability field used to generate the test data is shown in Fig.\ref{exact_eg1}.  The synthetic data $d$ is generated by selecting  the values of the states at a uniform $9\times 9$ sensor network
\begin{eqnarray*}
d_j=u(x_j)+\max_{j}\{|u(x_j)|\}\delta\xi_j,
\end{eqnarray*}
where $\delta$ dictates the relative noise level and $\xi_j$  is a  Gaussian random variable with zero mean and unit standard deviation.  In the following, unless  otherwise specified, we set $\delta =0.05$.

\subsubsection{Comparison of approximations}\label{sec:1}

We compare the posterior approximations obtained by three types of approaches:
\begin{itemize}
\item The {\it conventional MCMC}, or the direct MCMC approach based on the forward model evaluations.
\item The MCMC approach based on a prior-based DNN surrogate model evaluations.
\item The ADNN approach presented in Section \ref{sec:method}.
\end{itemize}
In our figure and results, we will use ``Direct" to denoted the conventional MCMC, ``DNN" to denoted the prior-based DNN approach, and ``ADNN" to denote the ADNN algorithm.  For the ADNN algorithm, we first construct a prior-based DNN surrogate $\mathcal{NN}^L$ using $N=50$ training points with 4 hidden layers and 40 neurons per layer.   Using this DNN model as low-fidelity model, we can construct and refine a multi-fidelity model $\mathcal{NN}^H$ via Algorithm \ref{alg:AMNN}.  Especially, when the error indicator $err(\tilde{z})$ exceeds the threshold $tol$, we choose $Q=10$ random points in a local set $B(\tilde{z},R)$ to refine the multi-fidelity NN surrogate $\mathcal{NN}^H$. Here, one hidden layer with $50$ neurons are used in $\mathcal{NN}^H$. In this example, the regularization parameter $\lambda$ is set to 0, i.e.,  no regularization is used.

\begin{figure}
\begin{center}
 \begin{overpic}[width=\textwidth,trim=20 0 20 15, clip=true,tics=10]{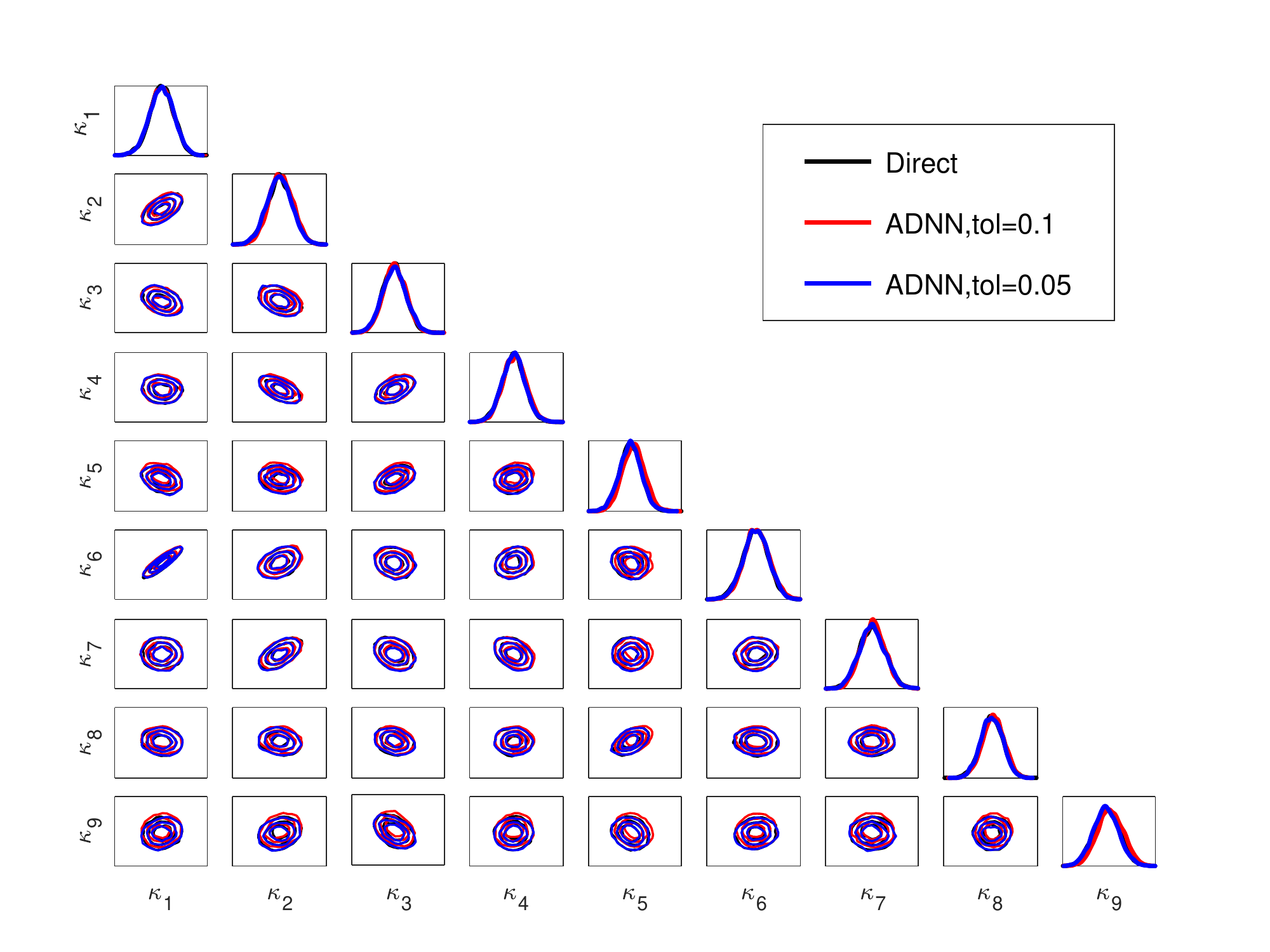}
 \end{overpic}
\end{center}
\caption{Example 1: One- and two-dimensional posterior marginals of the nine parameters. Black line: the Direct method.  Red line: the ADNN algorithm with $tol=0.1$; Blue line: the ADNN algorithm with $tol=0.05$. }\label{pos_contour_eg1}
\end{figure}

\begin{figure}
\begin{center}
  \begin{overpic}[width=0.32\textwidth,trim=20 0 20 15, clip=true,tics=10]{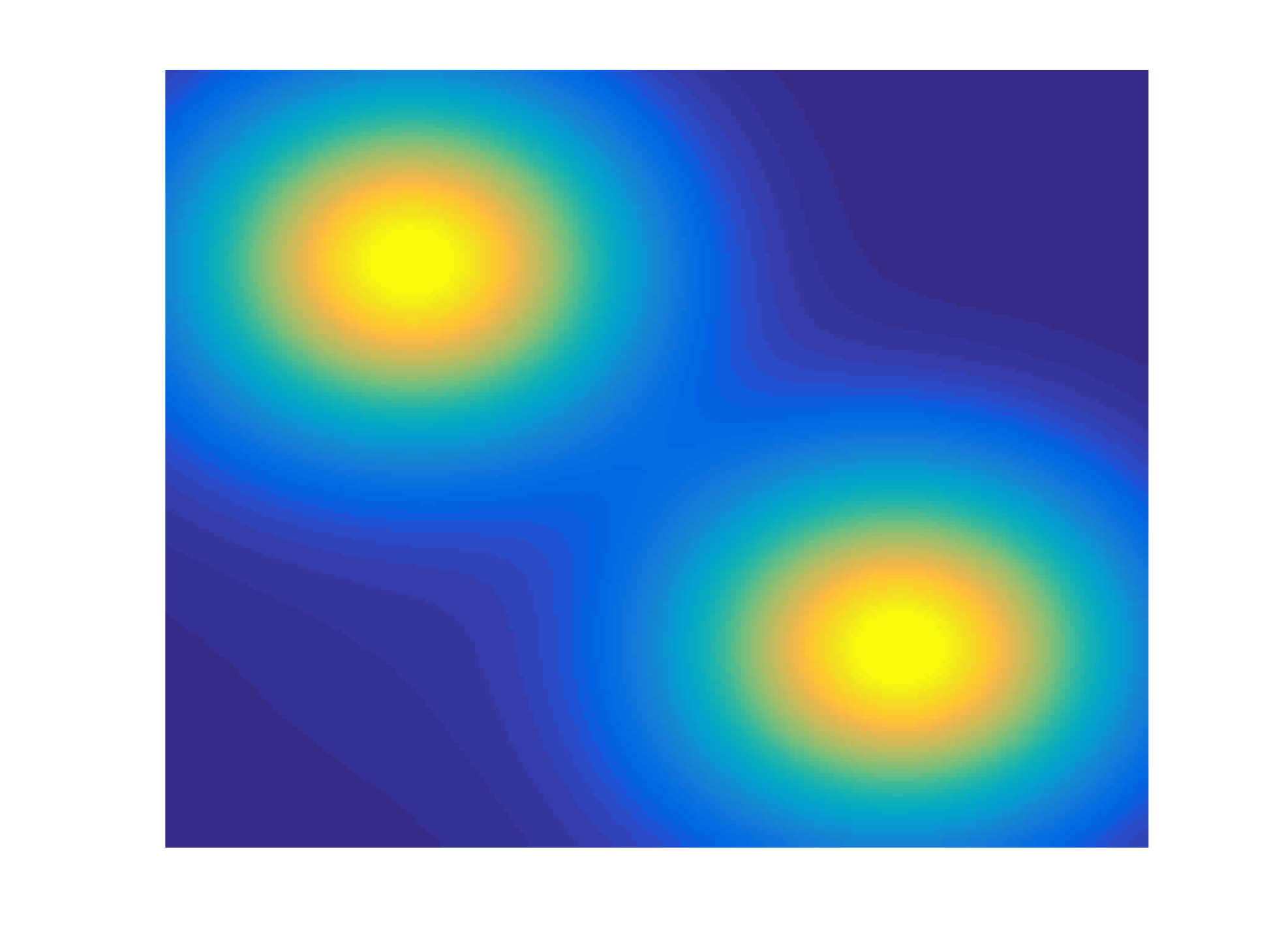}
  \end{overpic}
    \begin{overpic}[width=0.32\textwidth,trim= 20 0 20 15, clip=true,tics=10]{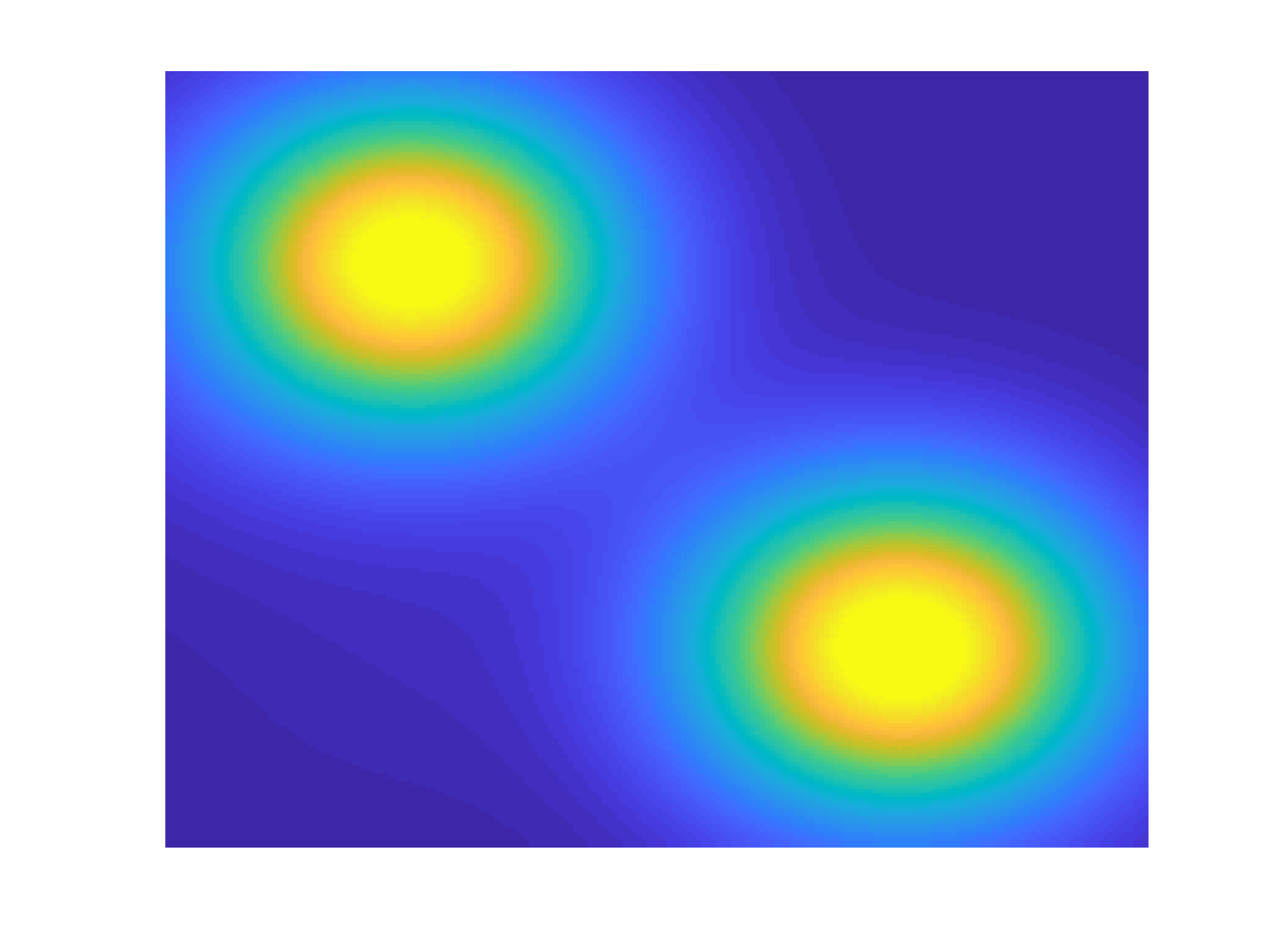}
  \end{overpic}
    \begin{overpic}[width=0.32\textwidth,trim= 20 0 20 15, clip=true,tics=10]{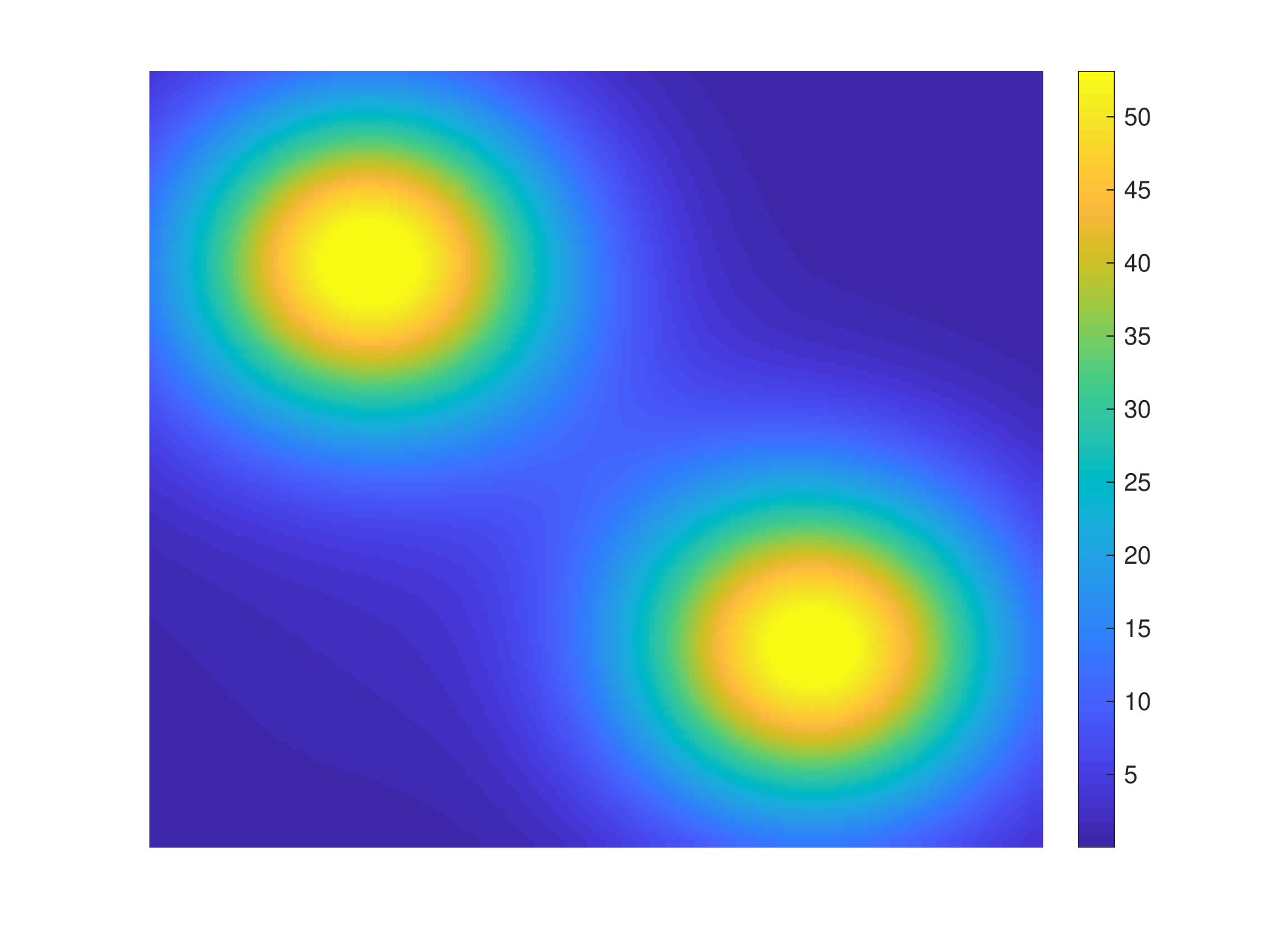}
  \end{overpic}
 \begin{overpic}[width=0.32\textwidth,trim=20 0 20 15, clip=true,tics=10]{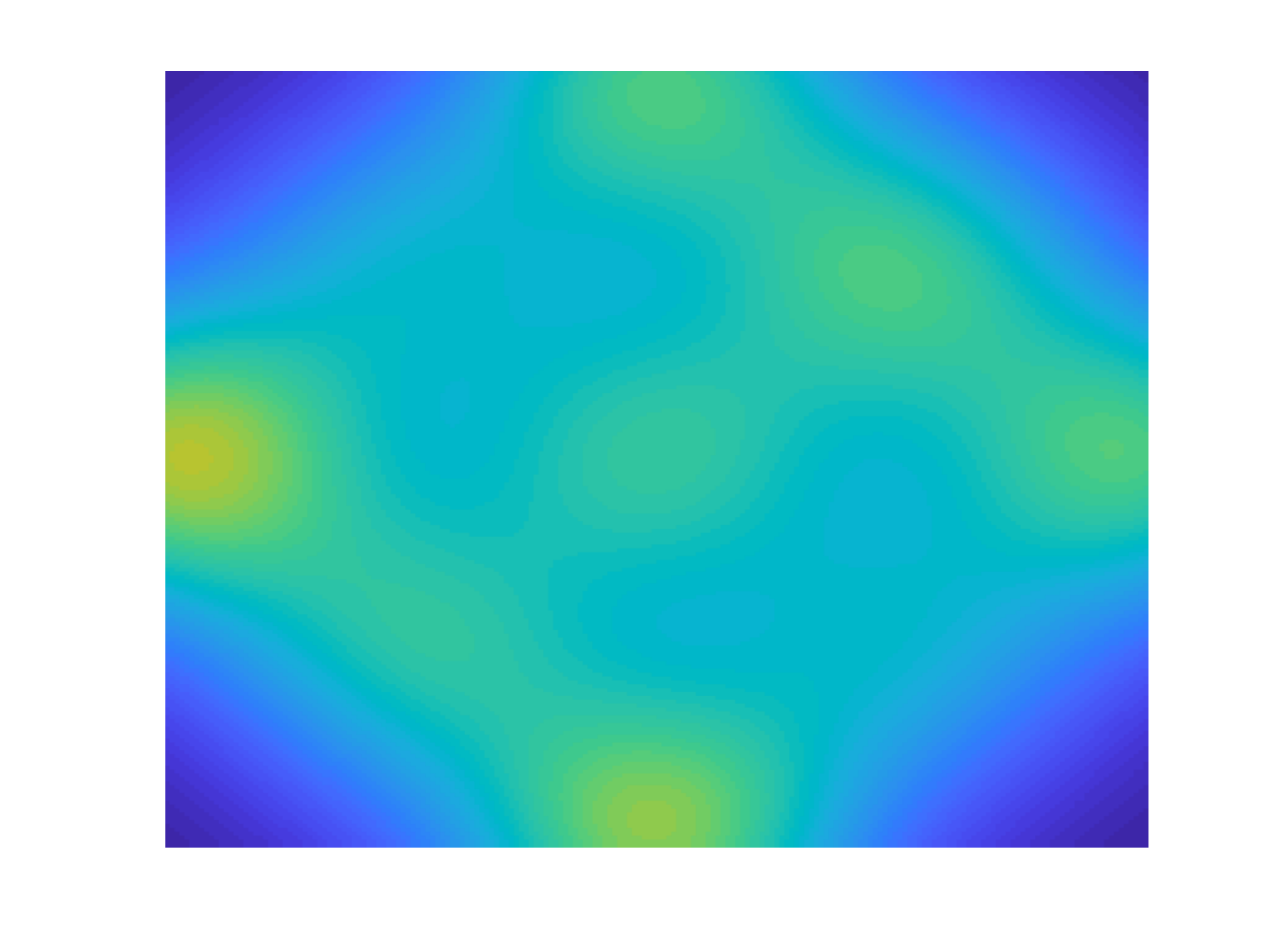}
   \put (43,-3) {\scriptsize Direct}
  \end{overpic}
    \begin{overpic}[width=0.32\textwidth,trim= 20 0 20 15, clip=true,tics=10]{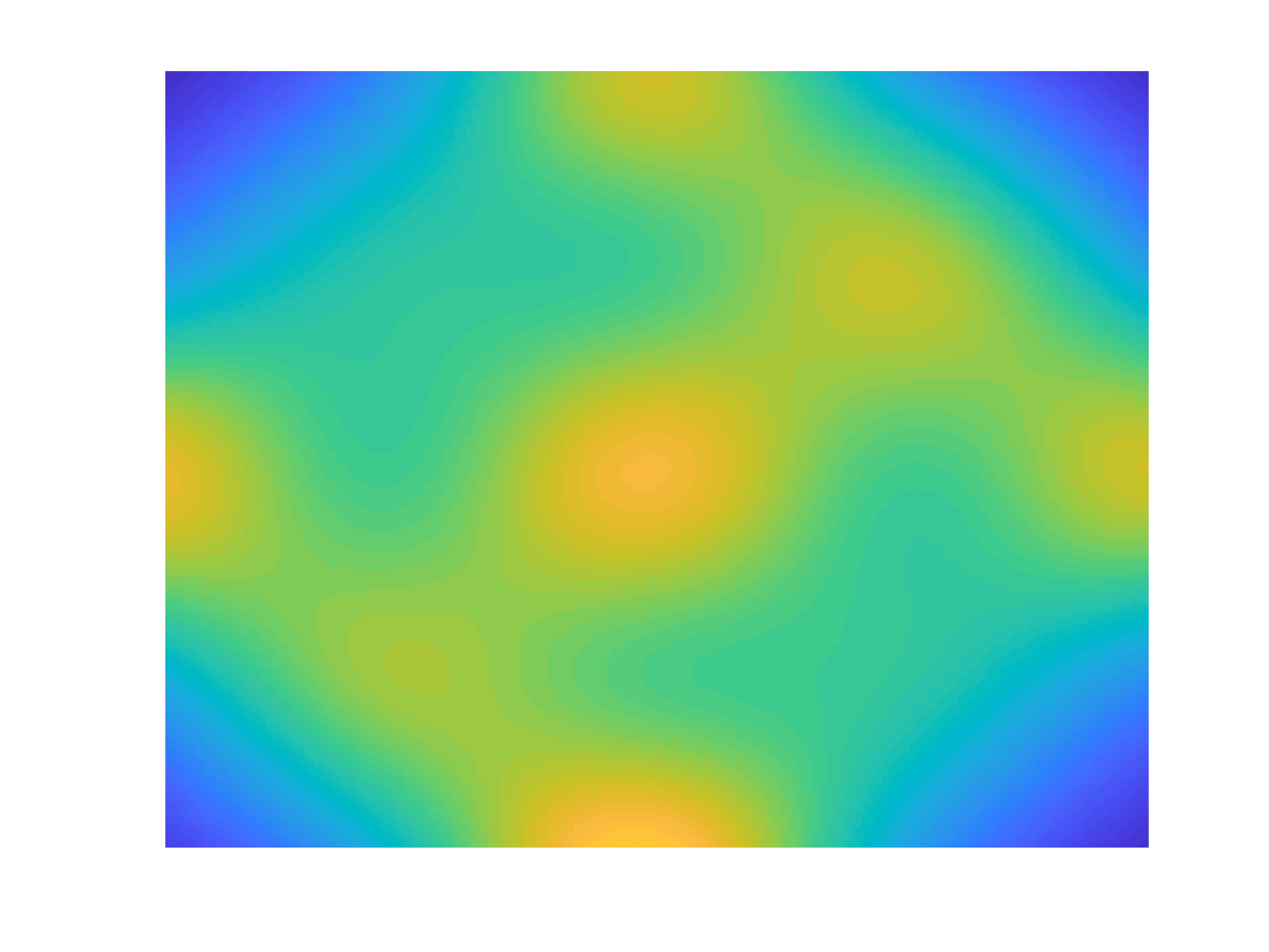}
      \put (28,-3) {\scriptsize ADNN, $tol=0.1$}
  \end{overpic}
    \begin{overpic}[width=0.32\textwidth,trim= 20 0 20 15, clip=true,tics=10]{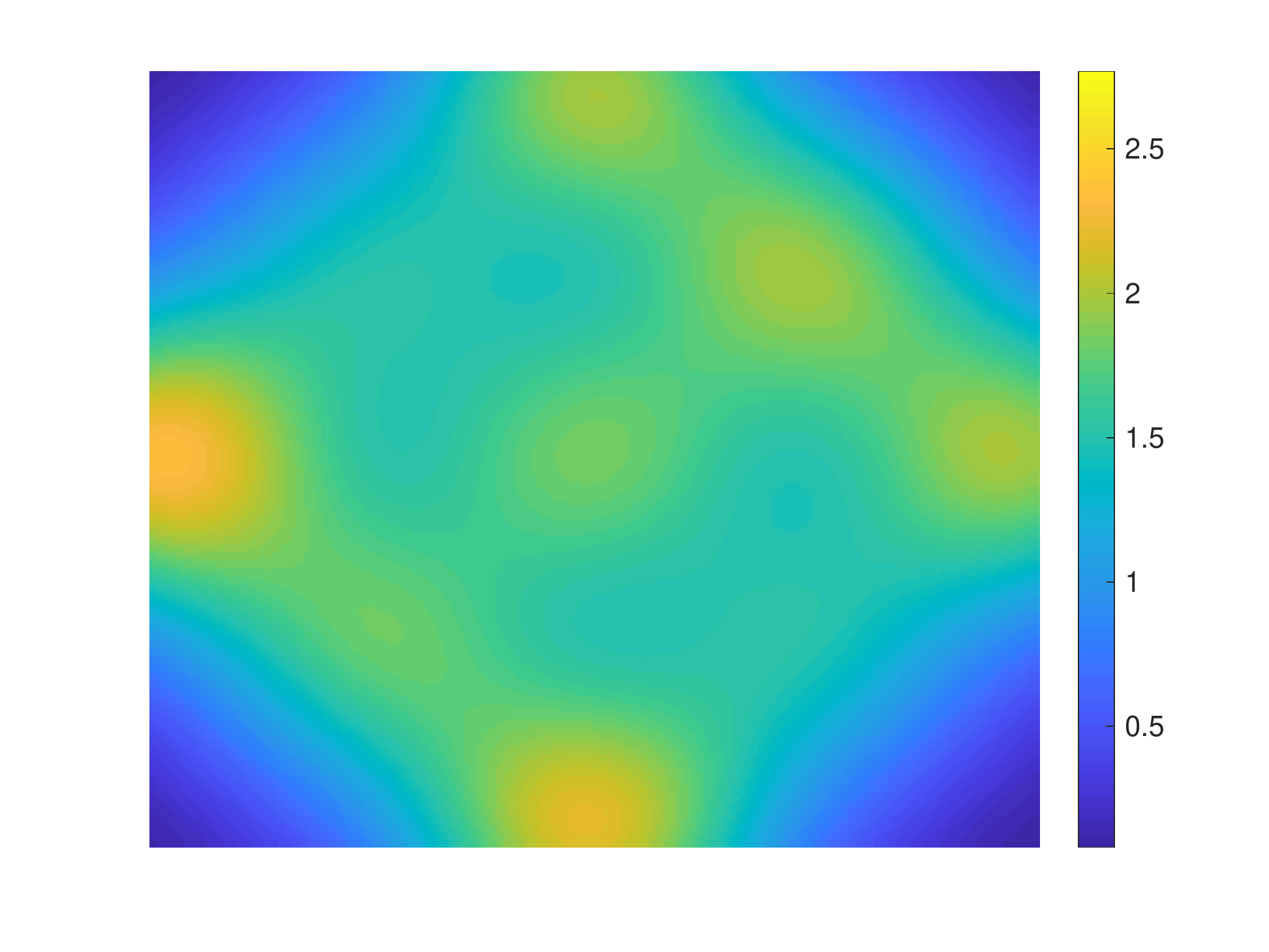}
      \put (23,-3) {\scriptsize ADNN, $tol=0.05$}
  \end{overpic}
\end{center}
\caption{Example 1:  Posterior mean (top) and  posterior standard deviation (bottom) arising from direct MCMC, ADNN ($tol=0.1$) and  ADNN ($tol=0.05$), respectively. The relative noise  level $\delta$ is $0.01$.}\label{adap_sol_eg1-1}
  \end{figure}

\begin{figure}
\begin{center}
 \begin{overpic}[width=0.32\textwidth,trim=35 10 45 5, clip=true,tics=10]{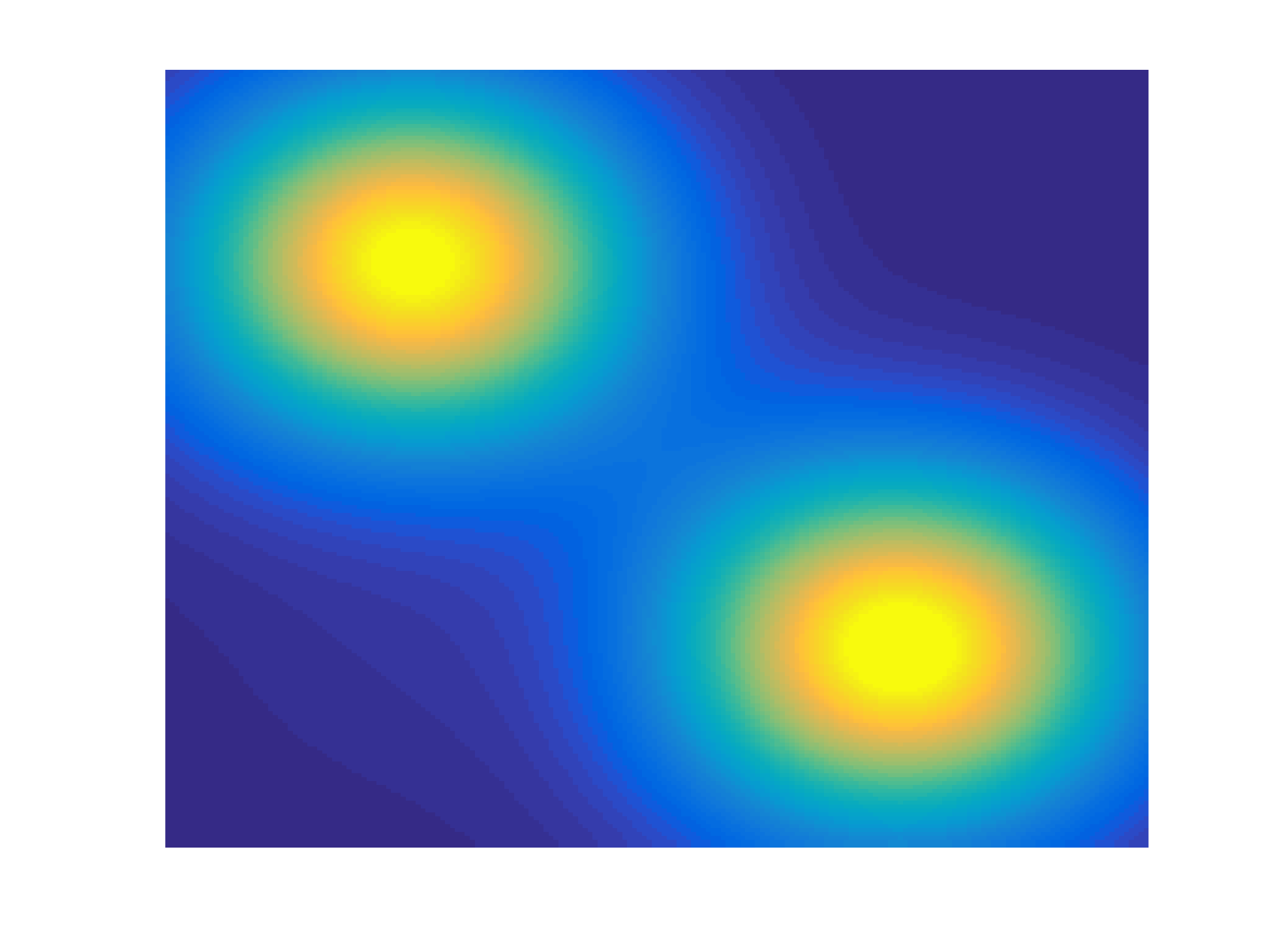}
  \end{overpic}
    \begin{overpic}[width=0.32\textwidth,trim= 35 10 45 5, clip=true,tics=10]{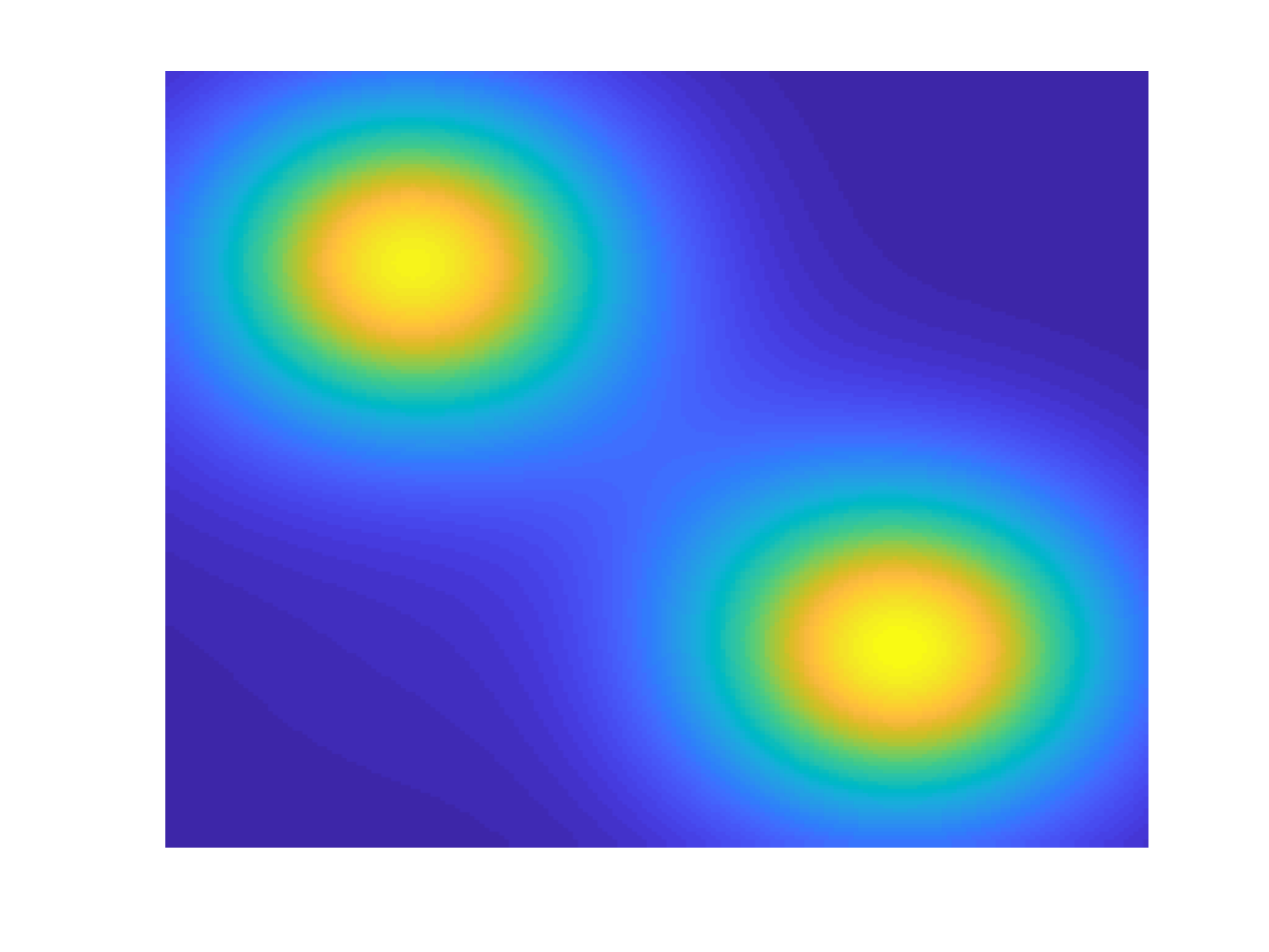}
  \end{overpic}
    \begin{overpic}[width=0.32\textwidth,trim= 35 10 45 5, clip=true,tics=10]{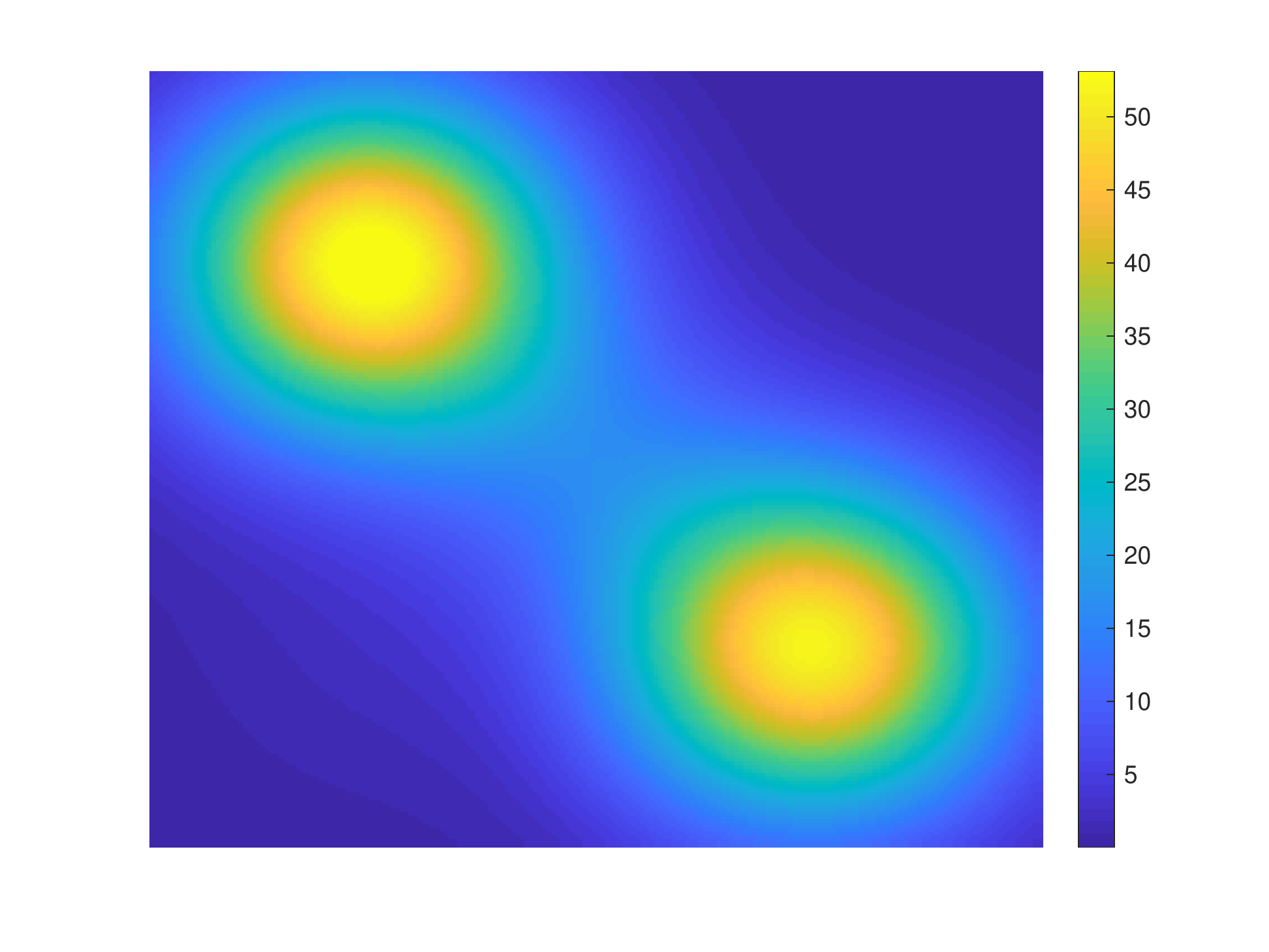}
  \end{overpic}
 \begin{overpic}[width=0.32\textwidth,trim=35 10 45 5, clip=true,tics=10]{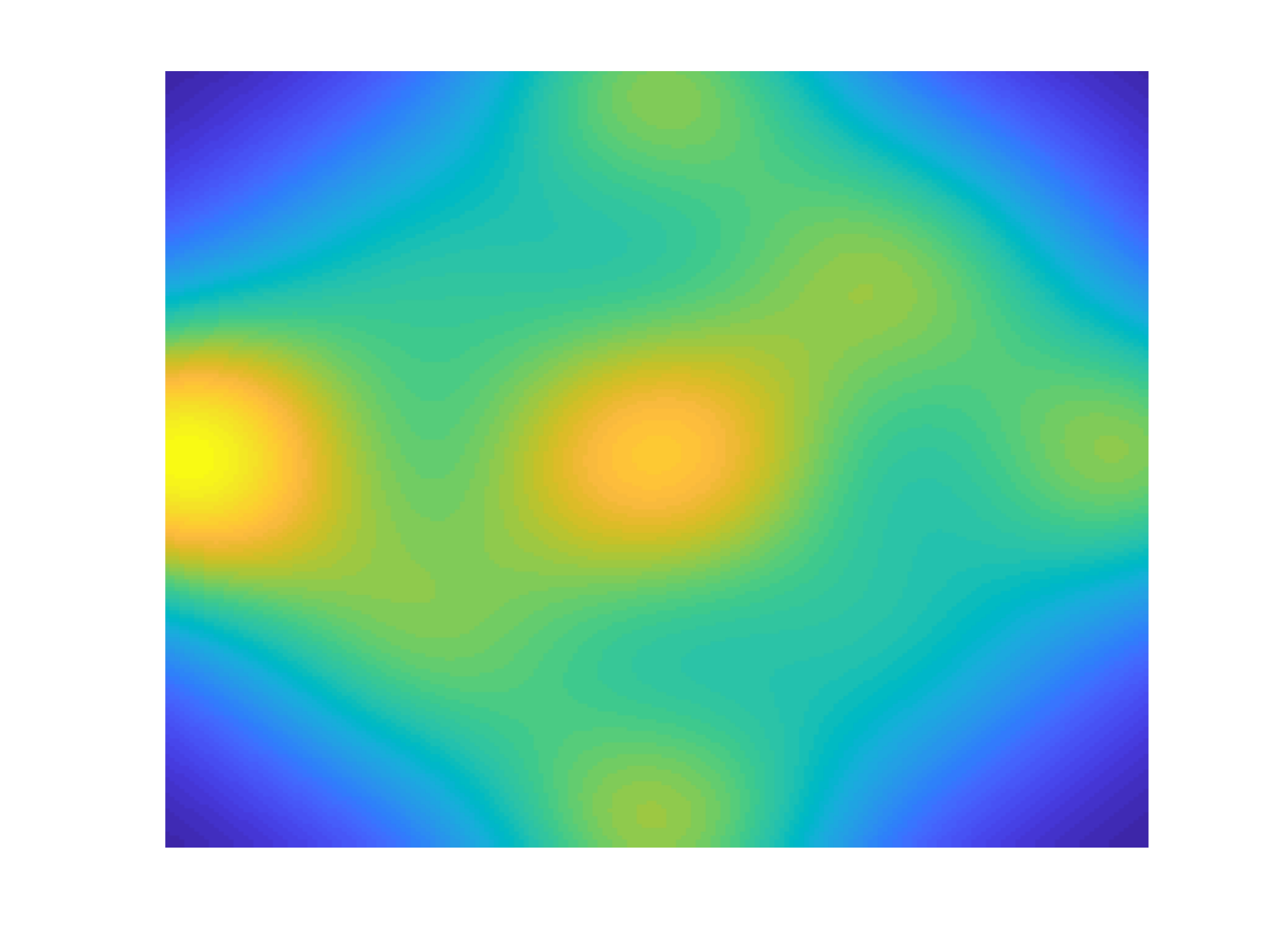}
 \put (43,-3) {\scriptsize Direct}
  \end{overpic}
    \begin{overpic}[width=0.32\textwidth,trim= 35 10 45 5, clip=true,tics=10]{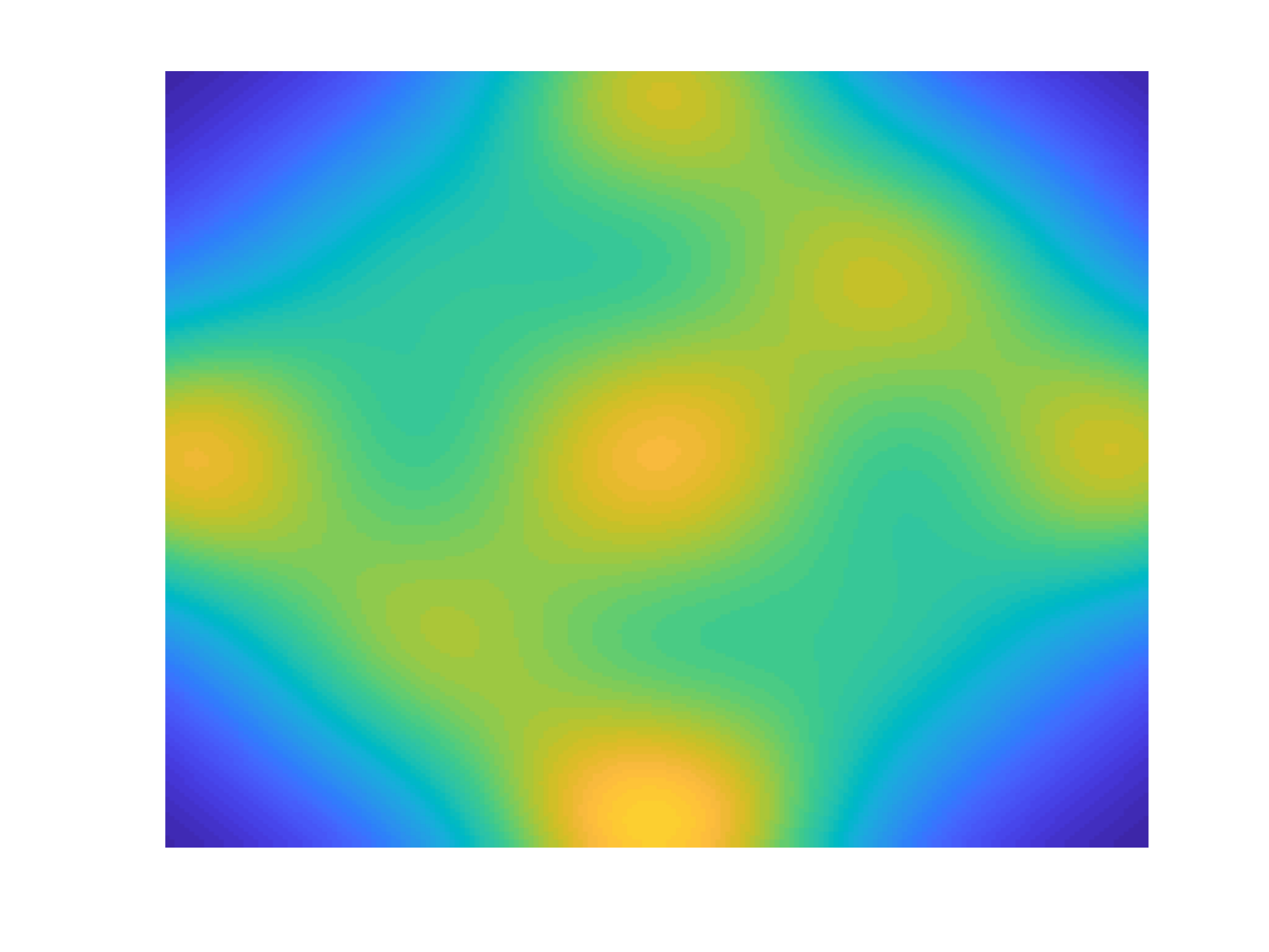}
       \put (28,-3) {\scriptsize ADNN, $tol=0.1$}
  \end{overpic}
    \begin{overpic}[width=0.32\textwidth,trim= 35 10 45 5, clip=true,tics=10]{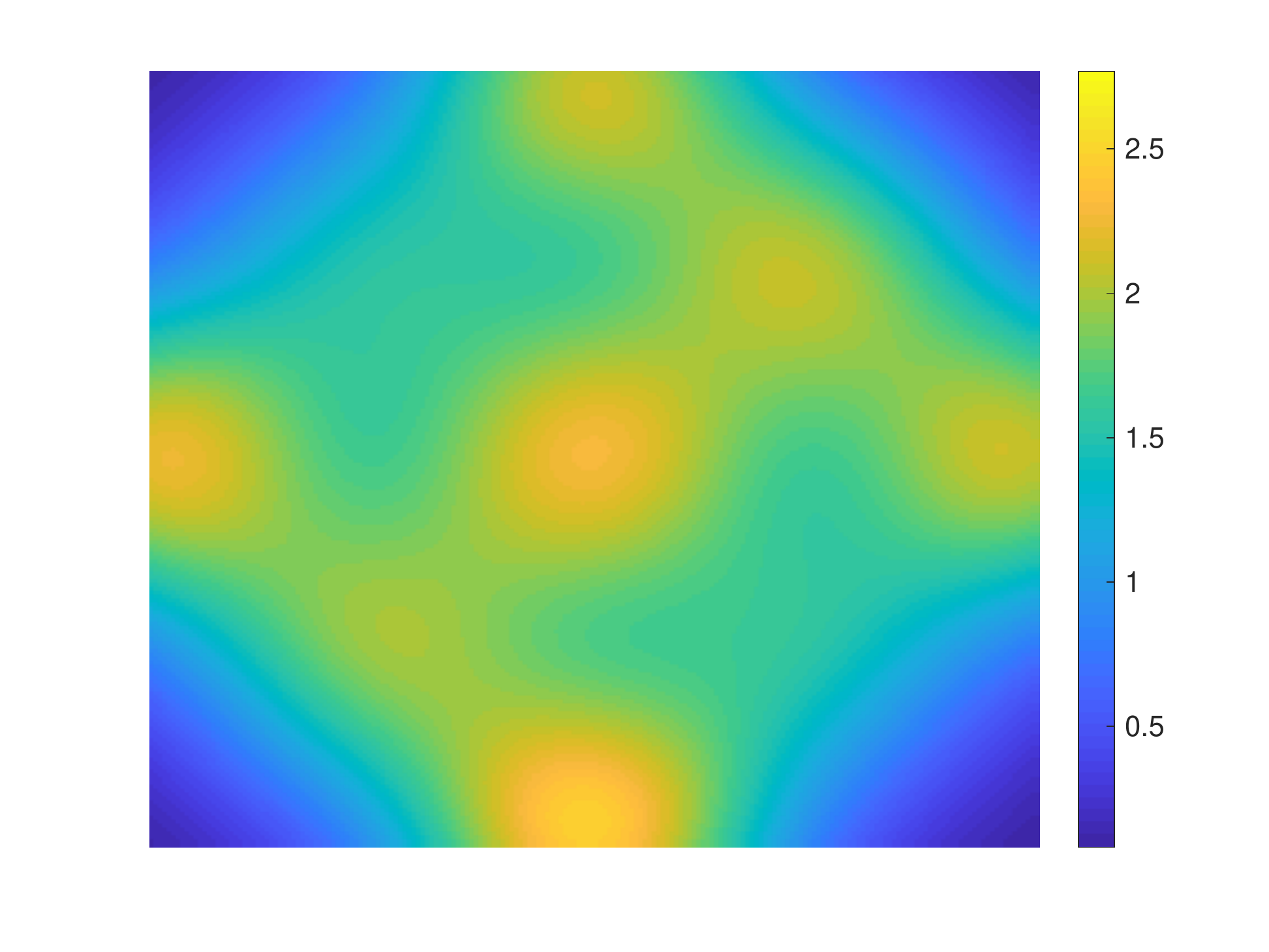}
     \put (23,-3) {\scriptsize ADNN, $tol=0.05$}
  \end{overpic}
\end{center}
\caption{Example 1:  Posterior mean (top) and   posterior standard deviation (bottom) arising from direct MCMC, ADNN ($tol=0.1$) and  ADNN ($tol=0.05$), respectively. The relative noise  level $\delta$ is $0.05$.}\label{adap_sol_eg1-2}
  \end{figure}

To assess the sampling accuracy of the ADNN algorithm, Fig. \ref{pos_contour_eg1} provides the marginal distributions of each component of the parameters, and the contours of the marginal distributions of each pair of components. The black lines represent the results generated by the direct MCMC approach based on the high-fidelity model evaluations (the reference solution), the red and blue lines represent results of the ADNN with $tol=0.1$  and $tol=0.05$, respectively.   The plots in Fig. \ref{pos_contour_eg1} suggest that the ADNN algorithm results in a good approximation to the reference solution.  \textcolor{black}{The posterior mean and posterior standard deviation obtained by the ADNN approach and the direct MCMC approach with $\delta=0.01$ and $0.05$ are shown in Fig. \ref{adap_sol_eg1-1} and Fig. \ref{adap_sol_eg1-2}, respectively.}  We observe that all algorithms produce similar estimates of mean in this test case. Furthermore, the results presented in Figs. \ref{adap_sol_eg1-1} and \ref{adap_sol_eg1-2}  show clearly that our method also produces comparable accuracy as Ref. \cite{Yan+Zhou19JCP}.

   \begin{figure}
\begin{center}
  \begin{overpic}[width=0.32\textwidth,trim=35 10 45 5, clip=true,tics=10]{figure/u_D_s05_eg1-eps-converted-to.pdf}
  \put (43,-3) {\scriptsize Direct}
  \end{overpic}
    \begin{overpic}[width=0.32\textwidth,trim= 35 10 45 5, clip=true,tics=10]{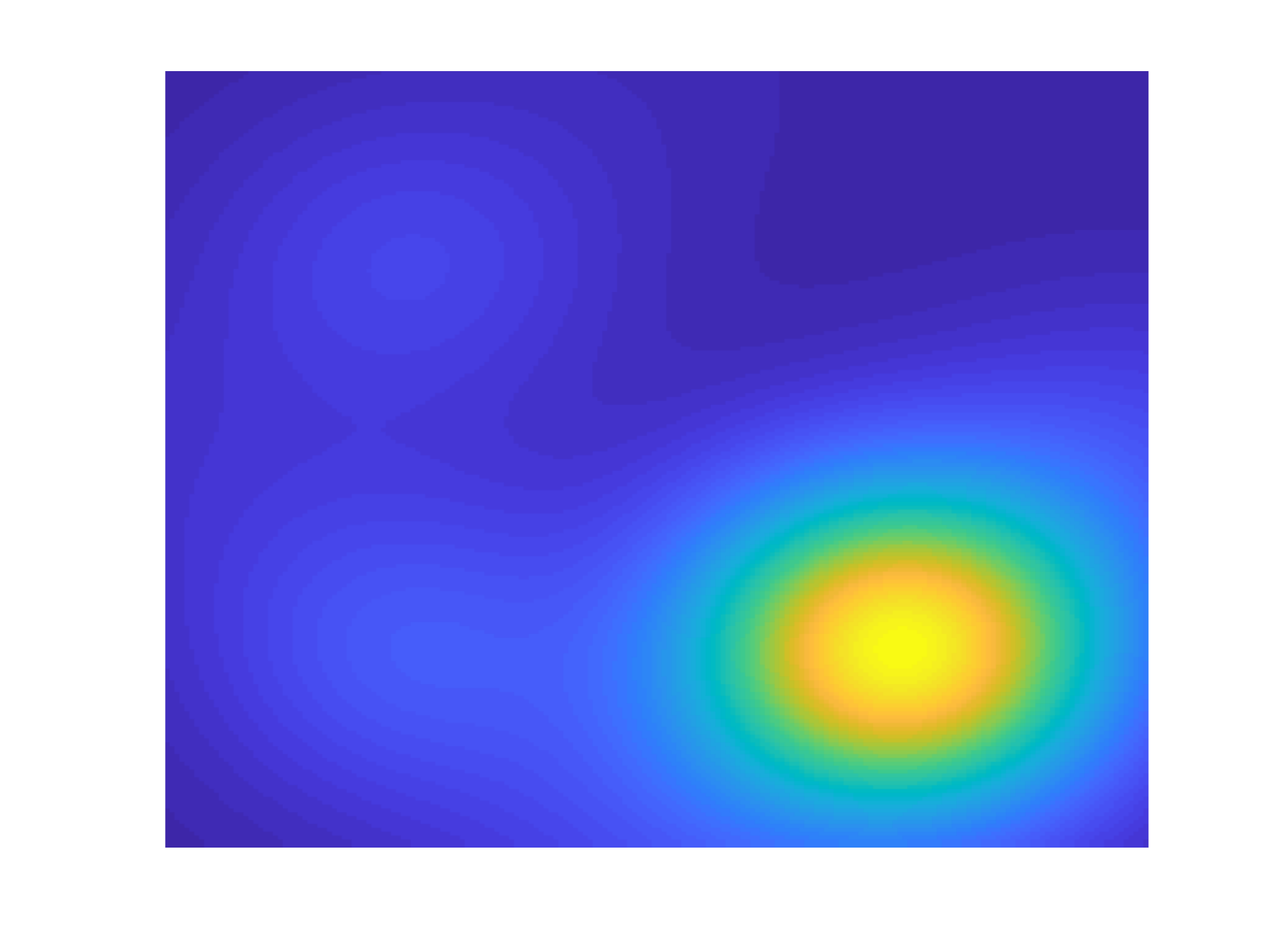}
    \put (30,-3) {\scriptsize DNN, $N=100$}
  \end{overpic}
    \begin{overpic}[width=0.32\textwidth,trim= 35 10 45 5, clip=true,tics=10]{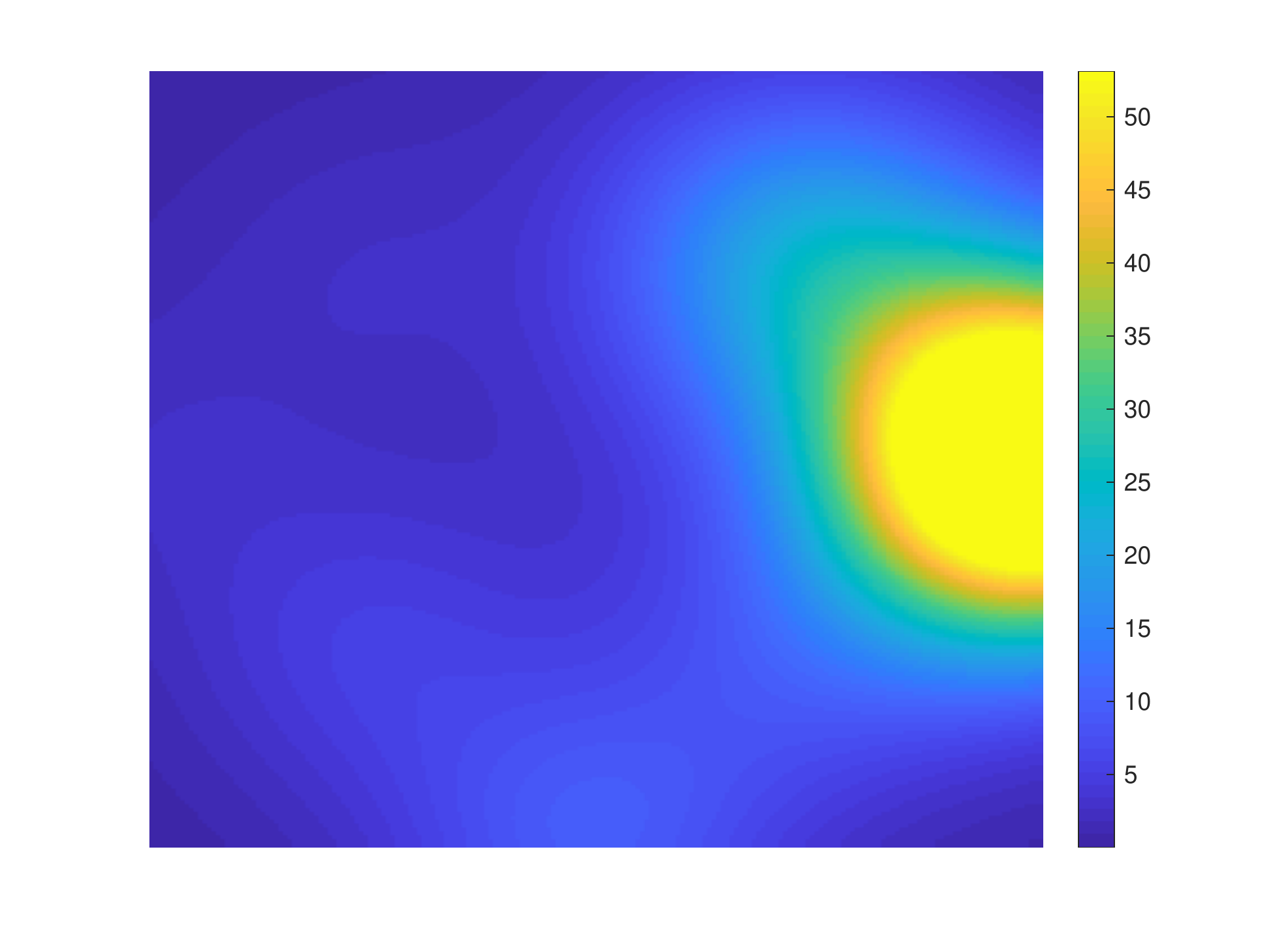}
    \put (30,-3) {\scriptsize DNN, $N=50$}
  \end{overpic}
\end{center}
\caption{Example 1:  Posterior mean arising from direct MCMC, prior-based DNN approach with $N=50$ and $N=110$, respectively.  The relative noise  level $\delta$ is $0.01$.}\label{pos_DNN_eg1}
\end{figure}

To verify the accuracy of our proposed algorithm, we also compute the posterior mean obtained by the prior-based DNN-approach described in Section \ref{sec:DNN}.  In this case,  the low-fidelity $\mathcal{NN}^L$ is built in advance and kept unchanged during MCMC computations. The posterior mean arising from the conventional MCMC and the DNN model with different sizes of the training dataset $N=\{50, 110\}$ are plotted in Fig. \ref{pos_DNN_eg1}.   Since the exact parameter is far from what is assumed in the prior, it is evident from the figure that the results using the prior-based DNN approach give a large error.  By comparing Figs. \ref{adap_sol_eg1-1} and \ref{pos_DNN_eg1}, we can obtain that the approximation results using ADNN are much more accurate than that of the prior-based DNN approach.

The computational costs, given by three mentioned algorithms, are presented in Table \ref{eg1_time}.  Moreover, we have also provide the results obtained by the adaptive multi-fidelity PC (AMPC) \cite{Yan+Zhou19JCP}, which is noted to be similar as the results from the ADNN approach. The main computational time in the DNN-based algorithm is \textcolor{black}{the offline high-fidelity model evaluations}. Upon obtaining  a trained DNN model $\mathcal{NN}^L$, the online simulation is very cheap as it does not require any high-fidelity model evaluations. For the ADNN, we do need the online forward model simulations to construct and refine the multi-fidelity DNN model $\mathcal{NN}^H$. Nevertheless, in contrast to $5\times 10^4$ high-fidelity model evaluations in the conventional MCMC, the number of \textcolor{black}{online high-fidelity model evaluations} for the ADNN with $N=\{50, 110\}$ are  300 and 180, respectively. As can be seen from Table \ref{eg1_time}, the ADNN algorithm  can significantly improve the approximation accuracy, yet without a dramatic increase in the computational time.  The results using the AMPC is displayed in the last column of Table \ref{eg1_time}. The number of high-fidelity model evaluations for the AMPC with an $N_C$-order ($N_C=2$) correction expansion are  1010, which is much larger than the ADNN.

\begin{table}[tp]
      \caption{Example 1. Computational times, in seconds, given by four different methods. }\label{eg1_time}
  \centering
  \fontsize{6}{12}\selectfont
  \begin{threeparttable}
    \begin{tabular}{ c ccccc}
  \toprule

   & \multicolumn{2}{c}{Offline}&\multicolumn{2}{c}{Online}\cr
\cmidrule(lr){2-3} \cmidrule(lr){4-5}

  \multirow{1}{*}{Method}  &$\text{$\#$ of \textcolor{black}{high-fidelity model evaluations}}$&CPU(s) &$\text{$\#$ of \textcolor{black}{ high-fidelity model evaluations}}$&CPU(s)     &\multirow{1}{*}{Total time(s)}\cr
  \midrule
    Direct                                     & $-$       & $-$         & 5$\times 10^4$    &$\text{\bf{1492.2}}$      & $\text{\bf{1492.2}}$   \cr
   DNN, $N=50$                 & 50    & 11.3     & $-$                       & 5.1                      &16.4     \cr
   DNN, $N=110$                & 110  & 12.7     & $-$                       &5.1                     & 17.8  \cr
   ADNN, $N=50$                & 50   & 11.3      & $300$                        & 31.9                         & 43.2  \cr
   ADNN, $N=110$             & $\text{\bf{110}}$ & $\text{\bf{12.7}}$ & $\text{\bf{180}}$    & $\text{\bf{18.4}}$       & $\text{\bf{31.1}}$   \cr
   AMPC, $N_C=2$       & $\text{\bf{110}}$ & $\text{\bf{2.9}}$ & $\text{\bf{1,010}}$                   & $\text{\bf{38.7}}$       & $\text{\bf{41.6}}$   \cr
    \bottomrule
      \end{tabular}
    \end{threeparttable}

\end{table}

 \subsubsection{The influence of tuning parameters $tol$ and $R$}

Intuitively, one would expect that the accuracy of the ADNN will improve as the value of the threshold $tol$ decreases. To verify this proposition,  we test several constant values choosing from  $tol \in \{0.1, 0.05\}$ and $N \in \{50, 110\}$. With these settings, we  run Algorithm \ref{alg:AMNN} using the ADNN model.   After discarding $2\times 10^4$ burn-in samples for each chain, we consider the evolution of the error as the chain lengthens; we compute an error measure at each step $rel(k)$   defined as
\begin{eqnarray*}
rel(k)= \frac{\|\bar{\kappa}-\kappa^{\dag}\|_{\infty}}{\|\kappa^{\dag}\|_{\infty}},
\end{eqnarray*}
where $\kappa^{\dag}$ are the ``true" posterior mean arising from direct MCMC,  and $\bar{\kappa}$ is the posterior  mean arising from ADNN.   Because the SGD for ADNN algorithm is random, we report the mean error over a size-10 ensemble of tests.

\begin{figure}
\begin{center}
  \begin{overpic}[width=0.5\textwidth,trim=15 0 20 15, clip=true,tics=10]{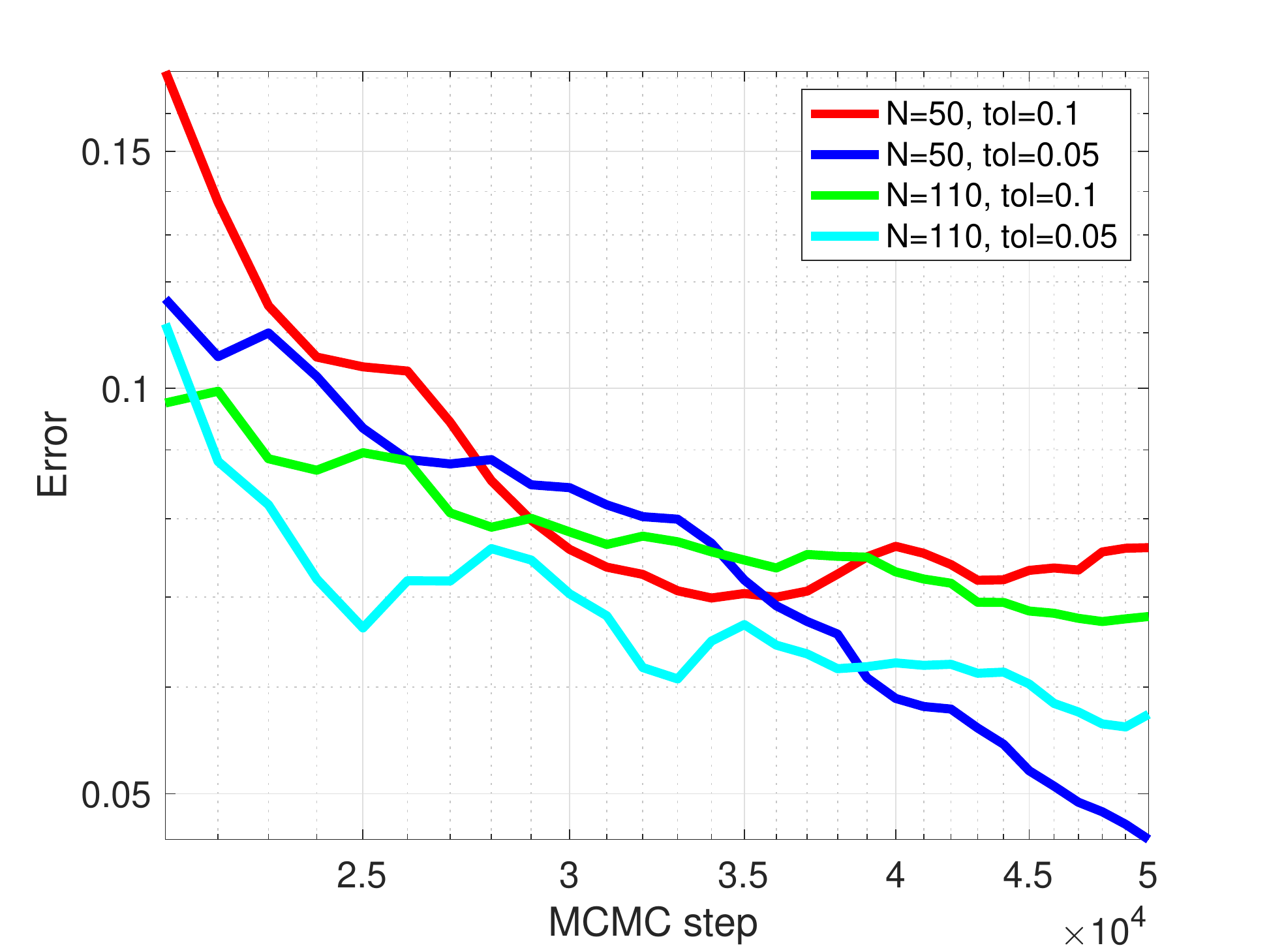}
  \end{overpic}
\end{center}
\caption{Example 1: The accuracy  errors $rel(\kappa)$ of ADNN using various numbers of the threshold $tol$.}\label{err_eg1_tol}
 \end{figure}

  \begin{figure}
\begin{center}
  \begin{overpic}[width=0.45\textwidth,trim=20 0 20 15, clip=true,tics=10]{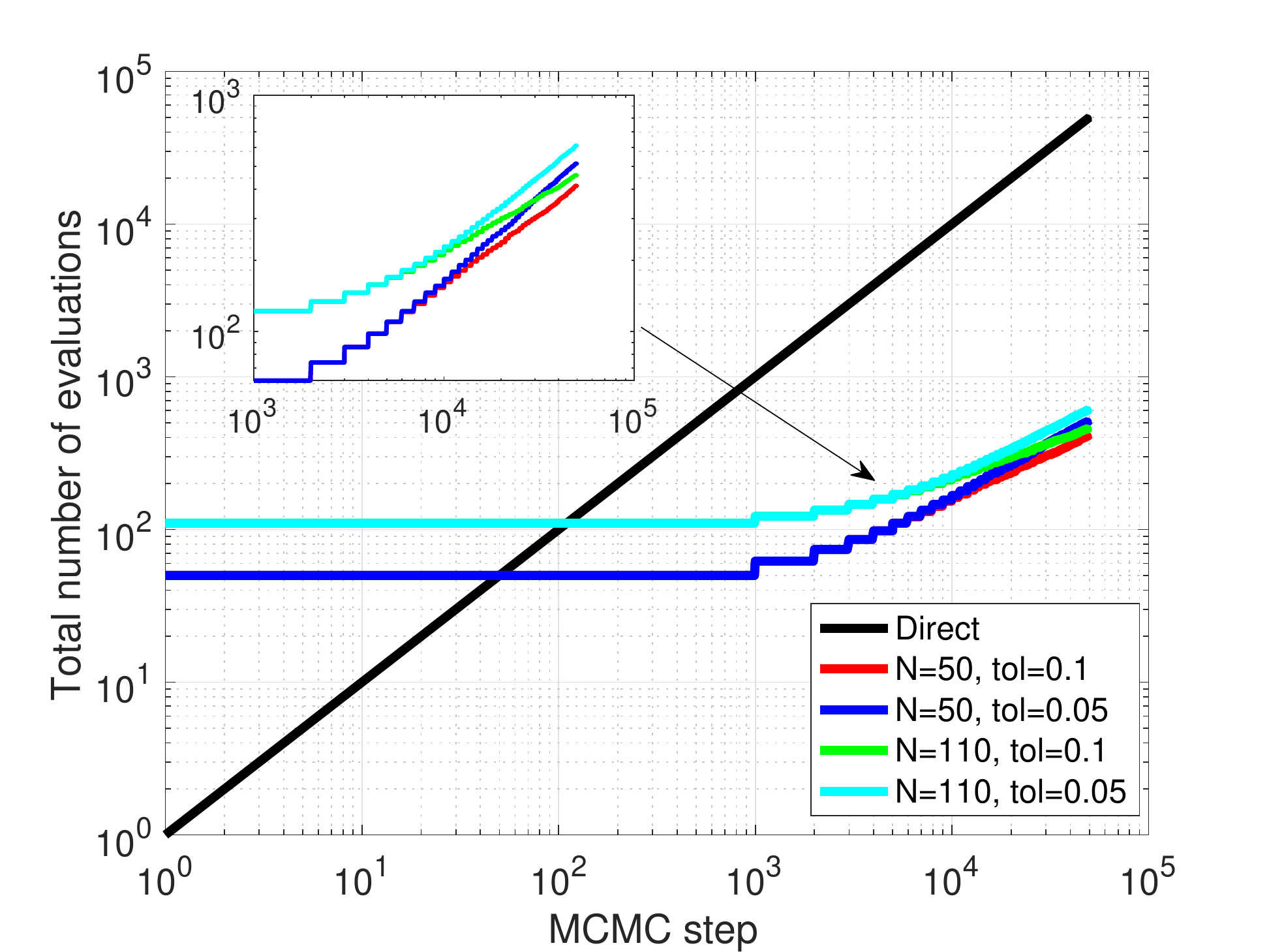}
  \end{overpic}
    \begin{overpic}[width=0.45\textwidth,trim= 20 0 20 15, clip=true,tics=10]{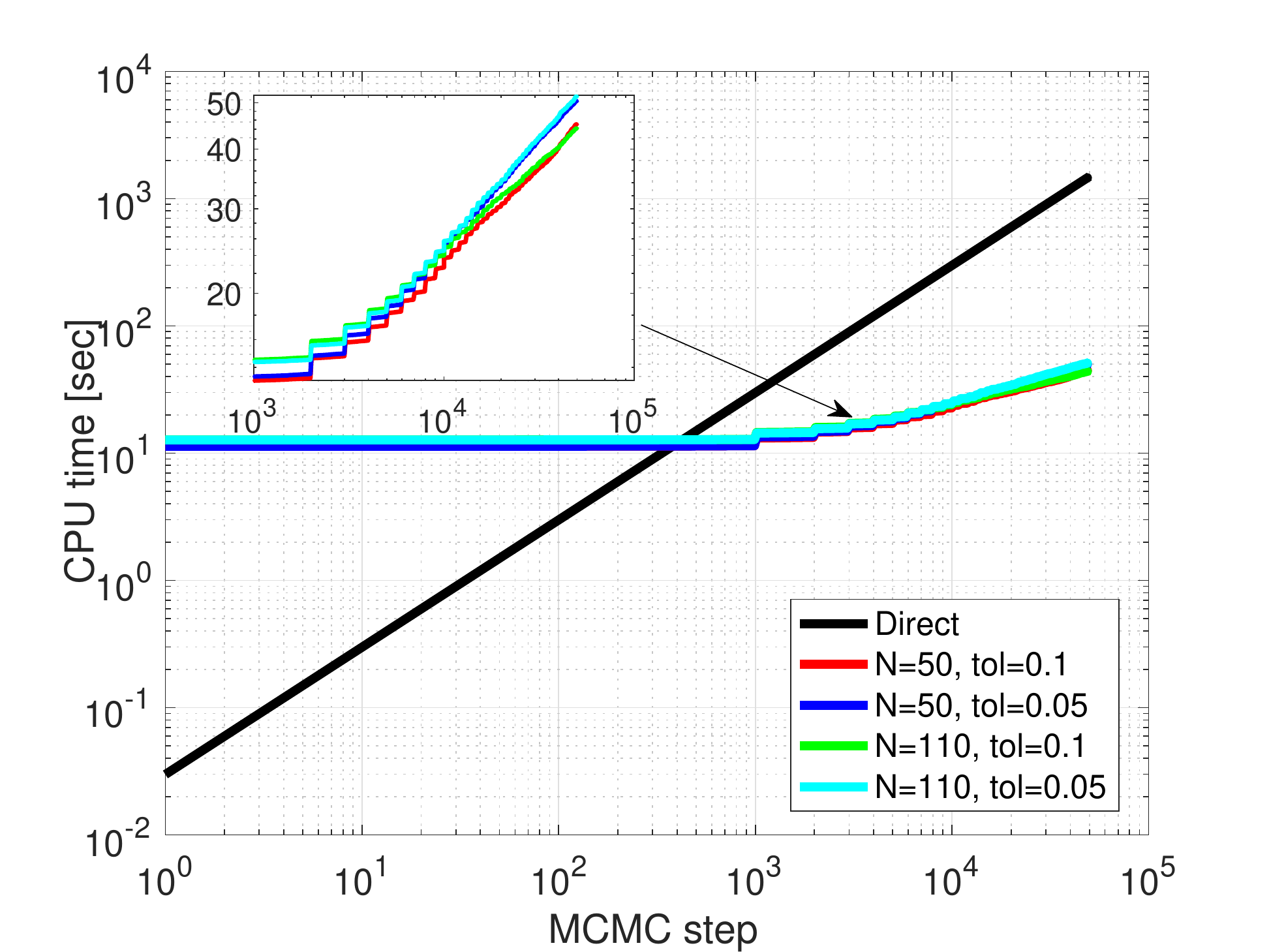}
  \end{overpic}
\end{center}
\caption{Example 1:  Total number of high-fidelity model evaluations and computational time for MCMC simulation; direct evaluation versus ADNN with various numbers of the threshold $tol$.}\label{cpu_eg1_tol}
  \end{figure}

    \begin{figure}
\begin{center}
\includegraphics[width=.45\textwidth]{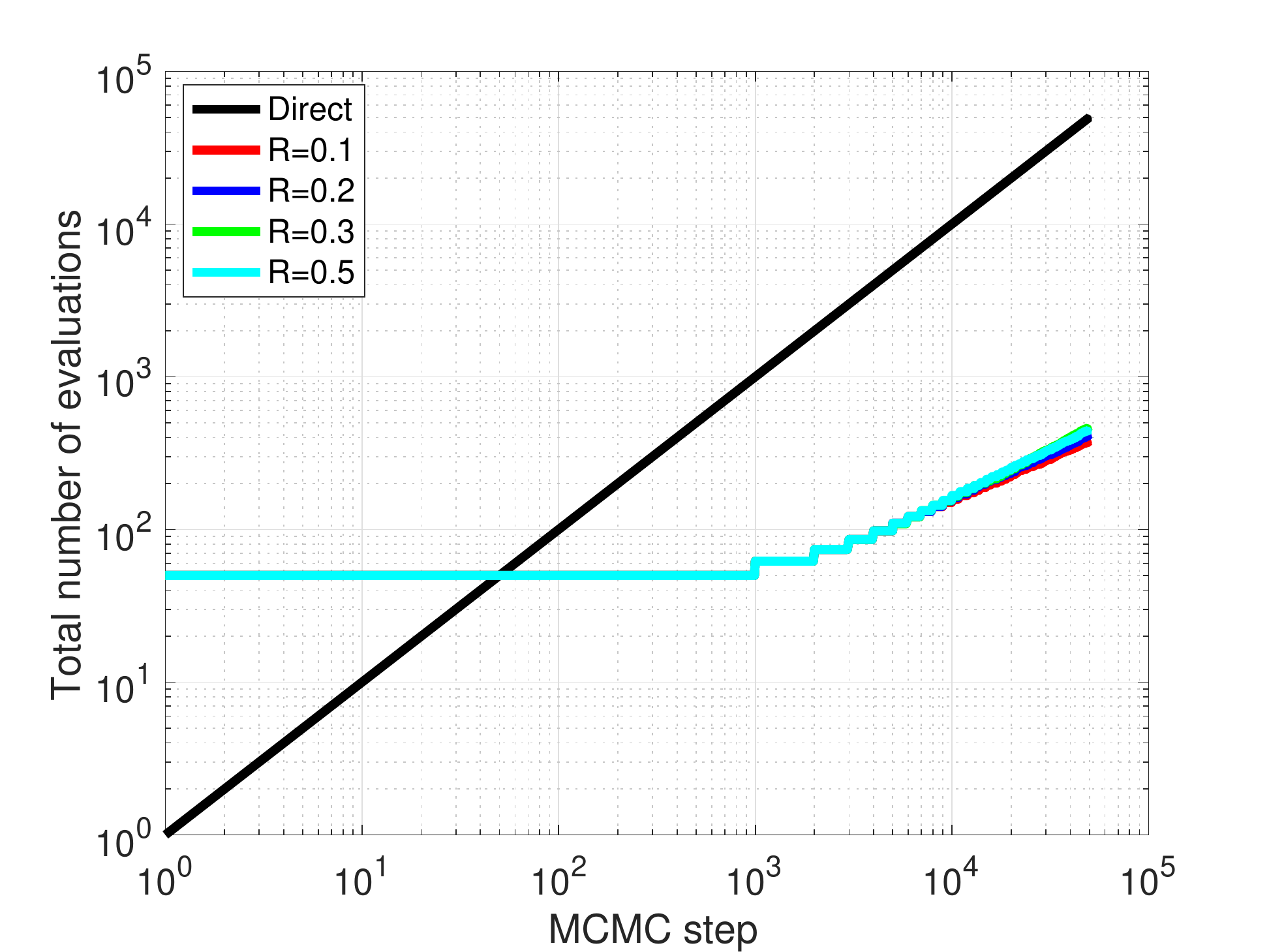}
\includegraphics[width=.45\textwidth]{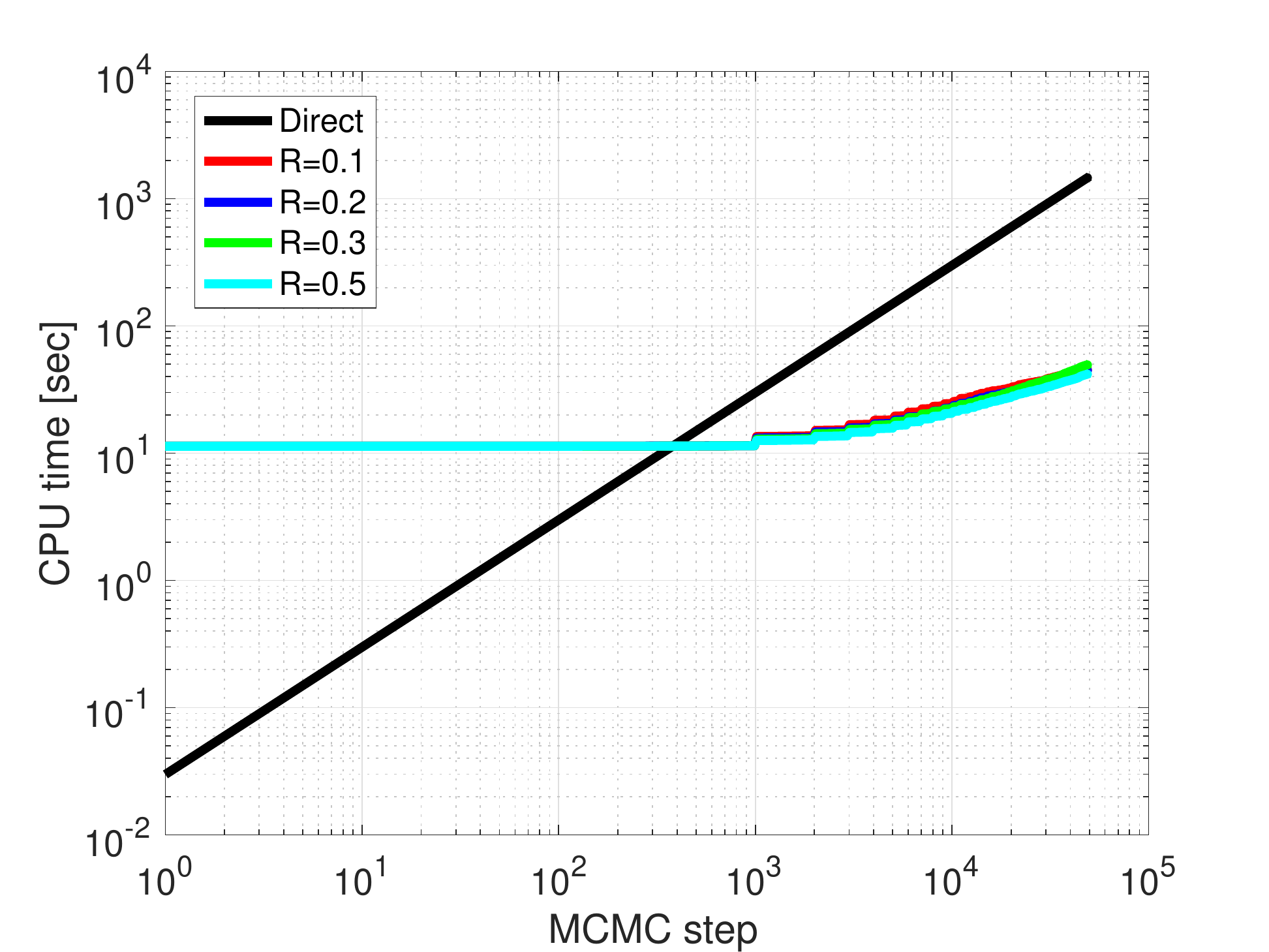}
\end{center}
\caption{Example 1: Total number of high-fidelity model evaluations and the computational time for MCMC simulation; direct evaluation versus ADNN with  various numbers of the radius $R$.}\label{cpu_eg1_R}
  \end{figure}

The accuracy comparison is given in Fig. \ref{err_eg1_tol}, which shows the evolution of the relative error with the number of MCMC steps. The corresponding computational costs are summarized in Fig.  \ref{cpu_eg1_tol}, which shows the total number of high-fidelity model evaluations and CPU times performed for any given number of MCMC steps.
As expected, the relative error decreases when threshold $tol$ is smaller; these values trigger more frequent refinements. When refinement is set to occur at a very low rate, the resulting chain is inexpensive, and in contrast, smaller values of $tol$ show increased cost and reduced errors.  However, even with $N=110, tol=0.05$, the speedup of the ADNN approach over the direct method is still quite dramatic.  The total computational  time required by the conventional MCMC grows linearly with the number of samples. On the other hand, the ADNN only requires a fixed amount of computational  time when training the prior-based DNN model.  For the model correction procedure, the per-sample cost is several orders of magnitude smaller than the direct evaluations, due to the fact that the updated models are highly correlated, and thus the online computation can be very efficient yet admits good approximation results.

Next, we investigate the sensitivity of the numerical results with respect to the radius $R$. The numerical results for Example 1, with $N=50, tol =0.1$ and various values for the radius $R$ are illustrated in Fig.  \ref{cpu_eg1_R}.  From the numerical results, we can conclude that the proposed scheme is relatively insensitive to the local radius $R$.

 \textcolor{black}{
\subsubsection{The influence of depth and width of $\mathcal{NN}^H$}
In this section, we investigate the influence of the depth \& width of $\mathcal{NN}^H$.  Since constructing a sufficiently accurate surrogate $\mathcal{NN}^L$ over the whole domain of the prior distribution is not necessary, and furthermore, even for a less accurate low-fidelity surrogate model, the ADNN approach can still obtain a good approximation to the reference solution due to the correlation of two models. Therefore, we focus only the influence of depth and width for $\mathcal{NN}^H$ due to the fact that the few high-fidelity data may yeild overfitting. Hence, we limit the ranges of the depth $L$ and width $d_k$ as $L\in \{1,2,3\}$ and $d_k \in \{10,20,30,50,70\}$, respectively. }

 \textcolor{black}{
In order to analyze the sensitivity of the proposed method with respect to the structures of the $\mathcal{NN}^H$, we first  train a prior-based DNN surrogate $\mathcal{NN}^L$ with 4 hidden layers and 40 neurons per layer using $N=50$ training points.  The relative error $rel(\kappa)$ (and the required number of the online high-fidelity model evaluations) of the ADNN approach for $tol=0.1, Q=10$ and $\delta =0.05$ for Example 1 are presented in Table \ref{rel_L}. As shown in this table, the computational results for ADNN with different depth $L\in\{1, 2\}$ and width $d_k\in\{20,30,50\}$ are almost the same.  In addition, the use of $L=3$ or $d_k=70$ yields a very inaccurate approximation to the reference solution, as we have only used very few high-fidelity data, i.e., $Q=10$, to train the multi-fidelity neural network $\mathcal{NN}^H$.  In order to improve this, one need to increase the number of train points $Q$. The corresponding results using various values of $Q$ are shown in Table \ref{rel_L}. As expected, the $rel(\kappa)$ decrease as the number of value $Q$ increase. However, the required online high-fidelity model evaluations also become larger.  To reduce the online computational cost which retain the accuracy of estimate results, a reasonable choice of the size of $\mathcal{NN}^H$ may be $L\in\{1,2\}$ and $d_k\in \{30,50\}$ in moderate dimensions.  These numerical results  also demonstrate the robustness of the ADNN approach. }

\begin{table}
  \centering
\begin{tabular}{|c|c|c|c|c|c|}
\hline
\diagbox{L}{$d_k$} & 10& 20 & 30 & 50 &70\\
\hline
1 & 0.5243 (550) & 0.0895 (570) & 0.0800 (390) & 0.0554 (370)  & 0.0817 (380) \\
\hline
2 & 0.1461 (430) & 0.0665 (410) & 0.0614 (380) & 0.0681 (340)  & 0.1396 (390)\\
\hline
3 & 0.5256 (570) & 0.3406 (570) & 0.1292 (440) & 0.0868 (480) & 0.1227 (470) \\
\hline
\end{tabular}
  \caption{Example 1:  The relative error $rel$ (the number of online high-fidelity model evaluations) obtained using  ADNN approach with $tol=0.1, \delta =0.05$ and $Q=10$.}\label{rel_L}
    \centering
\begin{tabular}{|c|c|c|c|c|c|}
\hline
\diagbox{Q}{$d_k$} & 10& 20 & 30 & 50& 70 \\
\hline
10 & 0.5256 (570) & 0.3406 (570) & 0.1292 (440) & 0.0868 (480)  & 0.1227 (470) \\
\hline
20 & 0.1421 (800) & 0.0932 (720) & 0.0781 (620) & 0.0743 (500)  & 0.0670 (580)\\
\hline
30 & 0.0653 (1210) & 0.0271 (880) & 0.0407 (910) & 0.0602 (640)  & 0.0586 (700) \\
\hline
\end{tabular}
  \caption{Example 1:  The relative error $rel$ (the number of online high-fidelity model evaluations) obtained using ADNN approach with $tol=0.1, \delta =0.05$ and $L=3$. }\label{rel_Q}
\end{table}

 \subsection{Example 2}

   \begin{figure}
\begin{center}
  \begin{overpic}[width=0.45\textwidth,trim=35 10 45 5, clip=true,tics=10]{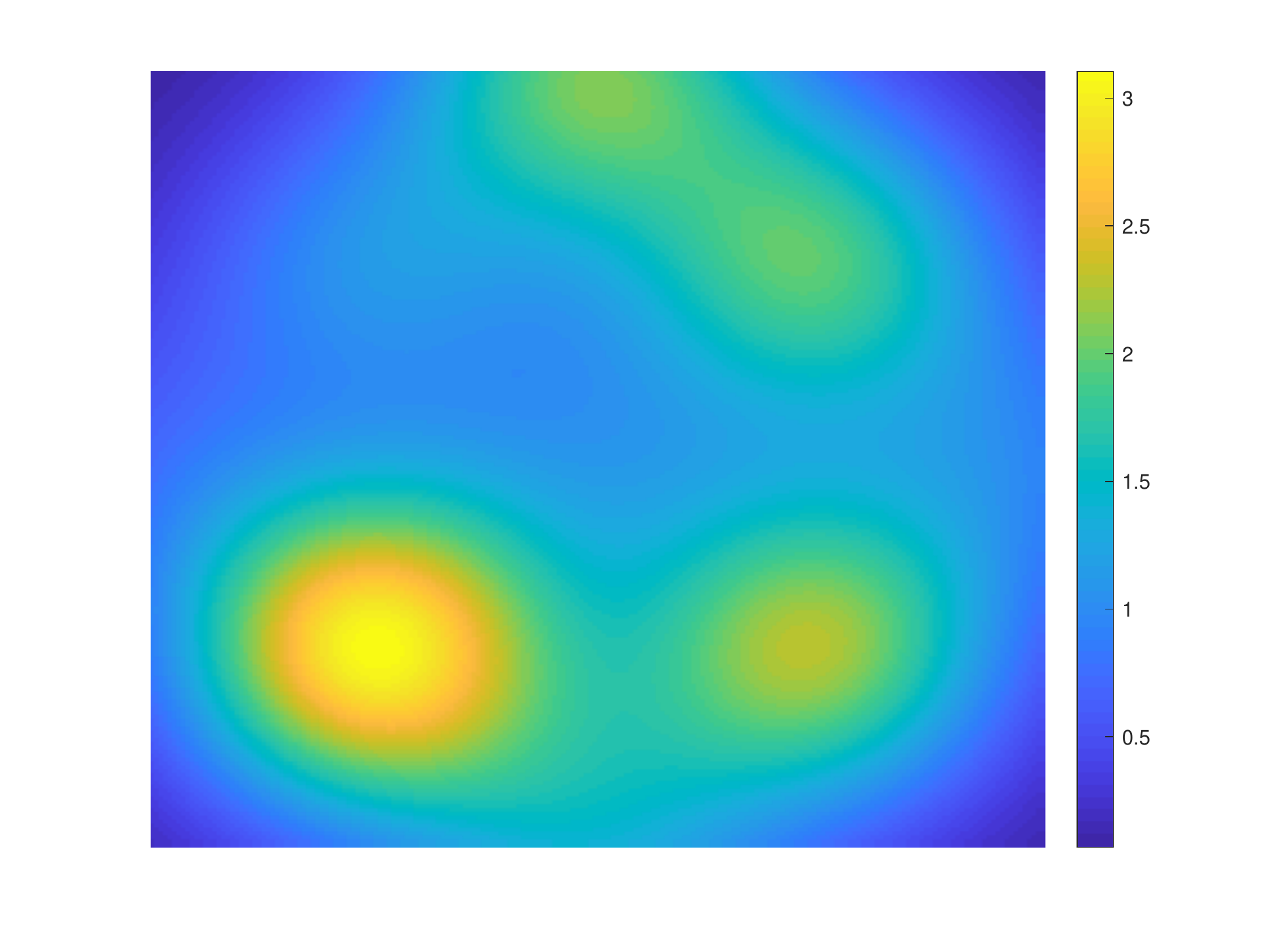}
    \put (35,-3) {\scriptsize exact solution}
  \end{overpic}
  \begin{overpic}[width=0.45\textwidth,trim=35 10 45 5, clip=true,tics=10]{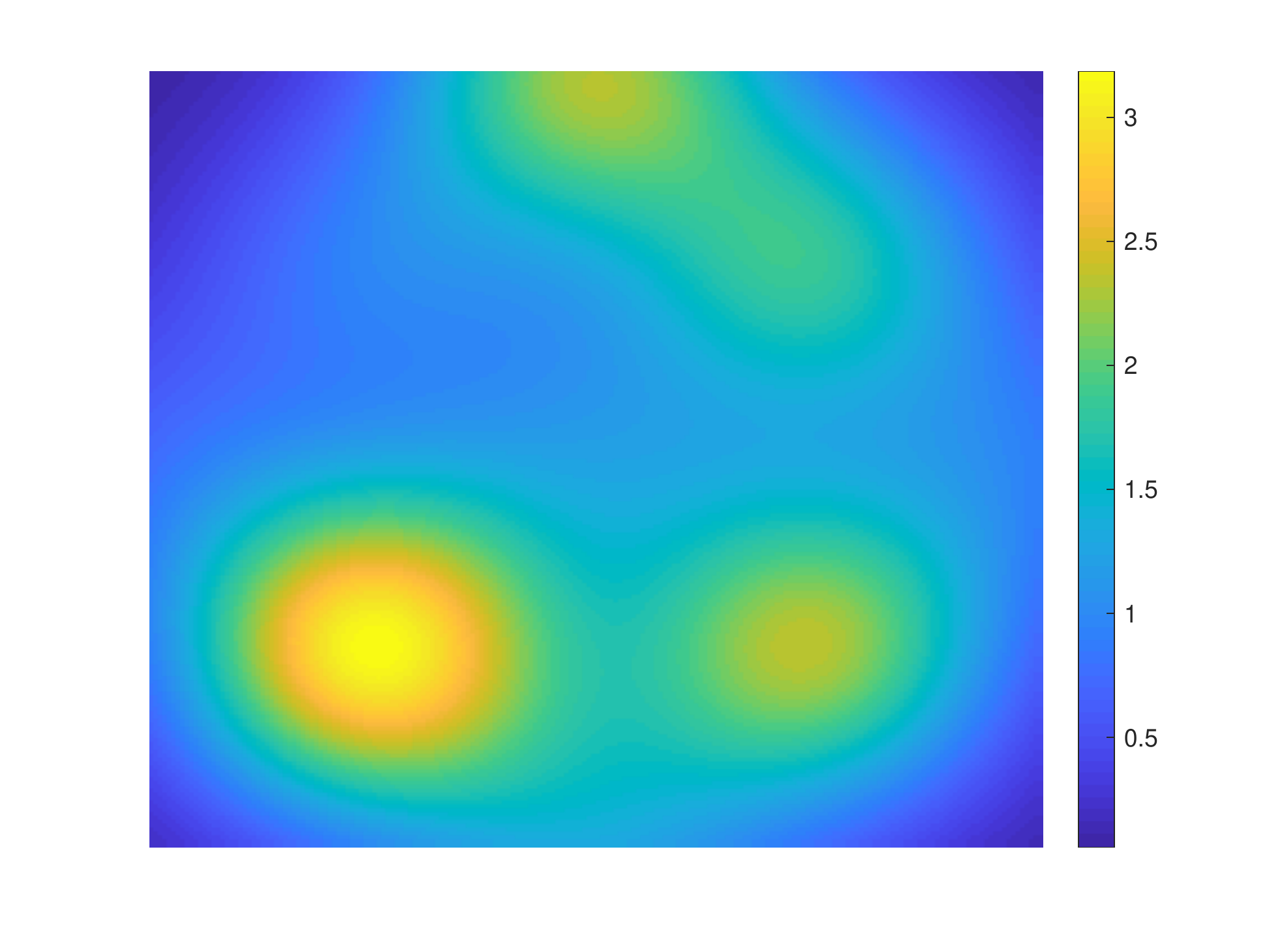}
    \put (35,-3) {\scriptsize reference solution}
  \end{overpic}
 \end{center}
\caption{Example 2.  Left: the true permeability used for generating the synthetic data set. Right: the reference solution arising the full model.}\label{exact_eg2}
  \end{figure}

As the second example,  the true parameter  is a draw from the prior distribution described in Example 1. In other words, we consider the best-case-scenario where our prior knowledge includes  the truth.  The exact permeability used for generating the synthetic data  and the reference solution arising the full model are displayed in Fig.\ref{exact_eg2}.

  \begin{figure}
\begin{center}
    \begin{overpic}[width=0.32\textwidth,trim= 35 10 45 5, clip=true,tics=10]{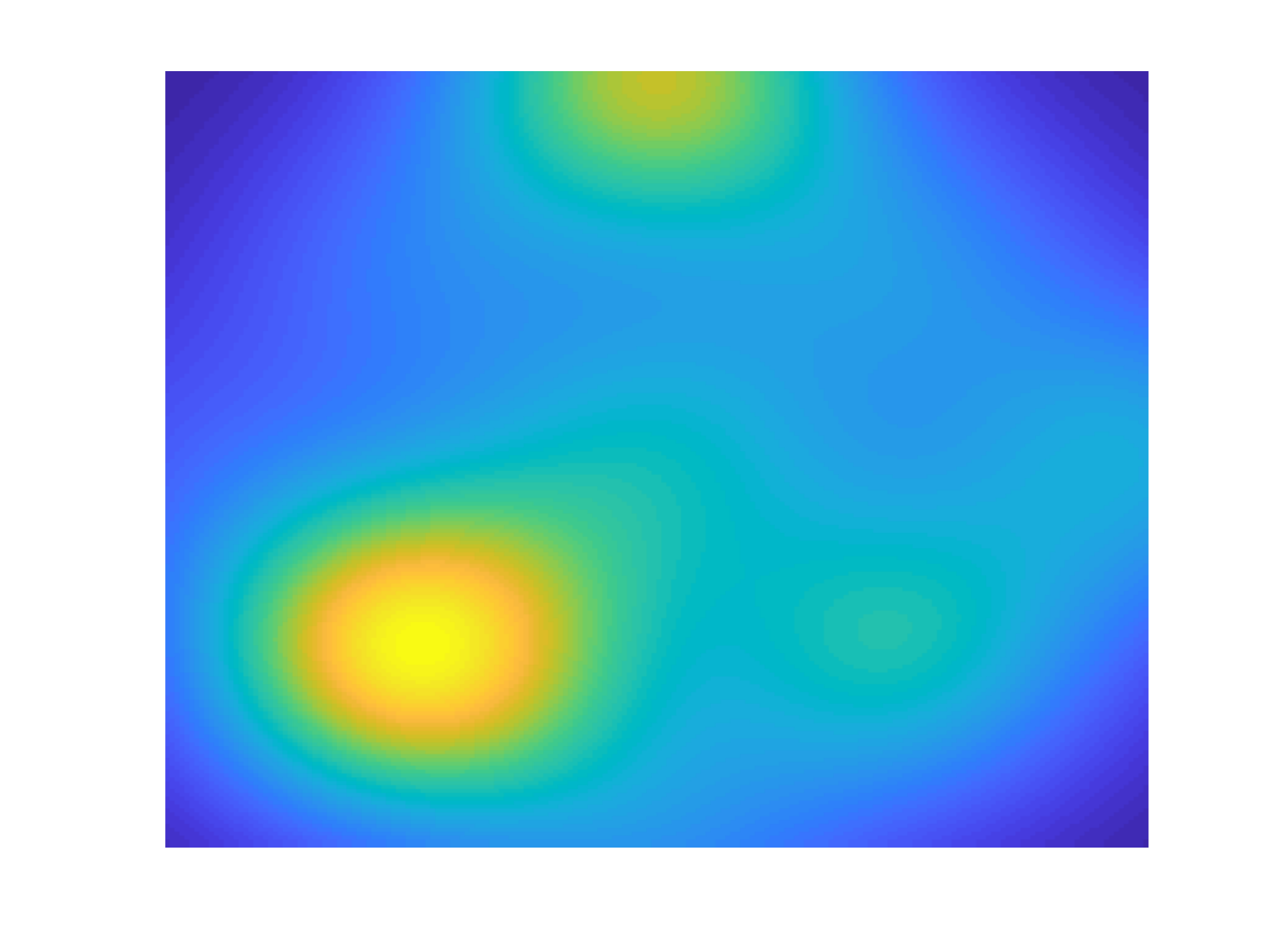}
  \end{overpic}
  \begin{overpic}[width=0.32\textwidth,trim=35 10 45 5, clip=true,tics=10]{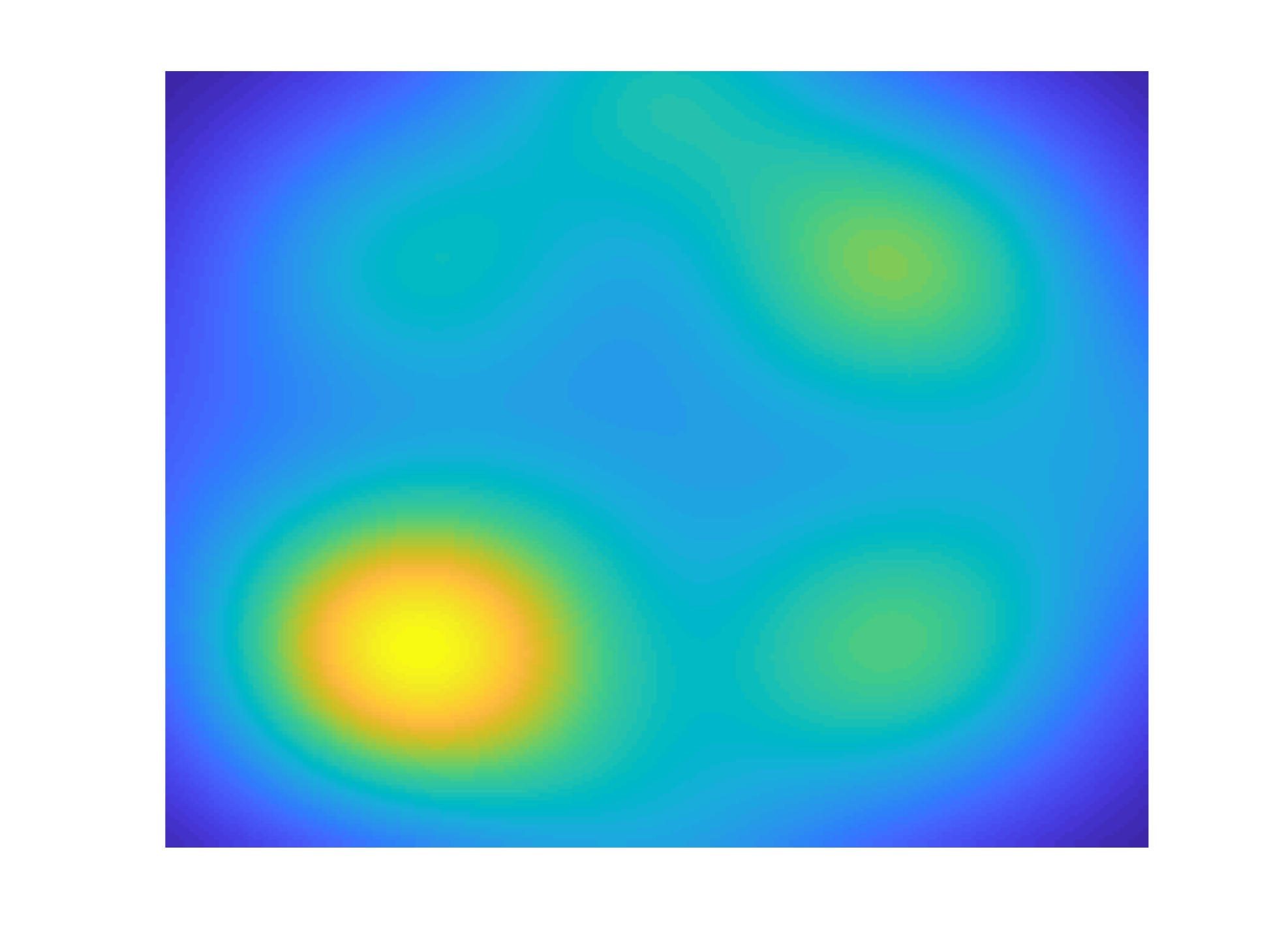}
  \end{overpic}
    \begin{overpic}[width=0.32\textwidth,trim=35 10 45 5, clip=true,tics=10]{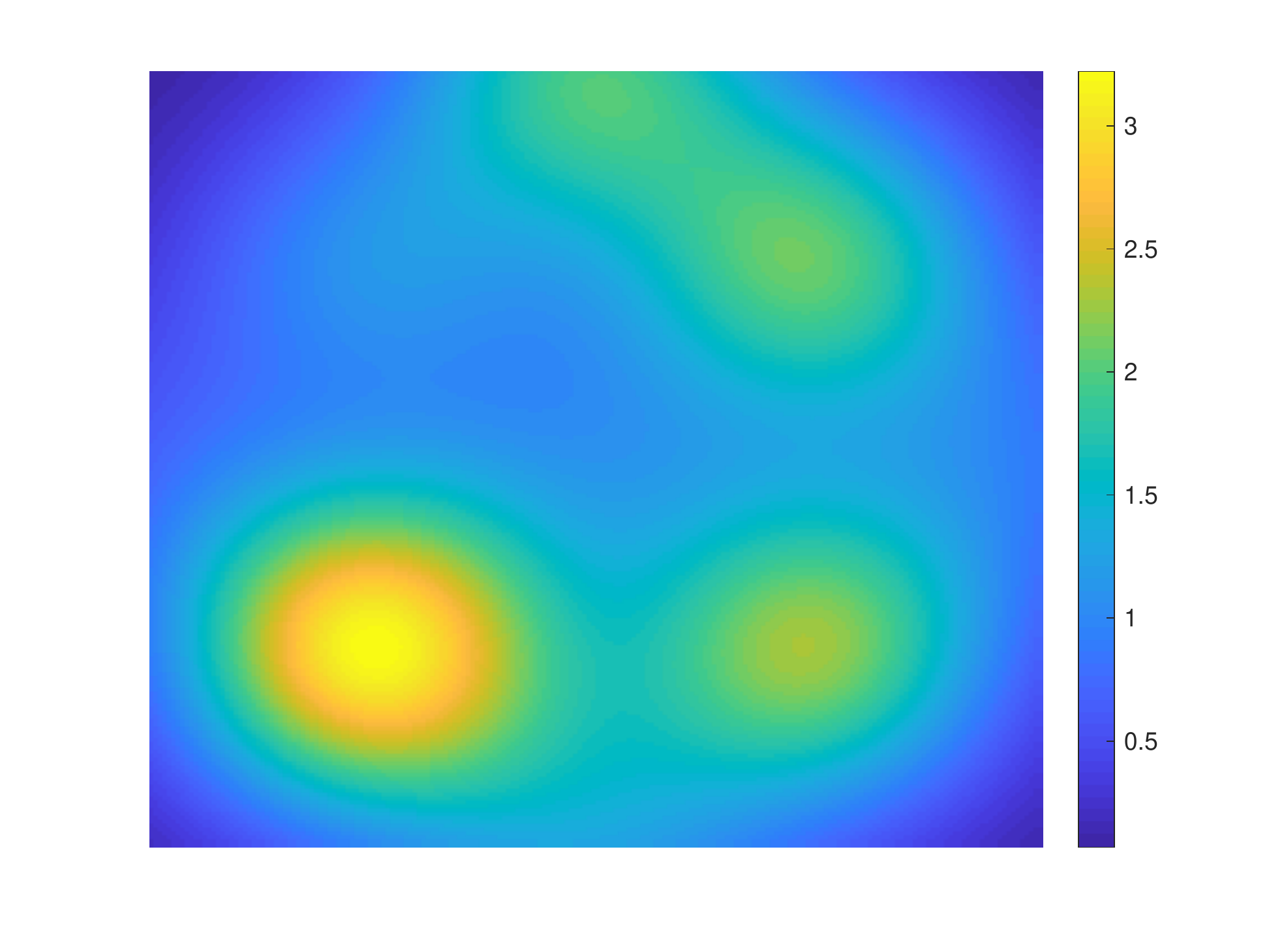}
  \end{overpic}
      \begin{overpic}[width=0.32\textwidth,trim= 35 10 45 5, clip=true,tics=10]{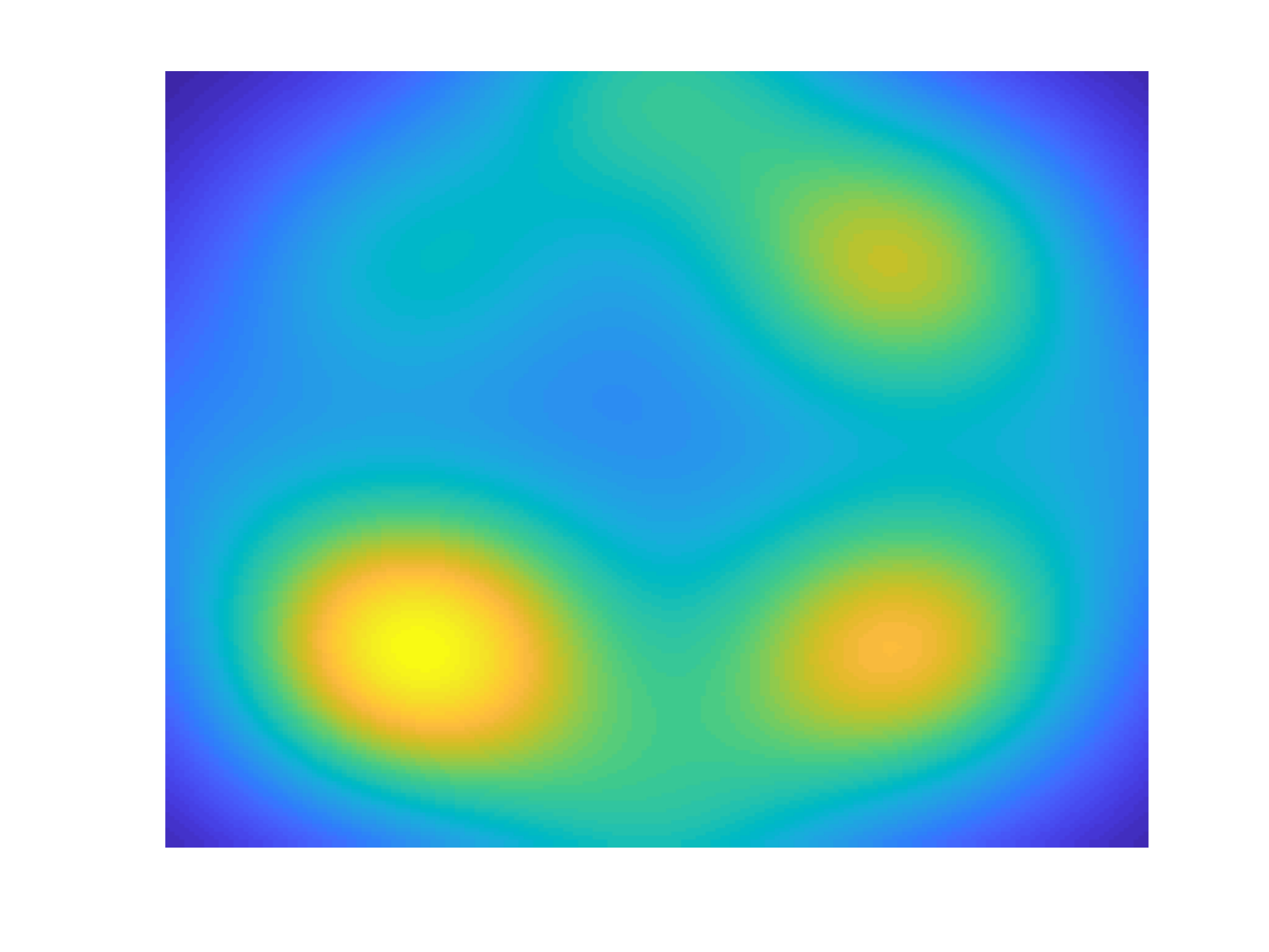}
     \put (40,-3) {\scriptsize DNN}
  \end{overpic}
    \begin{overpic}[width=0.32\textwidth,trim=35 10 45 5, clip=true,tics=10]{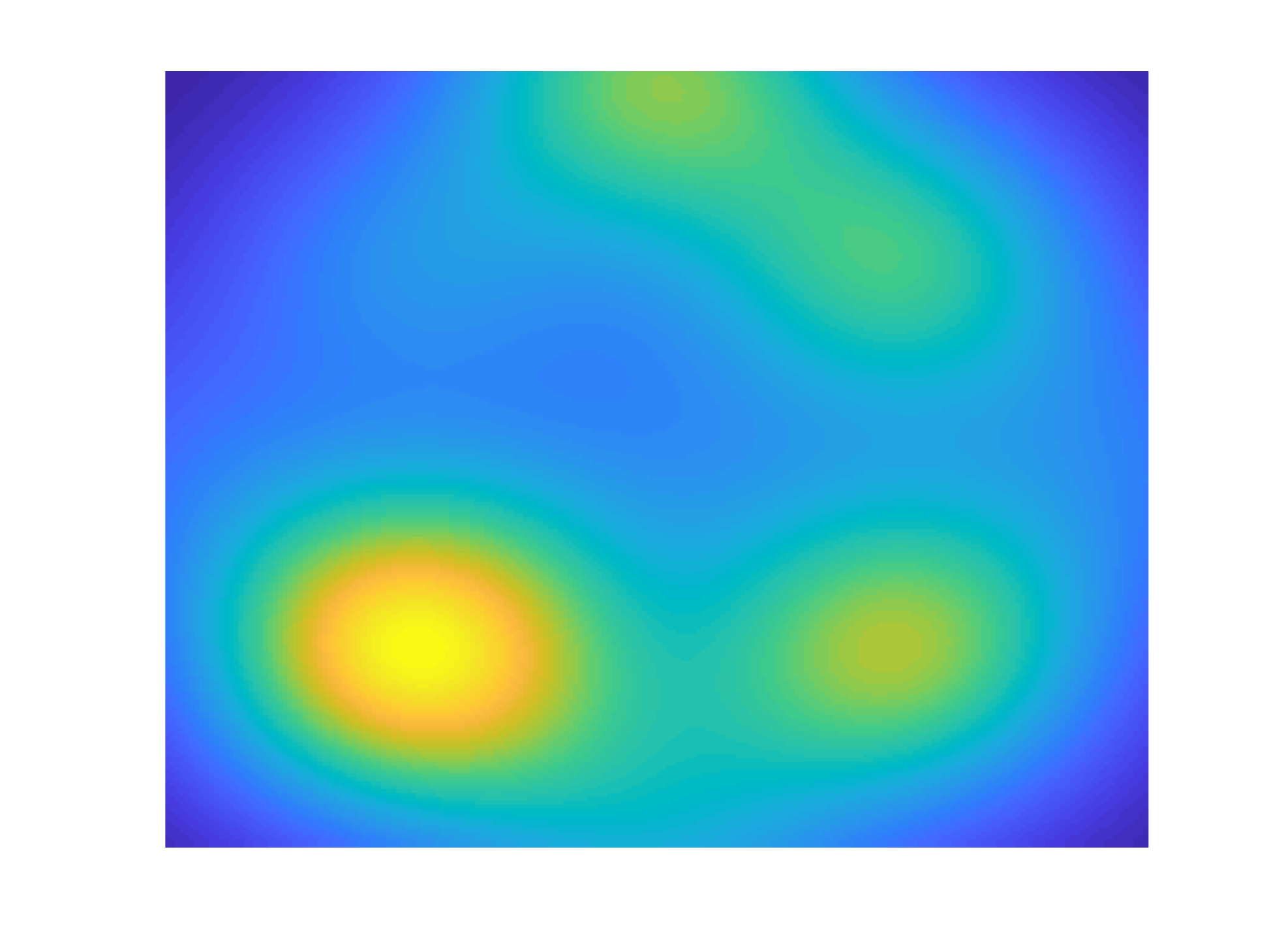}
      \put (28,-3) {\scriptsize ADNN $(tol=0.1)$}
  \end{overpic}
    \begin{overpic}[width=0.32\textwidth,trim=35 10 45 5, clip=true,tics=10]{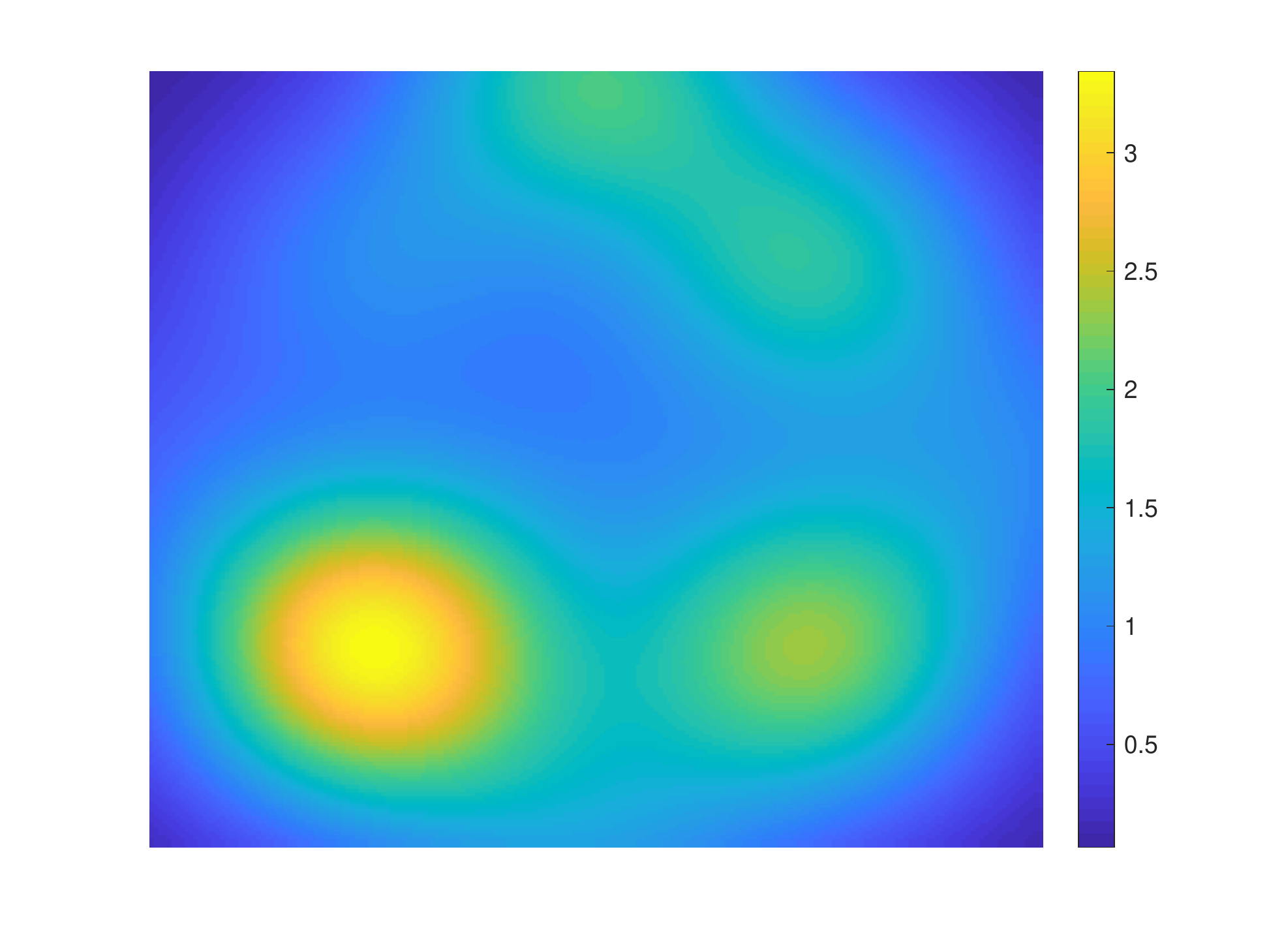}
      \put (23,-3) {\scriptsize  ADNN $(tol=0.05)$}
  \end{overpic}
\end{center}
\caption{Example 2:  Posterior mean arising from DNN approach and ADNN  with various value of $tol$. From top to bottom, the number of the training set $N$ is $50,110$ respectively.}\label{sol_eg2}
  \end{figure}

 \begin{figure}
\begin{center}
  \includegraphics[width=.6\textwidth]{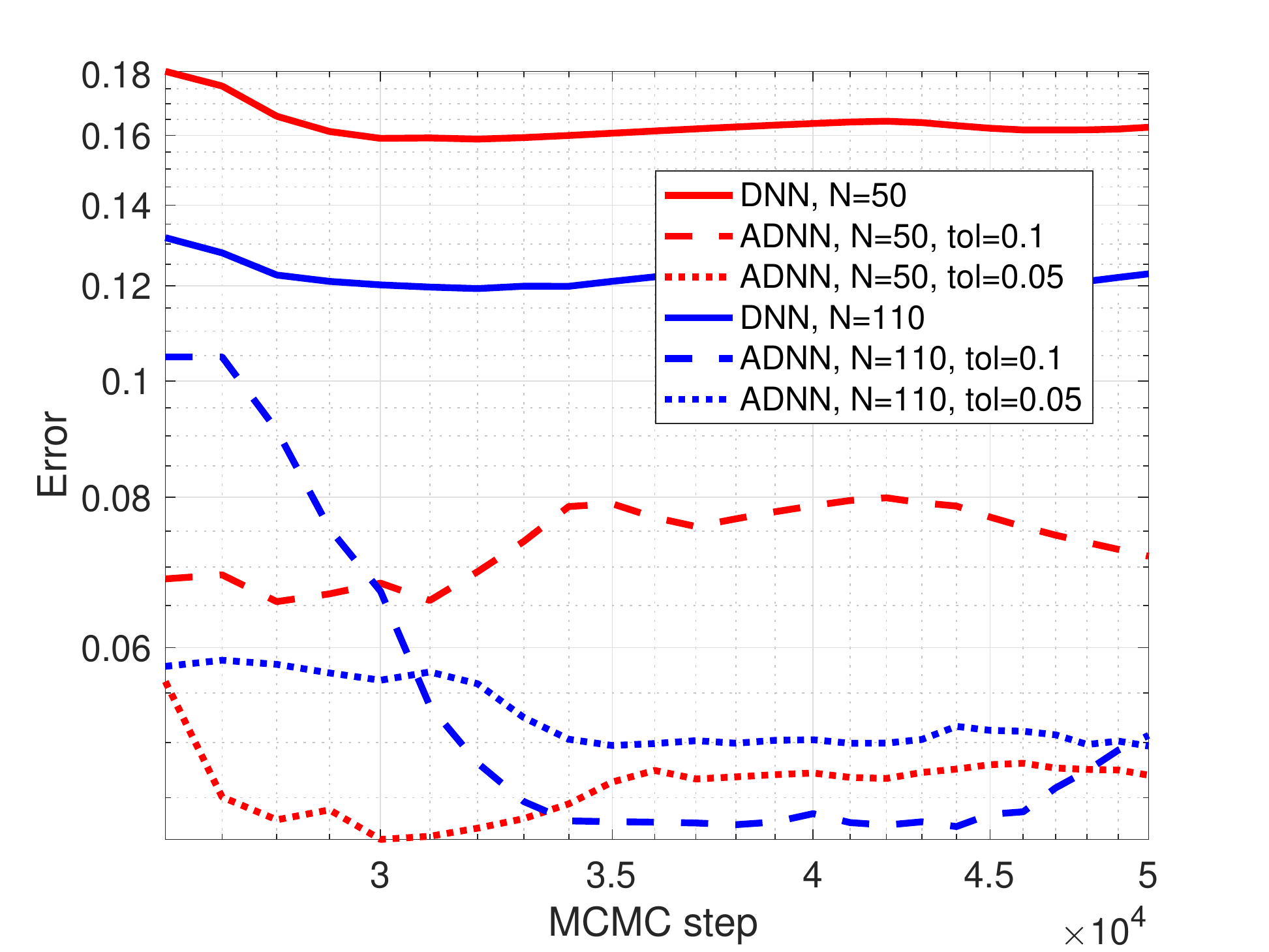}
\end{center}
\caption{Example 2: The accuracy errors $rel(\kappa)$ obtained using two different method.}\label{err_eg2}
\end{figure}

 \begin{figure}
\begin{center}
\includegraphics[width=.46\textwidth]{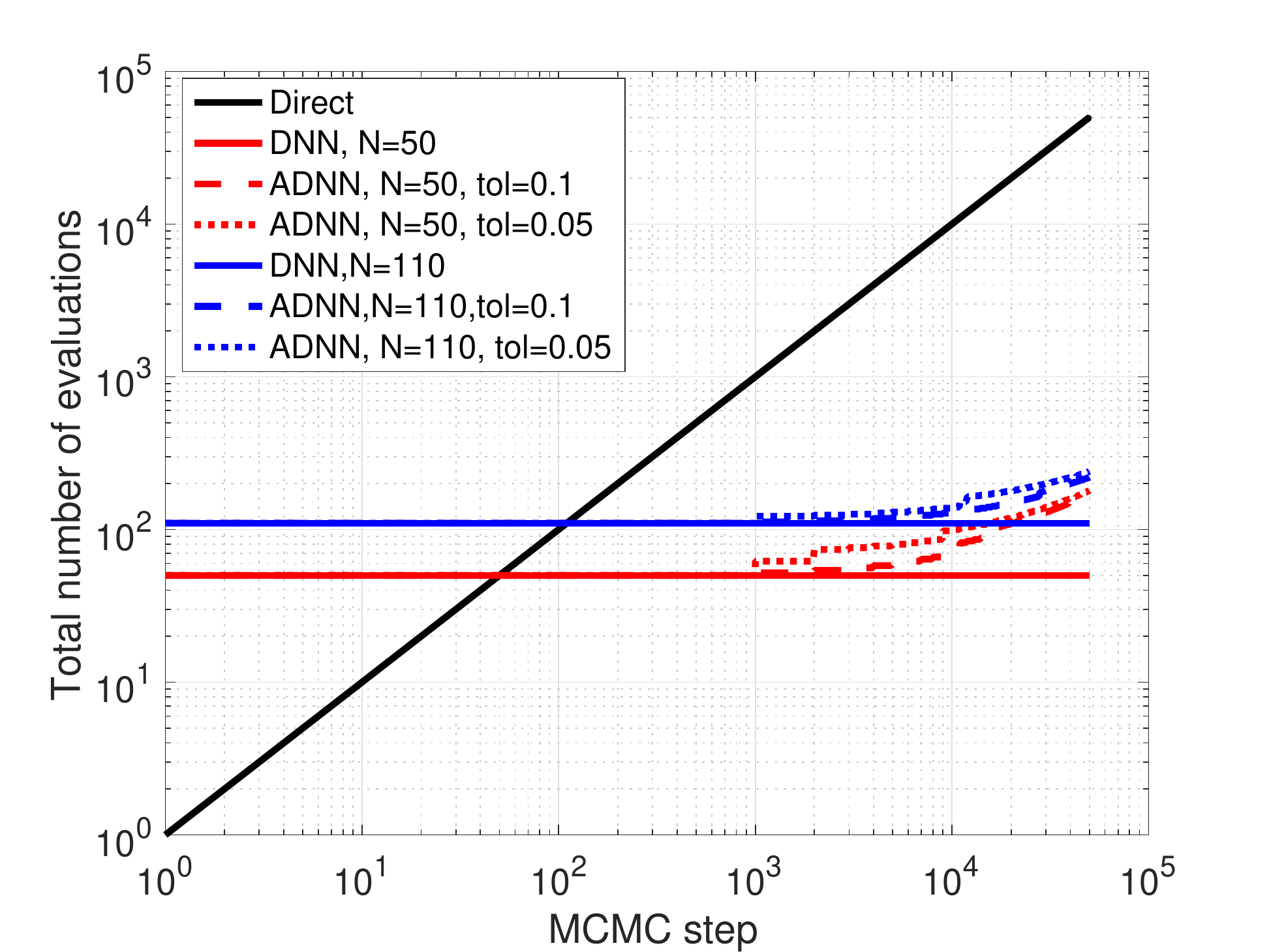}
\includegraphics[width=.46\textwidth]{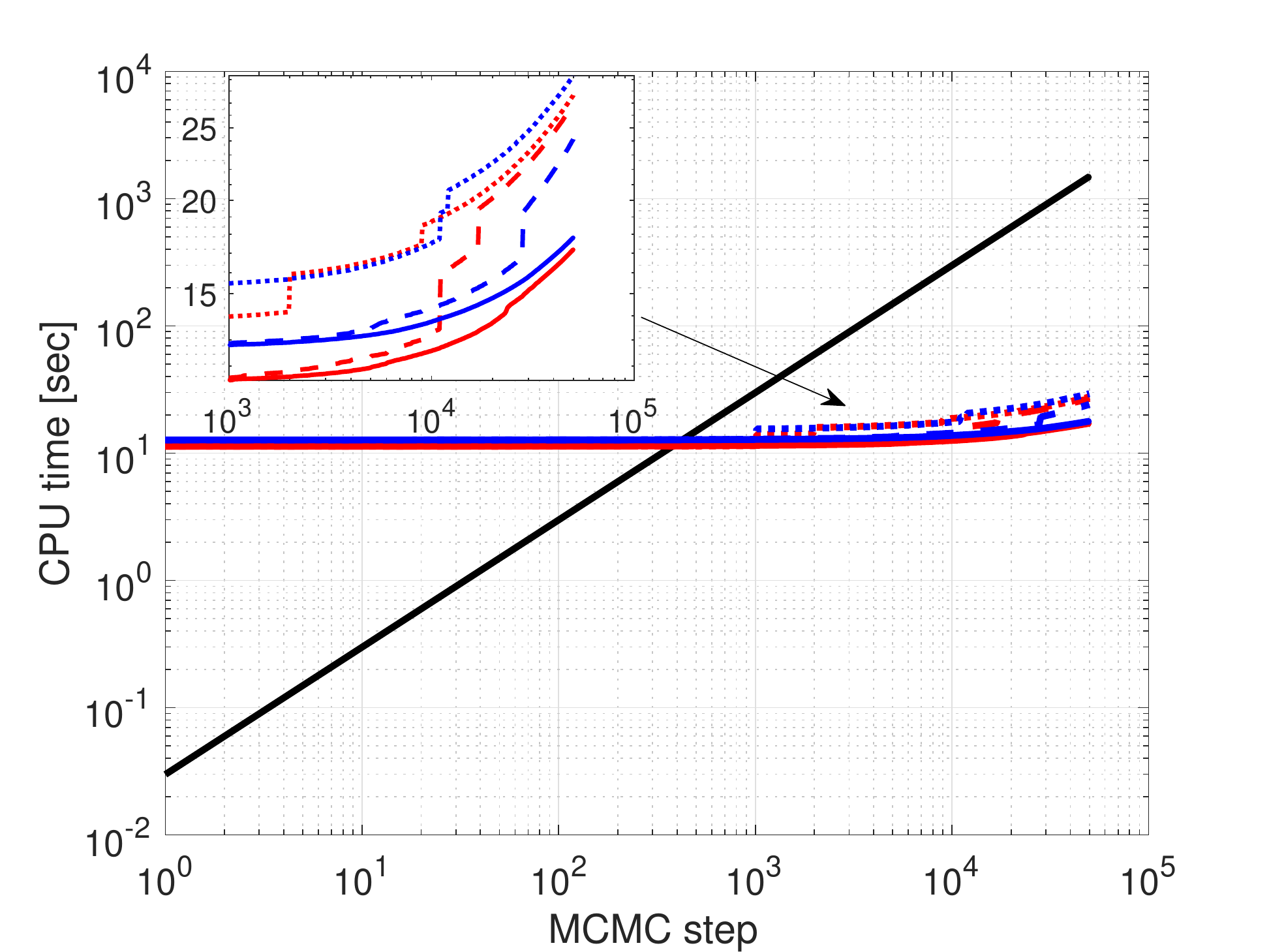}
\end{center}
  \caption{\small{Example 2: Total number of high-fidelity model evaluations and the computational time for MCMC simulation; direct evaluation versus prior-based DNN approach and ADNN.}}\label{cpu_eg2}
\end{figure}

Similar to the first example, we numerically investigate the efficiency of the ADNN approach.  Using the same setting as Example 1, we plot the posterior mean arising from  different methods using various values of  $N$.  The numerical  results obtained by DNN are shown in the left column of Fig.\ref{sol_eg2}.  \textcolor{black}{The corresponding  relative errors $rel(k)$ for the MCMC simulation are shown in Fig. \ref{err_eg2}.}  Compare with Fig. \ref{pos_DNN_eg1}, it can be seen that the numerical results obtained by DNN using $N=110$ training points agree with the reference solution. However, a smaller number of training points ($N=50$) results a poor estimate. The corresponding results obtained by ADNN are also shown in Figs.\ref{sol_eg2} and \ref{err_eg2}.  It is clearly shown that the ADNN approach results in a very good approximation to the reference solution. Even with a smaller  $N=50$ and a larger $tol=0.1$, the ADNN approach admits a rather accurate result.   \textcolor{black}{ The total number of high-fidelity model evaluations and the total computational time for MCMC simulation are summarized in Fig.\ref{cpu_eg2}. }  Again,  the online computational time required by ADNN and DNN is only a small fraction of that by the conventional MCMC.   It can also be seen from these figures that the ADNN offers a significant improvement in the accuracy, but does not significantly increase the computation time compared to the prior-based DNN approach.  This also confirms the efficiency of the ADNN algorithm for this best-case-scenario.

\subsection{Example 3: a high dimensional inverse problem}

  \begin{figure}
\begin{center}
  \begin{overpic}[width=0.45\textwidth,trim=35 10 45 5, clip=true,tics=10]{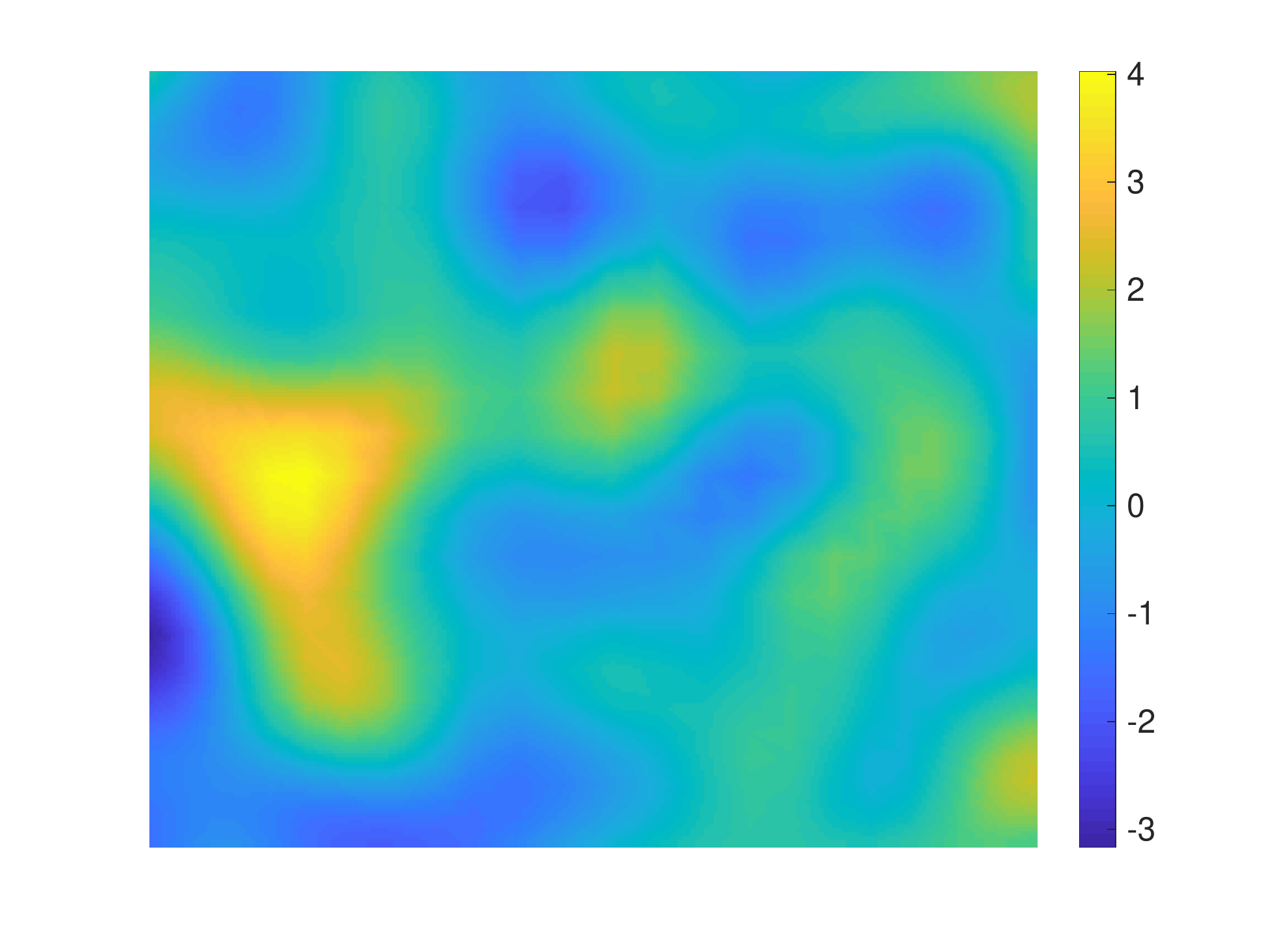}
   \put (40,-3) {\scriptsize $p(x)$}
  \end{overpic}
    \begin{overpic}[width=0.45\textwidth,trim=35 10 45 5, clip=true,tics=10]{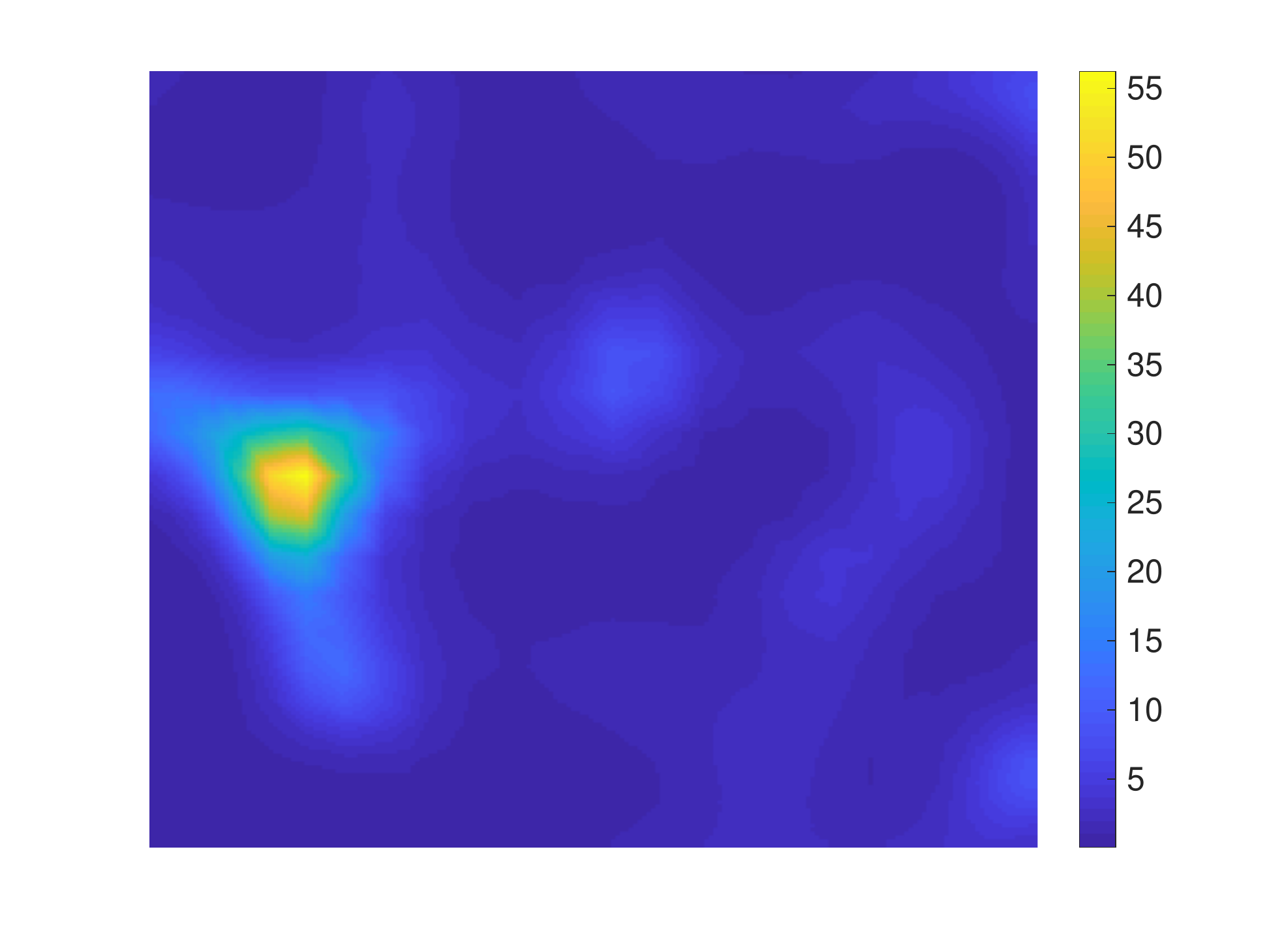}
     \put (43,-3) {\scriptsize $\kappa(x)$}
  \end{overpic}
 \end{center}
\caption{Example 3: The true solution $p(x)$ and $\kappa(x)$.}\label{exact_eg3}
  \end{figure}

In the last example, we consider the permeabilities as a random field.  Especially, the log-diffusivity field $\log\kappa(x):=p(x)$ is endowed with a Gaussian process prior, with mean zero and an isotropic  kernel:
\begin{equation*}
C(x_1,x_2)=\sigma^2 \exp\Big(-\frac{\|x_1-x_2\|^2}{2l^2}\Big),
\end{equation*}
for which we choose variance $\sigma^2=1$ and a length scale $l = 0.1$. This prior allows the field to be easily parameterized with a Karhunen-Loeve expansion:
\begin{equation}
p(x; z) \approx \sum^{n}_{i=1} z^i \sqrt{\lambda_i} \phi_i(x),
\end{equation}
where $\lambda_i$ and $\phi_i(x)$ are  the eigenvalues and eigenfunctions, respectively, of the integral operator on $[0,1]^2$ defined by the kernel $C$, and the parameter $z=(z^1,\cdots, z^n)$ are endowed with independent standard normal priors, $z^i \sim N(0,1)$. These parameters then become the targets of inference. In particular, we truncate the Karhunen-Loeve expansion at $n=111$ modes.  The true solution $p(x)$ and $\kappa(x)$  used to generate the test data are shown in Fig.\ref{exact_eg3}.  The measurement sensors of $u$ are evenly distributed over $\Omega$ with grid spacing 0.1, \textcolor{black}{i.e., $d \in \R^{81}$.  The observational errors are taken to be additive and Gaussian:
\begin{equation*}
d_j = u(x_j;z) +\xi_j, j=1,\cdots,81,
\end{equation*}
with $\xi_j \sim N(0,0.05^2)$. } In this example, three hidden layers and 150 neurons per layer are used in $\mathcal{NN}^L$, while 1 hidden with 150 neurons are used in $\mathcal{NN}^H$.  The regularization rate is set to $\lambda = 10^{-6}$.  If the refinement is set to occur, we choose $Q=50$ random points  to train the multi-fidelity DNN  $\mathcal{NN}^H$.

 \begin{figure}
\begin{center}
\begin{overpic}[height=6.4cm,width=4.4cm, trim= 35 10 45 5, clip=true,tics=10]{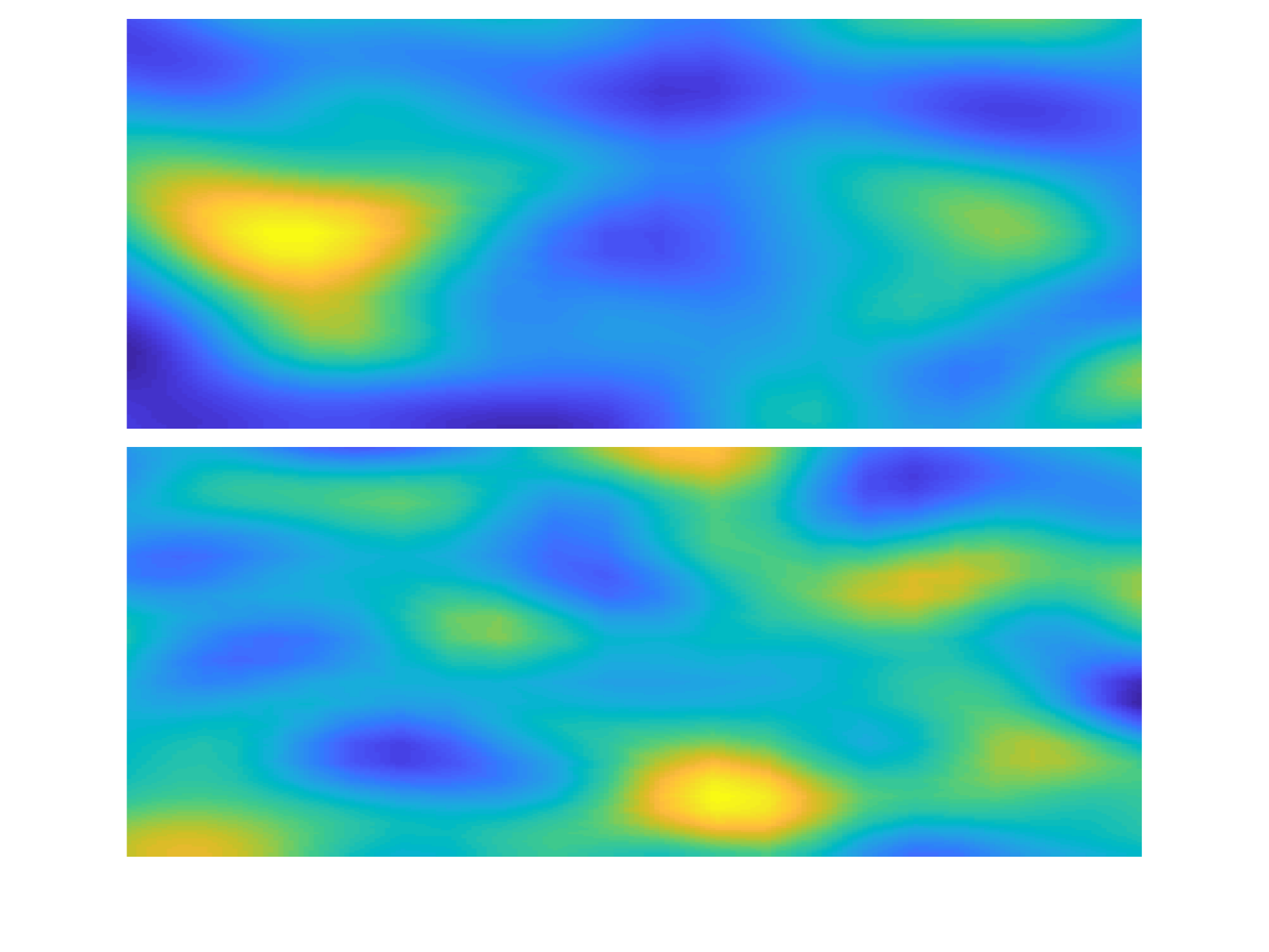}
         \put (43,-3) {\scriptsize Direct}
  \end{overpic}
    \begin{overpic}[height=6.4cm,width=4.4cm,trim= 35 10 45 5, clip=true,tics=10]{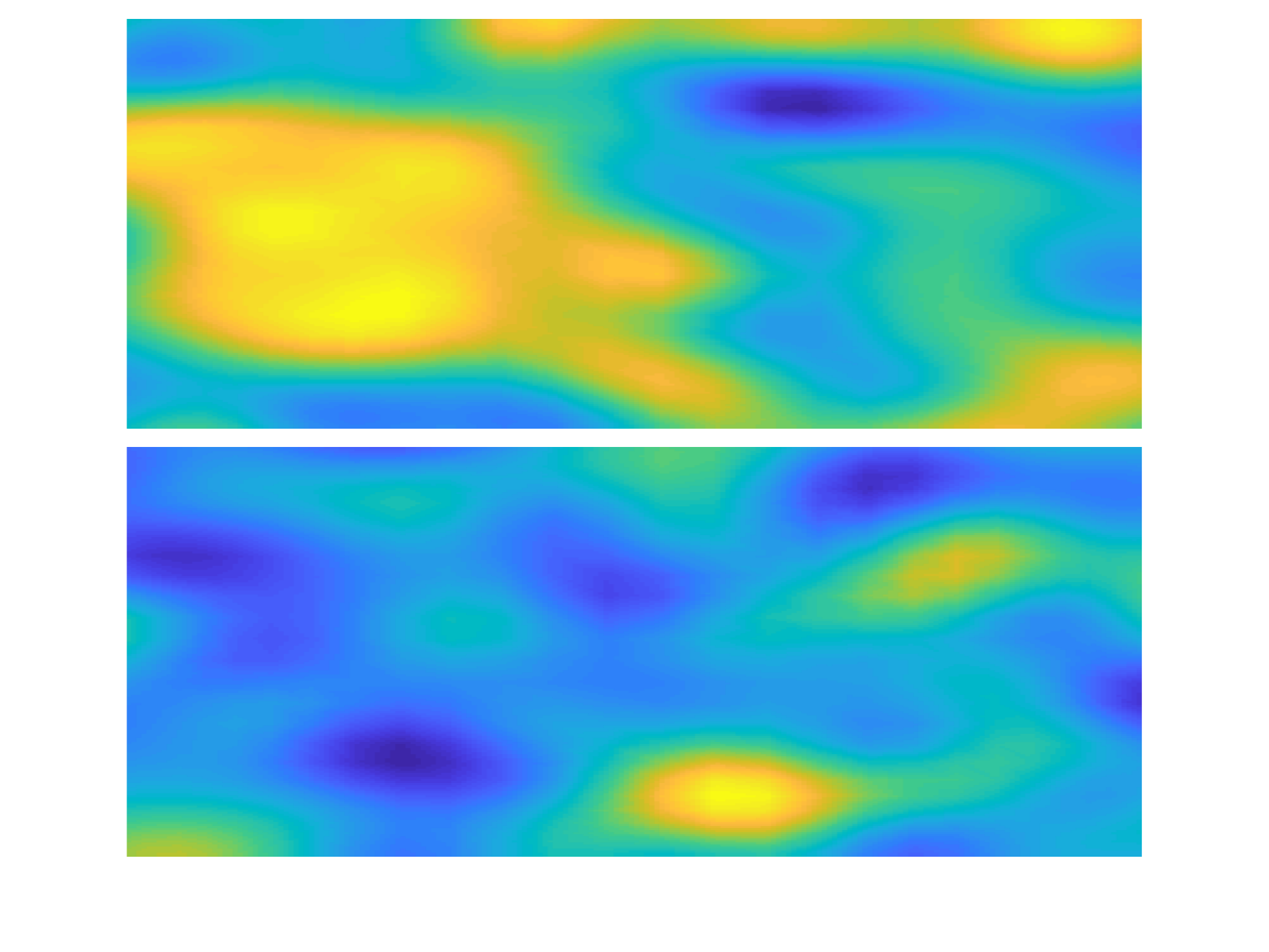}
         \put (43,-3) {\scriptsize DNN}
  \end{overpic}
  \begin{overpic}[height=6.4cm,width=4.4cm,trim=35 10 45 5, clip=true,tics=10]{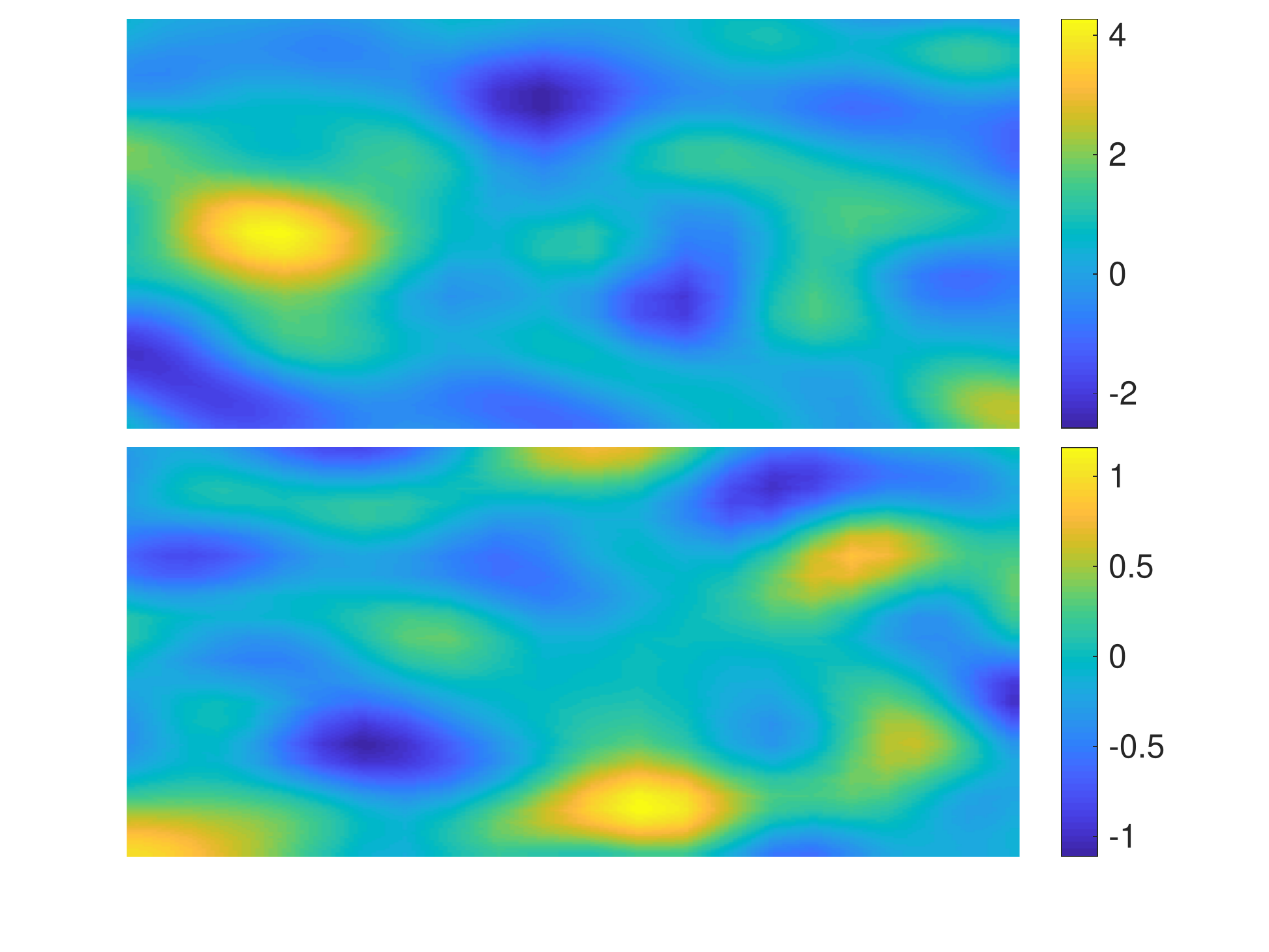}
       \put (40,-3) {\scriptsize ADNN}
  \end{overpic}
  \end{center}
\caption{Example 3:  Posterior mean (top) and  posterior standard deviation (bottom) of $p(x)$ arising from direct MCMC, prior-based DNN approach ($N=100$) and  ADNN ($tol=0.1$), respectively. }\label{mean-eg3}
  \end{figure}

   \begin{figure}
\begin{center}
       \begin{overpic}[height=6.4cm,width=4.4cm,trim= 35 10 45 5, clip=true,tics=10]{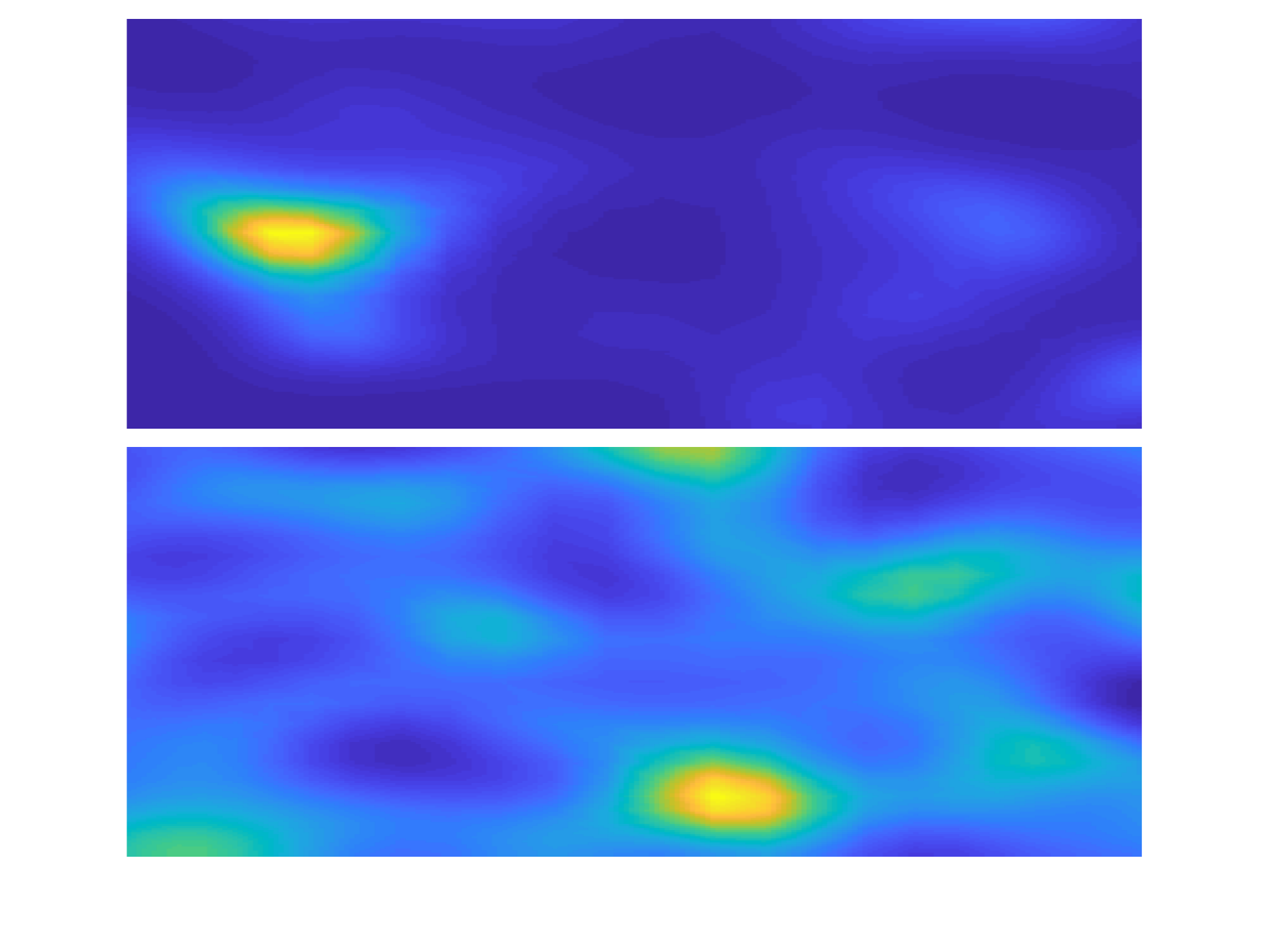}
          \put (43,-3) {\scriptsize Direct}
  \end{overpic}
    \begin{overpic}[height=6.4cm,width=4.4cm,trim= 35 10 45 5, clip=true,tics=10]{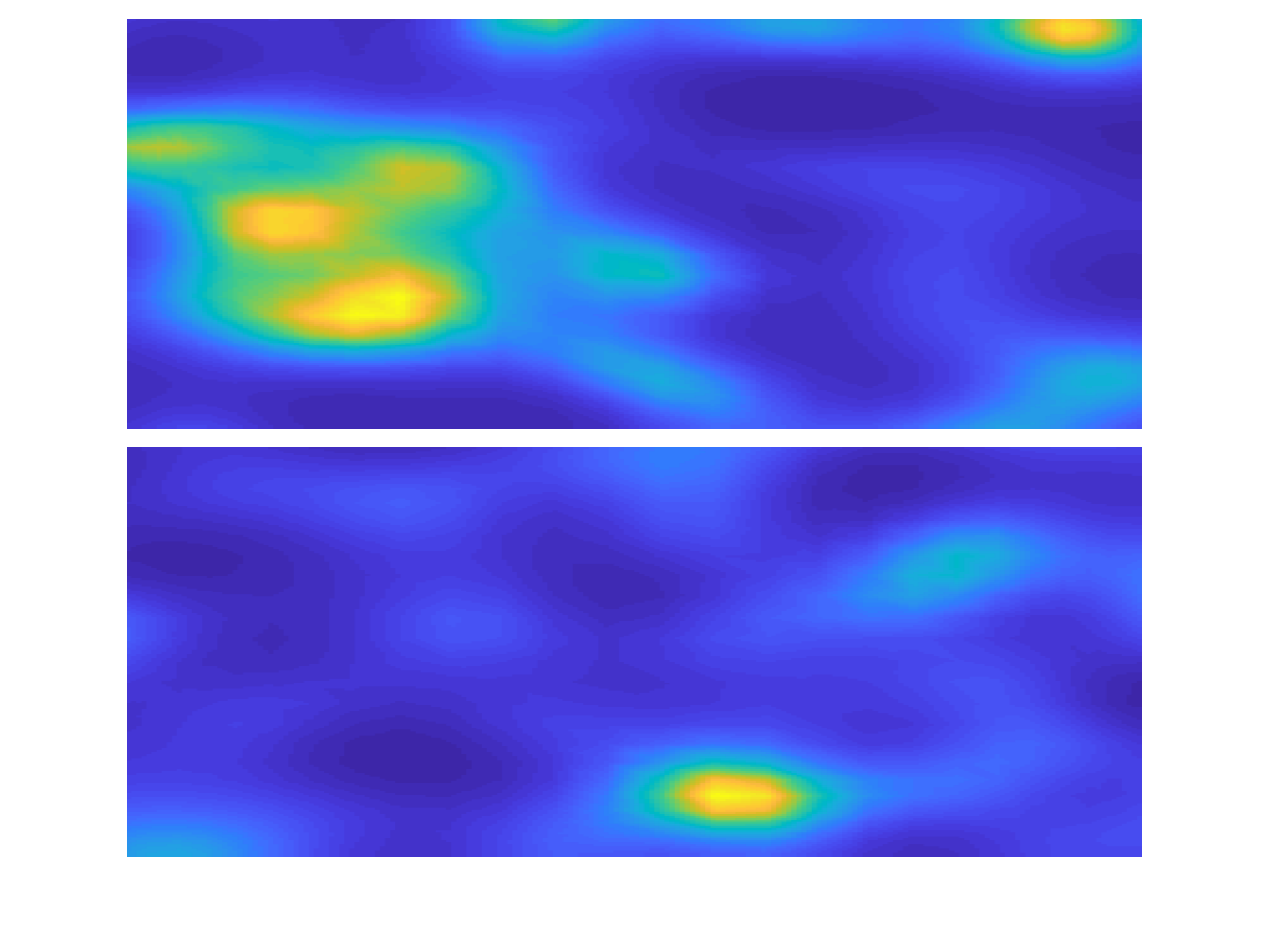}
         \put (43,-3) {\scriptsize DNN}
  \end{overpic}
  \begin{overpic}[height=6.4cm,width=4.4cm,trim=35 10 45 5, clip=true,tics=10]{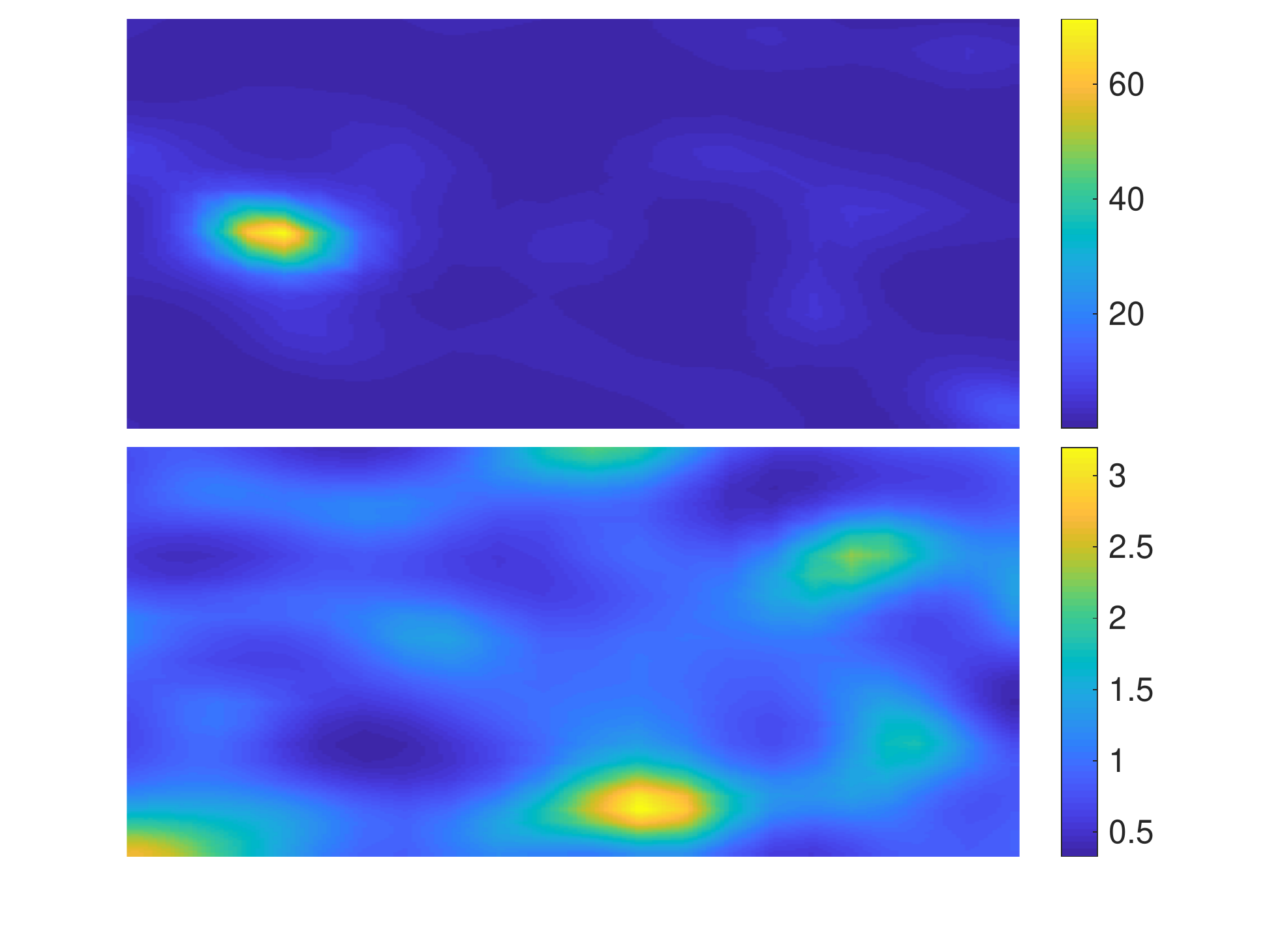}
       \put (40,-3) {\scriptsize ADNN}
  \end{overpic}
\end{center}
\caption{Example 3:  Posterior mean (top) and  posterior standard deviation (bottom) of $\kappa(x)$ arising from direct MCMC, prior-based DNN approach ($N=100$) and  ADNN ($tol=0.1$), respectively.  }\label{mean-eg3-2}
  \end{figure}

Figs. \ref{mean-eg3} and \ref{mean-eg3-2} plot the conditional mean arising from three different approaches.  As expected, a poor estimate is obtained by the prior-based DNN approach. The results are improved  with the ADNN algorithm.  The computational costs   for the different algorithms are shown in Table \ref{eg3_time}.  Building a DNN surrogate using $N=100$ training points  requires an offline CPU time of 26.8s, whereas its online evaluation requires 9.5s.    On the other hand, for the ADNN algorithm with  $tol=0.1$, the  offline and online CPU times are 26.8s and 88.9s, respectively.  This demonstrated that the ADNN can provide with much more accurate results, yet with less computational time.

\begin{table}[tp]
      \caption{Example 3. Computational times, in seconds, given by three different methods. }\label{eg3_time}
  \centering
  \fontsize{6}{12}\selectfont
  \begin{threeparttable}
    \begin{tabular}{ c ccccc}
  \toprule
   & \multicolumn{2}{c}{Offline}&\multicolumn{2}{c}{Online}\cr
\cmidrule(lr){2-3} \cmidrule(lr){4-5}
  \multirow{1}{*}{Method}  &$\text{$\#$ of high-fidelity evaluations}$&CPU(s) &$\text{$\#$ of high-fidelity evaluations}$&CPU(s)     &\multirow{1}{*}{Total time(s)}\cr
  \midrule
    Direct                           & $-$       & $-$         & 5$\times 10^4$    &5507.5      & 5507.5   \cr
   DNN                           & 100     & 26.8        & $-$                             & 9.5                         &36.3    \cr
   ADNN                         & 100    & 26.8         & 950                      & 88.9                       & 115.7  \cr
    \bottomrule
      \end{tabular}
    \end{threeparttable}
\end{table}

\section{Summary} \label{sec:summary}

We have presented an adaptive DNN-based surrogate modelling procedure for Bayesian inference problems.  In our computational procedure, we first construct offline a DNNs-based surrogate according to the prior distribution, and then, this prior-based DNN-surrogate will be adaptively \& locally refined online using only a few high-fidelity simulations. In particular, in the refine procedure, we construct a new shallow neural network that view the previous constructed surrogate as an input variable -- yielding a composite multi-fidelity neural network approach. The performance of the proposed strategy has been illustrated by three numerical examples. There are several potential extensions of the present scheme, which are currently under investigation.

Finally, we remark that although our approach was presented in the MCMC framework, the idea can be easily extended to other approaches such as the Sequential Monte Carlo approach \cite{Beskos2015sequential}  or optimization-based sampling approaches \cite{Bardsley2014SISC,Wang+Bardsley2017SISC, Wang+Cui2019scalable}.  The extension of the present algorithm to Ensemble Kalman inversion \cite{Iglesias+Law+Kody2013ensemble,Yan+Zhou19IJUQ} is also straightforward.


\end{document}